\documentclass[bpam]{ipart}

\pubyear{2024}
\volume{0}
\issue{0}
\firstpage{1}
\lastpage{1}

\usepackage{booktabs} 

\usepackage[ruled]{algorithm2e} 

\SetAlFnt{\small}
\SetAlCapFnt{\small}
\SetAlCapNameFnt{\small}
\SetAlCapHSkip{0pt}

\usepackage{wrapfig}
\usepackage{multirow}

\usepackage{amssymb,amsmath,amsthm}
\usepackage{todonotes}
\setuptodonotes{inline}

\theoremstyle{plain}
\newtheorem{thm}{Theorem}[section]

\theoremstyle{definition}

\theoremstyle{definition}
\newtheorem{rem}[thm]{Remark}

\newtheorem*{ass*}{Assumption}

\numberwithin{equation}{section}

\newcommand{\RR}{\mathbb{R}} 
\DeclareMathOperator{\im}{Im}

\DeclareMathOperator{\vf}{VF}
\DeclareMathOperator{\fk}{FK}
\DeclareMathOperator{\hk}{HK}
\DeclareMathOperator{\cg}{CG}
\DeclareMathOperator{\hg}{HG}
\DeclareMathOperator{\GG}{GG}

\renewcommand{\l}{$L^2$}

\begin{document}

\begin{frontmatter}

\title[Hodge Decomposition]{Topology-preserving Hodge Decomposition in the Eulerian Representation\protect\thanksref{T1}}
\thankstext{T1}{This paper is an extended version of a SIGGRAPH Asia conference paper.}

\begin{aug}

\author{\fnms{Zhe} \snm{Su}\ead[label=e1]{suzhe@msu.edu}}
\address{Department of Mathematics, Michigan State University, East Lansing, MI 48824, USA\\
\printead{e1}}
\author{\fnms{Yiying} \snm{Tong}
\thanksref{t2}
\ead[label=e2]{ytong@msu.edu}}
\thankstext{t2}{Corresponding author.}
\address{Department of Computer Science and Engineering, Michigan State University, East Lansing, MI 48824, USA\\
\printead{e2}}
\and
\author{\fnms{Guo-Wei} \snm{Wei}
\thanksref{t3}      
\ead[label=e3]{weig@msu.edu}}%
\thankstext{t3}{Corresponding author.} 
\address{Department of Mathematics, Michigan State University, East Lansing, MI 48824, USA\\
Department of Biochemistry and Molecular Biology, Michigan State University, East Lansing, MI 48824, USA\\
Department of Electrical and Computer Engineering, Michigan State University, East Lansing, MI 48824, USA\\
\printead{e3}}
\end{aug}
\received{\sday{1} \smonth{9} \syear{2024}}

\begin{abstract}
The Hodge decomposition is a fundamental result in differential geometry and algebraic topology, particularly in the study of differential forms on a Riemannian manifold. Despite extensive research in the past few decades,
topology-preserving Hodge decomposition of scalar and vector fields on manifolds with boundaries in the Eulerian representation remains a challenge due to the implicit incorporation of appropriate topology-preserving boundary conditions.      
In this work, we introduce a comprehensive 5-component topology-preserving Hodge decomposition that unifies normal and tangential components in the Cartesian representation. 
Implicit representations of planar and volumetric regions defined by level-set functions have been developed. 
Numerical experiments on various objects, including single-cell RNA velocity,  validate the effectiveness of our approach,  confirming the expected rigorous \l-orthogonality and the accurate cohomology.
\end{abstract}


\begin{keyword}
\kwd{vector field decomposition}
\kwd{cohomology}
\kwd{topological Laplacian}
\kwd{Cartesian grids}
\end{keyword}


\end{frontmatter}



%
%

%
%


%
%


\section{Introduction}

As a hallmark of the 20th century's mathematics, de Rham-Hodge theory is a mathematical framework that combines aspects of differential geometry, algebraic topology, and partial differential equation (PDE). It provides a powerful set of tools for studying the properties of differential forms on smooth manifolds and is particularly important in the study of Riemannian geometry, Hodge theory, and algebraic topology. Some of its key concepts include de Rham cohomology, a tool for studying the topology of smooth manifolds using differential forms, and  Hodge theory, which provides a way to study the de Rham cohomology groups using the tools of Riemannian geometry. Mathematically, the Hodge decomposition theorem states that on an oriented compact Riemannian manifold, any differential form can be uniquely decomposed into an exact form, a coexact form, and a harmonic form. This is a crucial result that bridges the gap between topology (de Rham cohomology) and analysis (harmonic forms).

In science and engineering, the Hodge decomposition of vector fields into components with specific topology-preserving properties has a broad range of applications. For instance, it plays a crucial role in computational fluid dynamics~\cite{yang2021clebsch}, flow processing and visualization \cite{sawhney2020monte}, geometric modeling~\cite{wang2021computing}, spectral data analysis~\cite{keros2023spectral}, and machine learning~\cite{su2024hodge}. The general form of the decomposition based on Hodge's foundational work~\cite{hodge1989theory} applies to differential forms (covariant antisymmetric tensor fields). It decomposes the space of differential $k$-forms into a direct sum of \l-orthogonal subspaces. As an essential special case, vector fields on a compact domain in the 2D or 3D Euclidean space can be orthogonally decomposed into divergence-free and curl-free components, as in the Helmholtz-Hodge decomposition.

\begin{figure}
	\centering
	\includegraphics[trim=6cm 7cm 6cm 7cm, clip, width=0.12\textwidth]{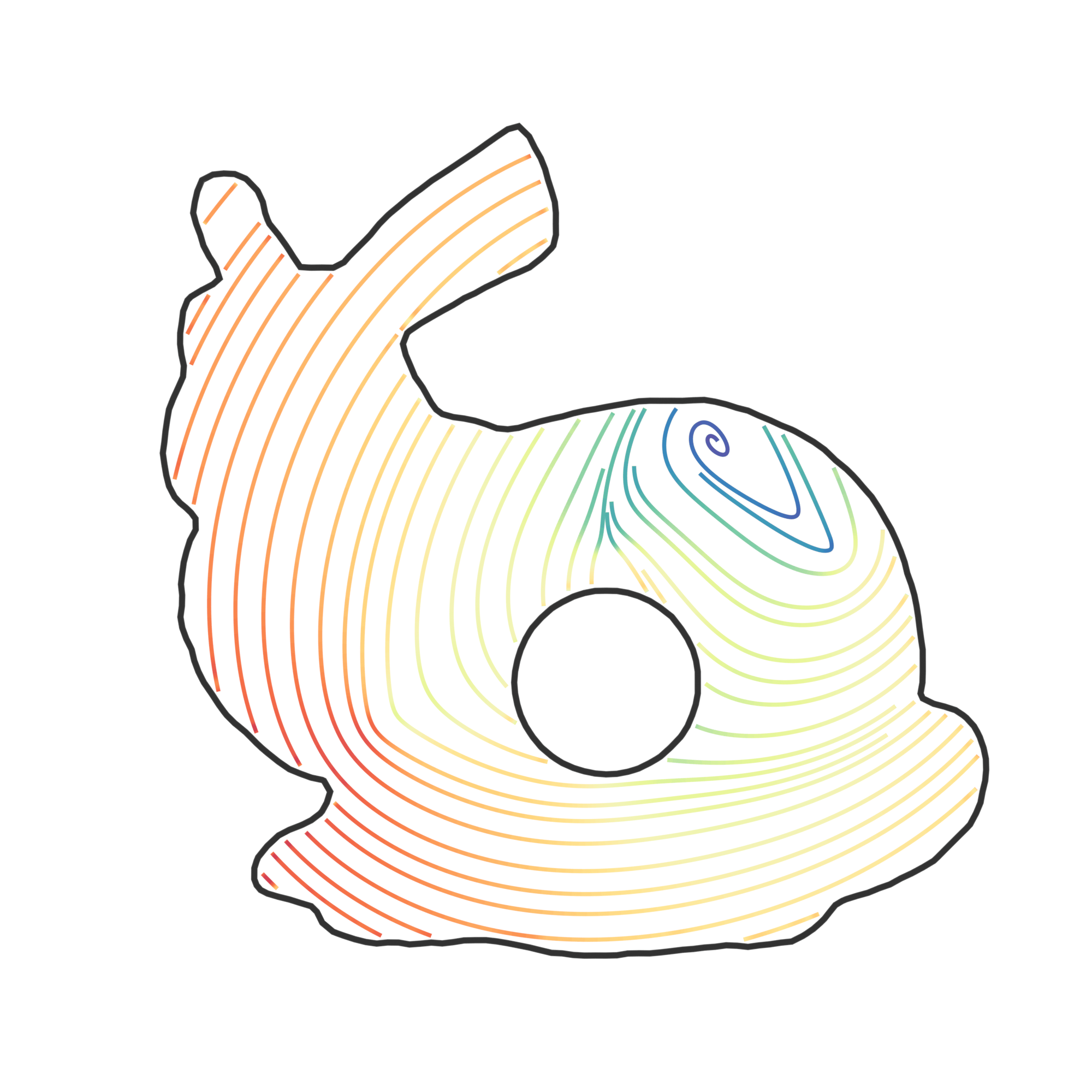}
 \raisebox{1cm}{\Large$=$}
	\includegraphics[trim=6cm 7cm 6cm 7cm, clip,width=0.12\textwidth]{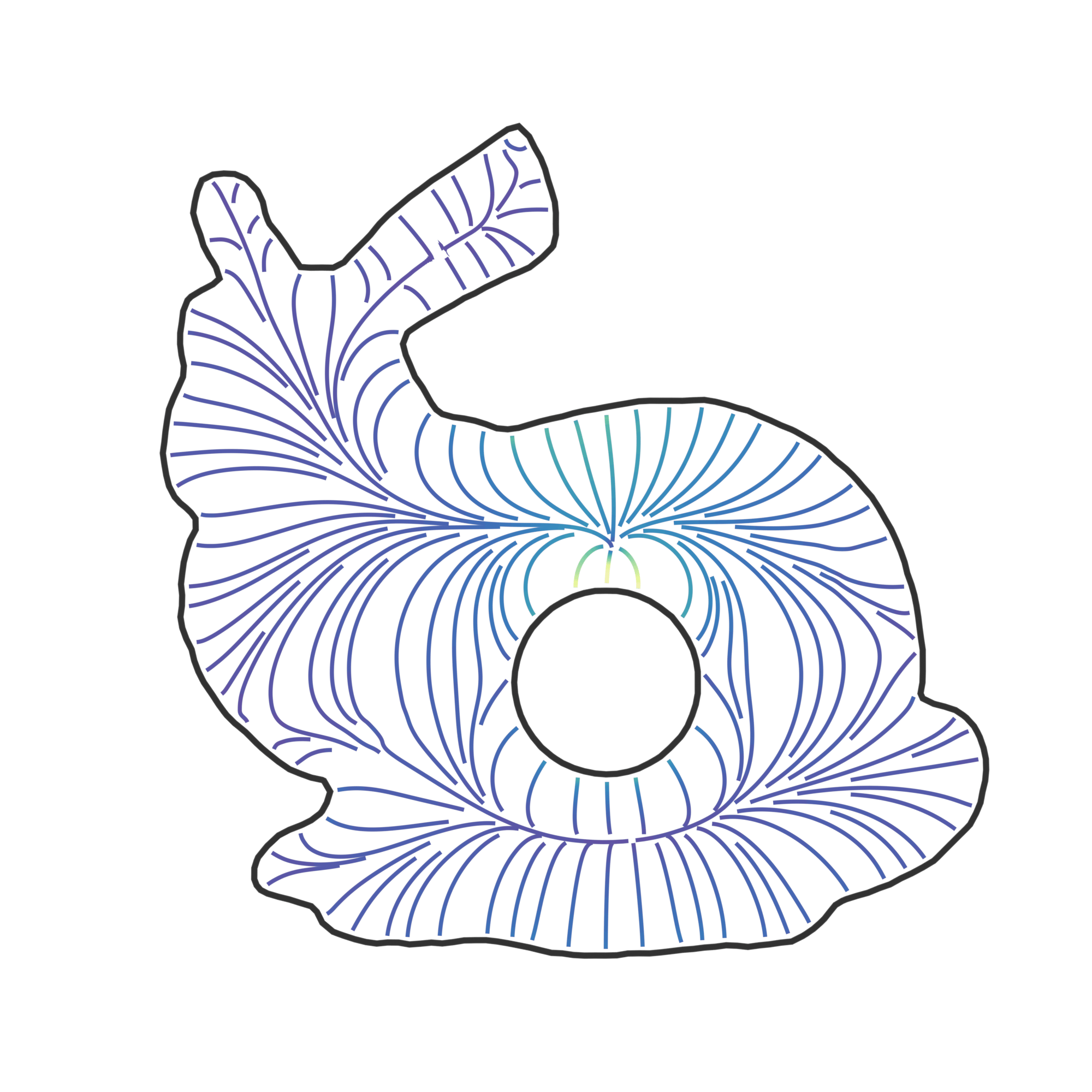}
  \raisebox{1cm}{\Large$+$}
	\includegraphics[trim=6cm 7cm 6cm 7cm, clip,width=0.12\textwidth]{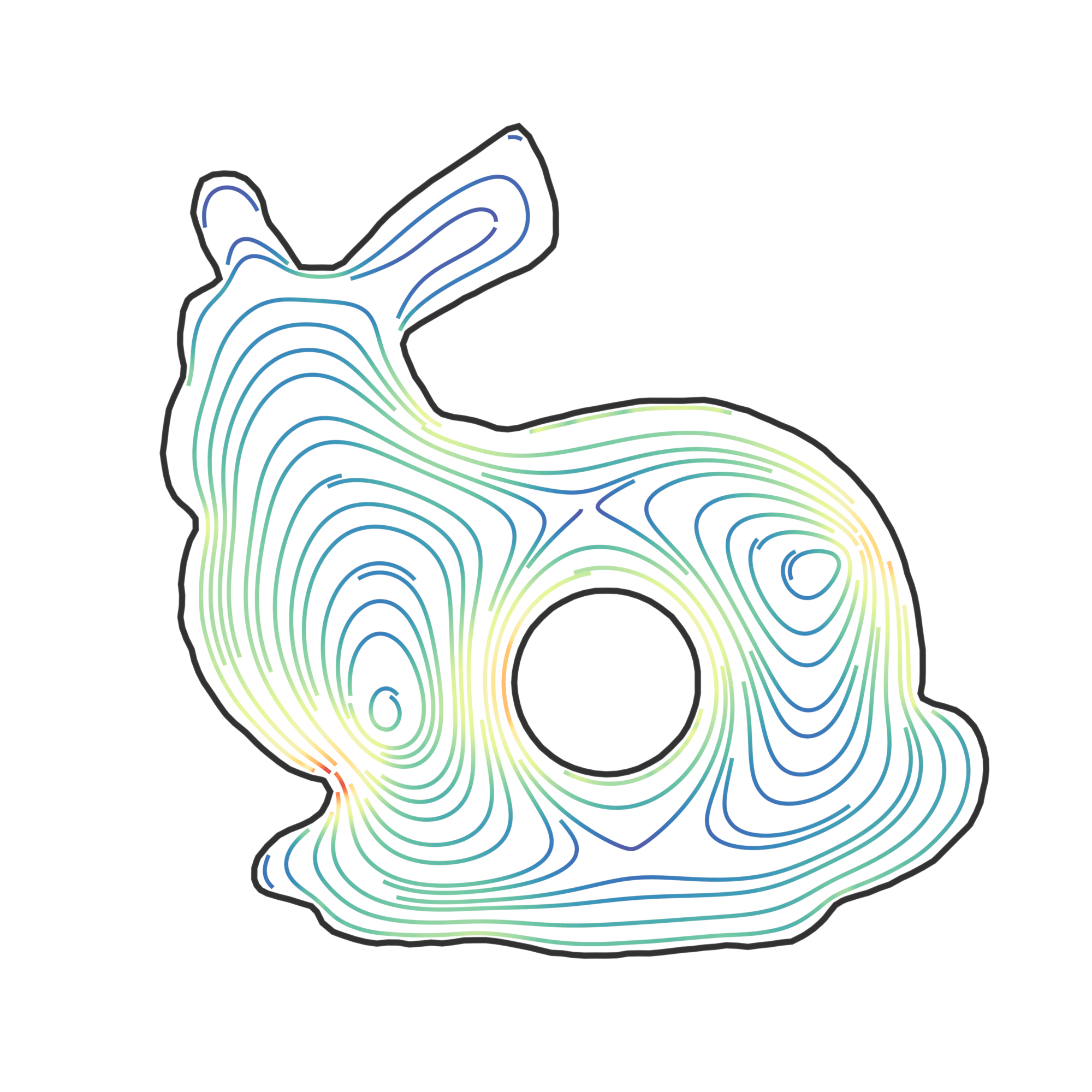}
  \raisebox{1cm}{\Large$+$}
	\includegraphics[trim=6cm 7cm 6cm 7cm, clip,width=0.12\textwidth]{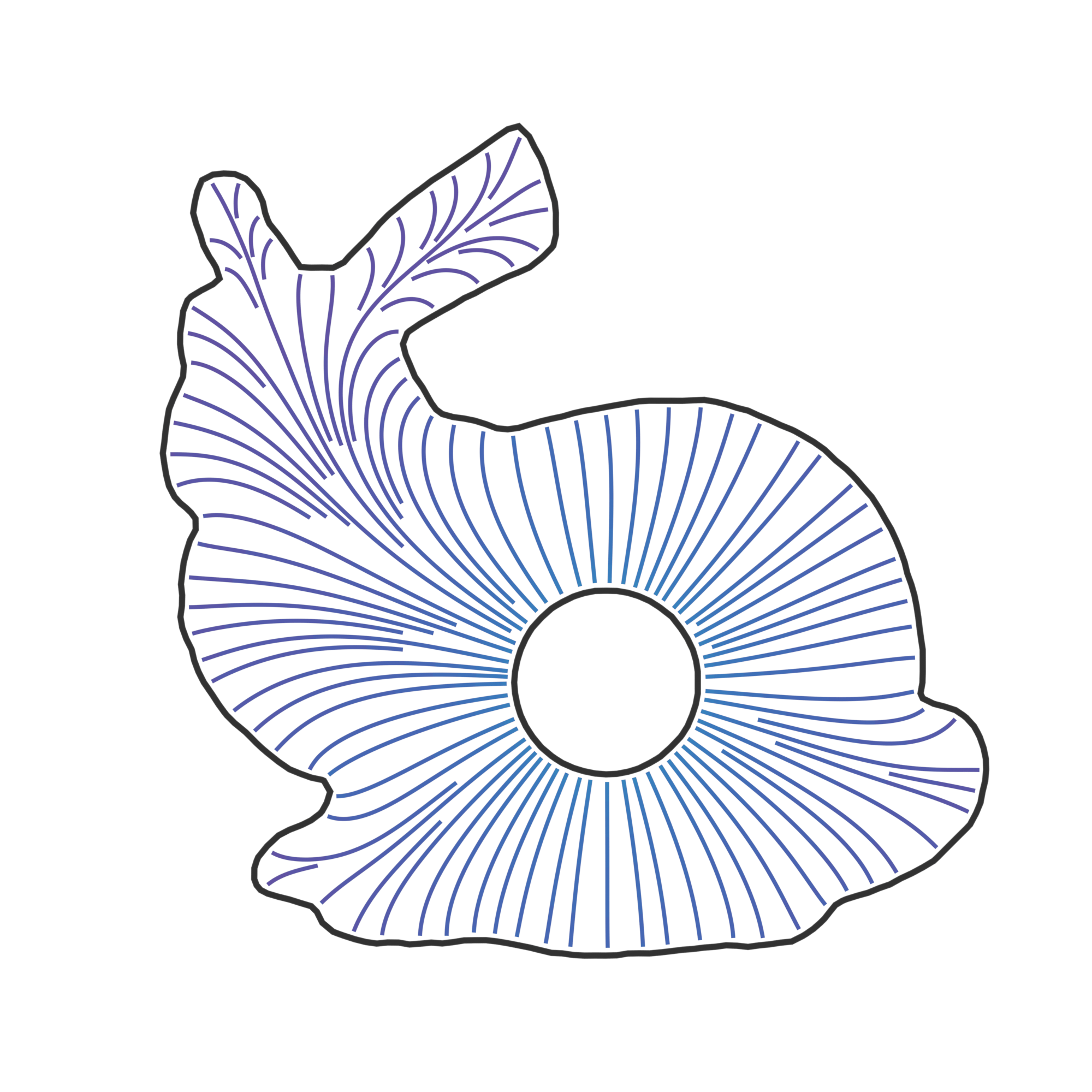}
  \raisebox{1cm}{\Large$+$}
	\includegraphics[trim=6cm 7cm 6cm 7cm, clip,width=0.12\textwidth]{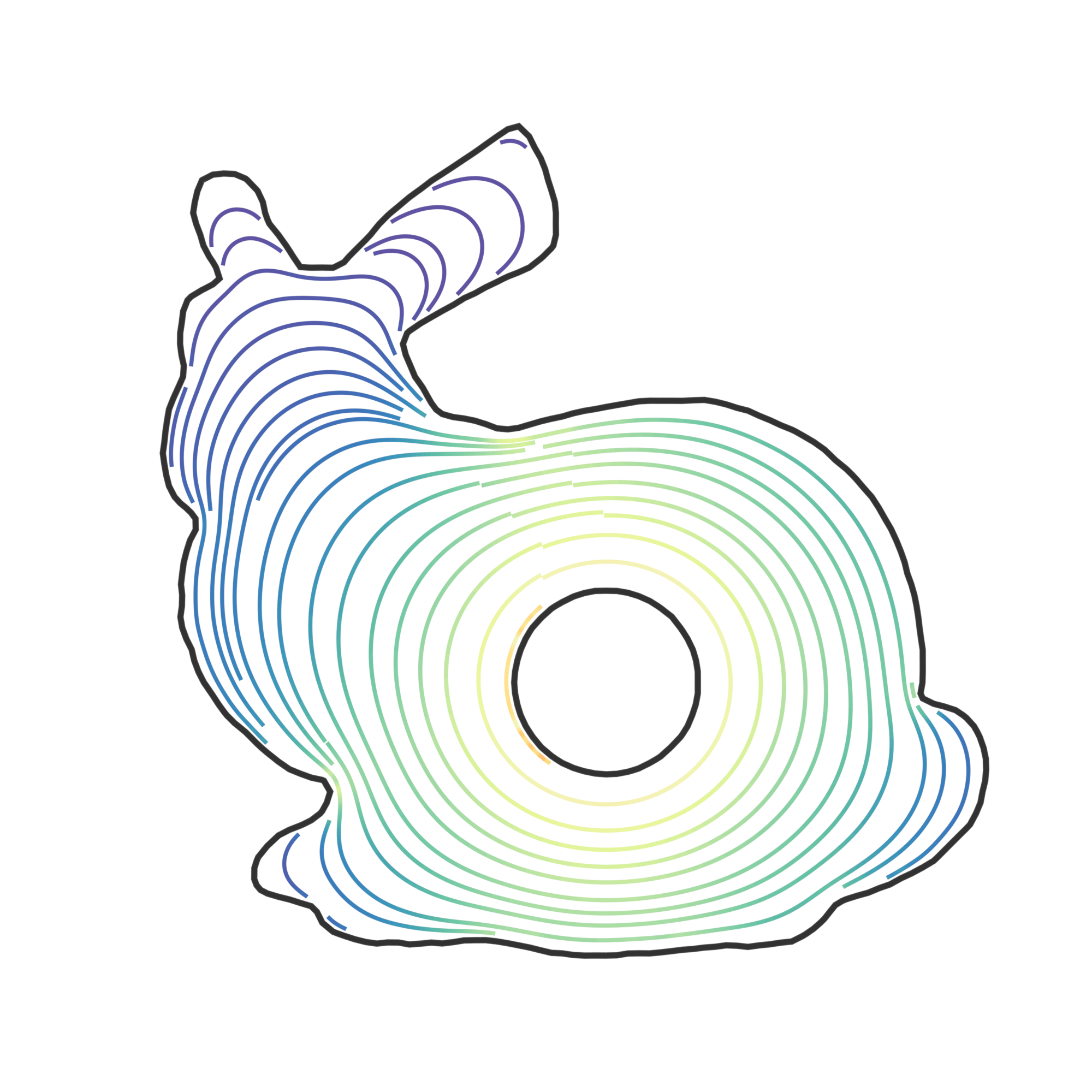}
  \raisebox{1cm}{\Large$+$}
	\includegraphics[trim=6cm 7cm 6cm 7cm, clip,width=0.12\textwidth]{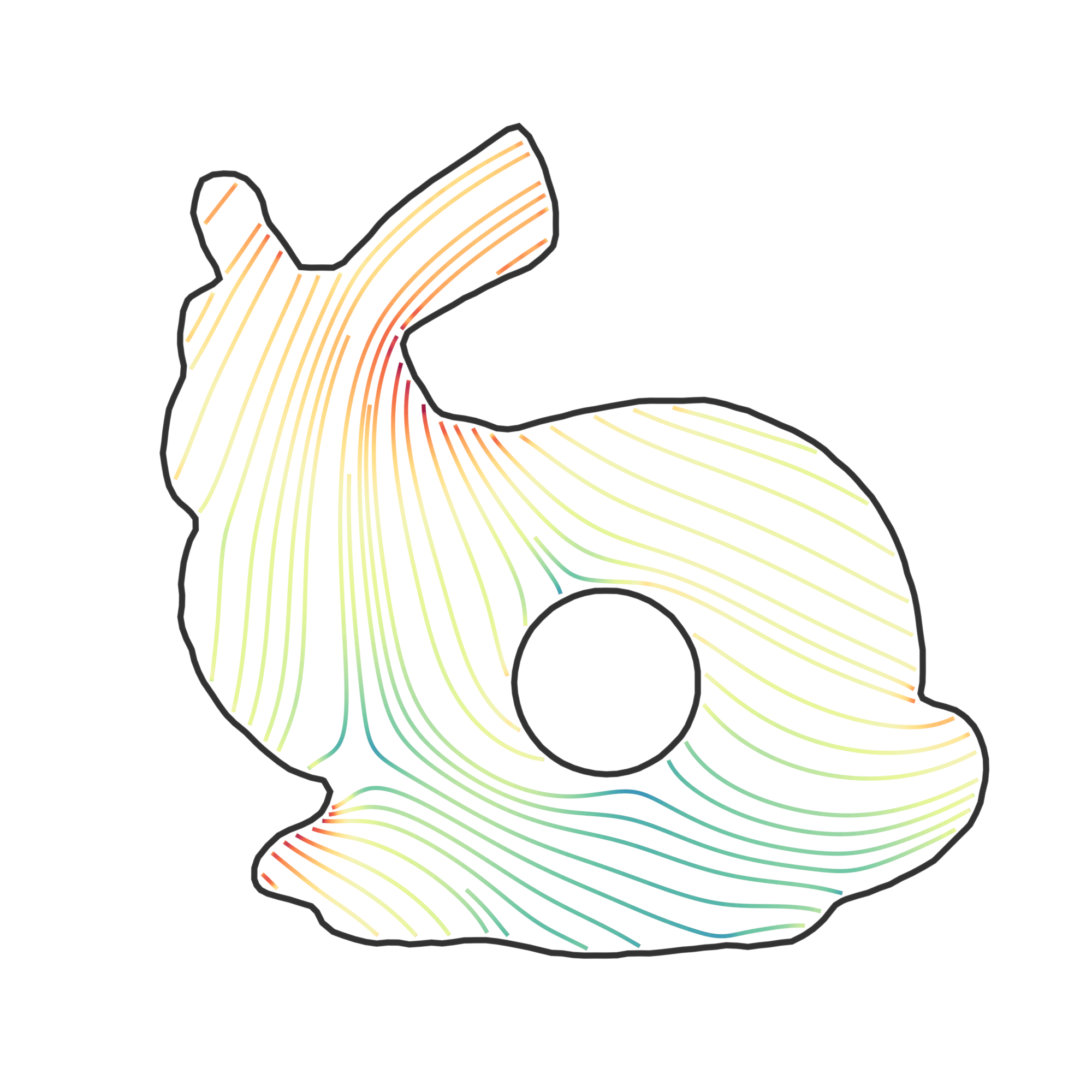}
	\caption{5-component Hodge decomposition. From left to right: the original vector field, the normal gradient field, the tangential curl field, the normal harmonic field, the tangential harmonic field, and the curly gradient field.}
	\label{fig.hd.2d.bunny.stepwise}
\end{figure}

On closed manifolds, the classical Hodge decomposition \cite{hodge1989theory} involves Laplacians, which are second-order linear differential operators with finite-dimensional kernels called harmonic spaces. The dimensions of these kernels correspond to the counts of topological features of the underlying manifolds, leading to implementations through linear systems with \emph{low} rank deficiency. However, on manifolds with boundary, the situation becomes subtle. The spaces of harmonic fields are infinite-dimensional, resulting in linear systems with substantial rank deficiency. To mitigate this issue, specific boundary conditions, such as the normal (Dirichlet) and tangential (Neumann) boundary conditions, have been introduced~\cite{shonkwiler2009poincare}. Under these conditions, kernels are again finite-dimensional, with correspondences to the topology of the manifold and its boundary.

The domains for vector field processing through finite element type approaches in geometric modeling and computer graphics\cite{poelke2016boundary, tong2003discrete} are frequently represented by simplicial meshes, e.g., the unstructured surface or volume meshes resulting from either Lagrangian or Eulerian formulations in fluid simulations. On the other hand, structured meshes such as Cartesian grids are advantageous in other types of discretization, including finite-volume, finite-difference, and pseudo-spectral methods~\cite{fedkiw2001visual, stam2023stable, liu2015model}, as well as for
\begin{wrapfigure}[10]{r}{0.29\columnwidth}
\vspace*{-5mm}\hspace*{-1mm}\centering
\includegraphics[trim=4cm 1cm 4cm 0.5cm, clip, width=0.31\columnwidth]{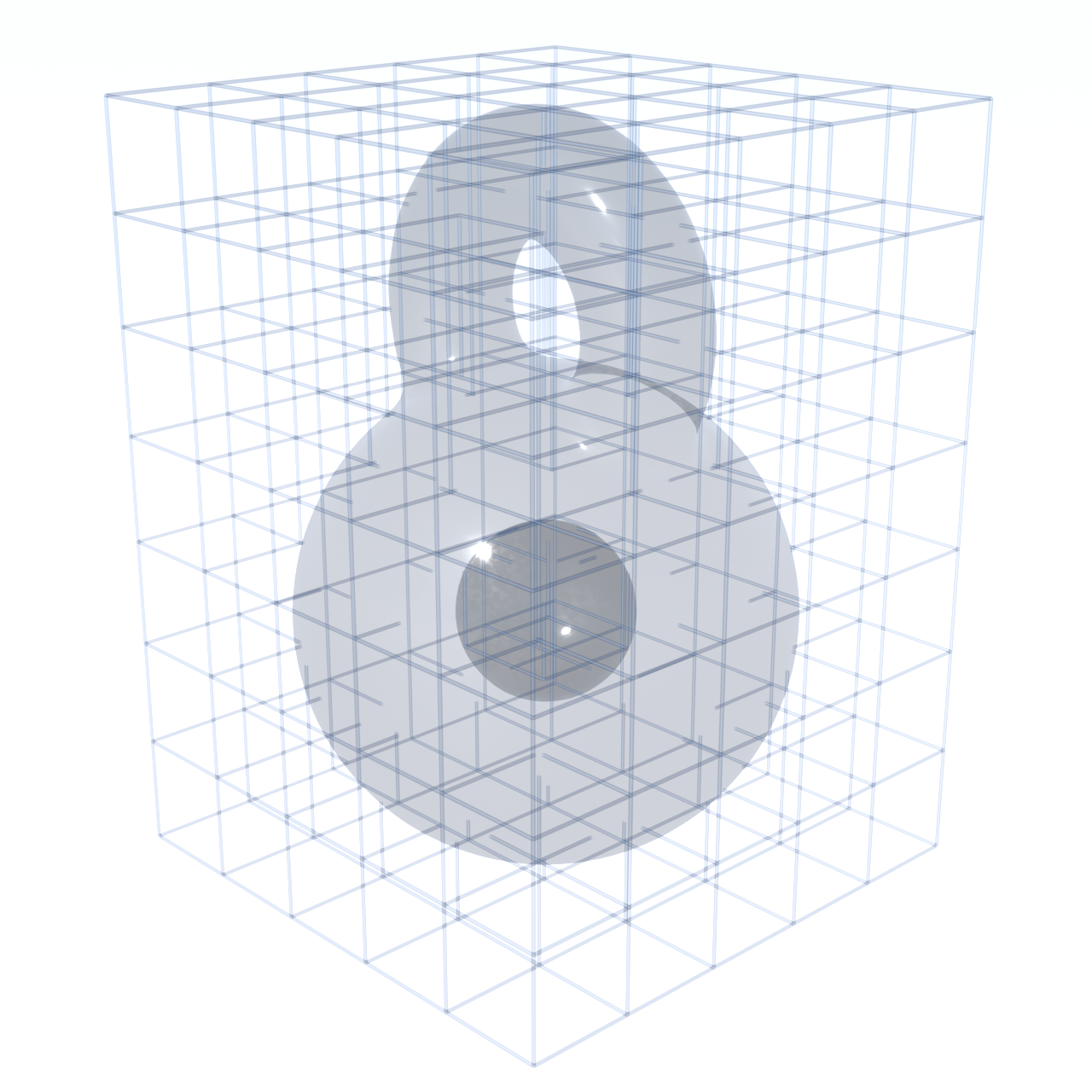}%
\end{wrapfigure}
convolutional neural networks on voxelized data.
However, for the complete
$5$-component Hodge decomposition of vector fields on compact domains, the existing methods \cite{poelke2017hodge,zhao20193d,razafindrazaka2019consistent}
rely exclusively on simplicial mesh discretization frameworks,
where the underlying manifolds are modeled either by
triangular or tetrahedral meshes, or deformed grids. To facilitate the direct application of the $5$-component Hodge decomposition on a compact domain embedded within a regular Cartesian grid, we propose linear systems directly assembled based on a given level-set function, without resorting to tessellation and data conversion.

\subsection{Related work}

In a board sense, the development of Hodge Laplacians is part of the ongoing effort on persistent topological  Laplacians (PTLs) that include persistent topological Laplacians defined on point cloud data \cite{wang2020persistent,liu2023algebraic,grbic2022aspects} and on manifolds \cite{chen2021evolutionary,su2024persistent}. 
PTLs overcome some of the limitations of persistent homology, a vital tool in topological data analysis (TDA) \cite{zomorodian2004computing,edelsbrunner2008persistent,liu2022biomolecular}. These approaches, formulated in terms of topological deep learning \cite{cang2017topologynet},  have had tremendous success in analyzing biomolecular data\cite{nguyen2019mathematical,nguyen2020review,nguyen2020mathdl} and forecasting emerging viral variants \cite{chen2022persistent}.  

The Hodge decomposition for smooth manifolds with boundary has been thoroughly discussed in \cite{schwarz2006hodge} with up to $4$ components. A refinement of the decomposition to five terms is presented in \cite{shonkwiler2013poincare}, which offers insights on further splitting the cohomology groups (which correspond to the kernels of Laplacians) into portions derived from the interior and boundary of the manifold. An elementary exposition of this $5$-component Hodge decomposition for compact domains in $\RR^3$ in terms of vector and scalar fields can be found in \cite{cantarella2002vector}.

In the discrete case, the Helmholtz-Hodge decomposition (HHD) has been implemented in various methods. In \cite{polthier2000variational}, a global variational approach was used for HHD of piecewise constant vector fields on triangulated surfaces, by $L^2$-projection to the spaces of curl-free and divergence-free components. It was then extended to the 3D case on tetrahedral meshes \cite{tong2003discrete, polthier2003identifying}, and to the 2D case on regular grids \cite{guo2005efficient}. A meshless algorithm of HHD was introduced in \cite{petronetto2009meshless} for point vectors in $\RR^2$. See \cite{bhatia2012helmholtz,de2016vector} for a comprehensive discussion on the theory and practice of HHD for 2D and 3D vector fields.

The complete $5$-component discrete Hodge decomposition was subsequently introduced for piecewise constant vector fields, for surface triangle meshes \cite{poelke2016boundary}, and for tetrahedral meshes \cite{poelke2017hodge}. The framework aligns well with the smooth case, accurately capturing the topological structures. The decomposition can be done more efficiently when using Whitney bases representations for differential forms, such as Nedelec elements on edges, and Raviart-Thomas elements on faces~\cite{razafindrazaka2019consistent,zhao20193d}. The generalization and unification of these finite-element methods in the context of exterior calculus were proposed in FEEC (finite-element exterior calculus)~\cite{arnold2018finite}, which studied the Hilbert complexes formed by existing elements in the preservation of cohomology, designed novel serendipity elements, and provided rigorous error analysis. 

The gradient, curl, and divergence operators, which are essential to the decomposition, can still be defined even on abstract simplicial complexes through boundary operators on such complexes~\cite{lim2020hodge}. Moreover, the nilpotent property of boundary operators still induces a cohomology, which enables a decomposition. Weights associated with simplices can be used to define $L^2$-norms similar to the treatment in DEC. However, the notion of correspondences to the topology of a domain including its boundary is lost unless the simplicial complex comes from a simplicial mesh, which then reduces the general complex back to DEC. Nevertheless, implementations based on abstract complexes such as PyDEC~\cite{bell2012pydec} produce the valid discretization for both DEC and (linear-order) FEEC.

All of these works focused on either simplicial meshes, or cubicle meshes that are deformed to conform to the boundaries, whereas we propose the decomposition on domains enclosed by level-sets in the Cartesian grids.

\subsection{Contributions}

We aim to fill the gap for the topology-preserving $5$-component Hodge decomposition on Cartesian grids with 2D/3D domains defined by level set functions. Unlike existing methods, we do not require an explicit tessellation of the domain into simplicial meshes, eliminating the reliance on high-quality meshing tools and streamlining the decomposition. Moreover, we maintain the pairwise \l-orthogonality among the components while preserving the consistency with the topology.

In the journal extension, we further performed extensive empirical error analysis on multiple new examples for the spectra of discrete Laplacians, demonstrating both the correct multiplicity of zeros determined by the topology and the convergence of the leading nonzero eigenvalues. The accuracy of the decomposition was also tested on both 2D and 3D examples for vector fields with analytic expressions, providing numerical evidence for the validity of our topology-preserving Hodge decomposition. The discussion on the crucial role played by topology in the decomposition is also further elaborated. Moreover, as a downstream application for our methods, more examples on scRNA velocity decomposition have been included.

\section{Hodge decomposition in the smooth case}
\label{sec.hodgeDec.smooth}
Before describing the discretization, we briefly review the continuous theory. Let $M$ be an $m$-dimensional smooth, orientable, compact manifold with boundary. As an entity that can be integrated on $k$-submanifolds, a differential $k$-form on $M$ is an antisymmetric covariant tensor field of rank $k$ defined on $M$. The space of all differential $k$-forms on $M$ is denoted as $\Omega^k.$ The differential $d$, also called exterior derivative, is the unique $\RR$-linear mapping from $k$-forms to $(k\!+\!1)$-forms satisfying the Leibniz rule with respect to the wedge product $\wedge$ (antisymmetric tensor product) and the nilpotent property $dd = 0$. This operator generalizes and unifies several differential operators in vector analysis, including gradient, curl, and divergence. A differential form $\omega\in\Omega^k$ is \emph{closed} if $d\omega = 0$ and \emph{exact} if there is a $(k\!-\!1)$-form $\zeta\in\Omega^{k-1}$ such that $\omega = d\zeta$. Every exact form is closed due to $dd = 0$. The integral of an exact form $d\omega$ over an oriented $k\!+\!1$-submanifold $S\subset M$ with boundary $\partial S$ can be reduced to a boundary integral of $\omega$, according to Stokes' theorem, a generalization of the Newton-Leibniz rule,
\begin{align}\label{eq.stokesthm}
	\int_S d\omega= \int_{\partial S}\omega.
\end{align}
Given a Riemannian metric $g$ on $M$,  any $k$-form may be identified with a unique $(n\!-\!k)$-form through the Hodge star $\star:\Omega^k\to\Omega^{n-k}$, the unique linear operator satisfying
\begin{equation}\label{eq.innerprod.p}
	\omega\wedge\star\eta = \langle\omega,\eta\rangle_g\;\mu_g,
\end{equation}
where $\omega,\eta\in \Omega^k$, $\langle\cdot,\cdot\rangle_g$ denotes the pointwise inner product induced by $g$ on $\Omega^k,$ and $\mu_g$ is the volume form induced by $g$. By integrating Eq.~\eqref{eq.innerprod.p}, we obtain the Hodge $L^2$-inner product on the space of $k$-forms $\Omega^k$ \begin{equation}\label{eq.l2innerprod.forms}
	(\omega, \eta) = \int_M\omega\wedge\star\eta.
\end{equation}
With the differential $d$ and the Hodge star $\star$, the codifferential operator $\delta: \Omega^k\to\Omega^{k-1}$ can be defined as
\begin{equation}
	\delta = (-1)^{m(k-1)+1}\star d\star.
\end{equation}
The codifferential $\delta$ is also nilpotent, $\delta\delta = 0$. A differential form $\omega\in\Omega^k$ is \emph{coclosed} if $\delta\omega = 0$ and \emph{coexact} if there is a $(k+1)$-form $\eta\in\Omega^{k+1}$ such that $\omega = \delta\eta$. The Hodge Laplacian is then defined as $\Delta = d\delta + \delta d: \Omega^k\to\Omega^k$. A $k$-form $\omega$ is \emph{harmonic} if $\Delta\omega = 0$.

\subsection{Hodge decomposition for closed manifolds}
The classical Hodge decomposition theorem states that the space of differential $k$-forms can be decomposed orthogonally as
\begin{align}\label{eq.hd.closedMfld}
	\Omega^k = d\Omega^{k-1}\oplus\delta\Omega^{k+1} \oplus \mathcal{H}^k_{\Delta}(M),
\end{align}
where the space $\mathcal{H}^k_{\Delta}(M)$ is the finite-dimensional kernel of $\Delta$. The orthogonality can be established from
\begin{align}\label{eq.hd.duality}
(d\alpha,\beta)-(\alpha,\delta\beta)=\int_{\partial M}\alpha\wedge\star\beta,
\end{align}
for any $\alpha\in \Omega^{k-1}$ and $\beta \in \Omega^k.$
When $\partial M =\emptyset,$ the above indicates any exact (coexact) form is orthogonal to coclosed (closed resp.) forms.
The dimension of $\mathcal{H}^k_{\Delta}(M)$ is determined through the Hodge isomorphism with the $k$-th de Rham cohomology, $H^k_{dR}(M) = \ker d^k/\im d^{k-1},$ the quotient space of closed $k$-forms modulo exact $k$-forms. It follows from the de Rham theorem and Poincar\'{e} duality that the de Rham cohomology group $H^k_{dR}(M)$ is isomorphic to the $(m\!-\!k)$-singular homology group, whose dimension is given by the $(m\!-\!k)$-th Betti number $\beta_{m-k},$ where $\beta_0, \beta_1, \beta_2$ provide the numbers of connected components, tunnels, and enclosed cavities respectively. Therefore, $\dim \mathcal{H}^k_{\Delta}(M)=\beta_{m-k}.$ Note that for closed manifolds, the Hodge duality implies $\dim \mathcal{H}^k_{\Delta}(M)=\dim \mathcal{H}^{m-k}_{\Delta}(M),$ and thus $\beta_k=\beta_{m-k}.$

\subsection{Hodge decomposition for manifolds with boundary}
{\bf Orthogonality.} With a nonempty boundary, the orthogonality among the above subspaces is no longer guaranteed unless boundary conditions are imposed to the r.h.s of Eq.~\ref{eq.hd.duality}. Two such subspaces can be defined: $\Omega^k_n$ and $\Omega^k_t$ for the homogeneous normal and tangential boundary conditions, respectively,
\begin{align}
	\Omega^k_n = \{\omega\in\Omega^k\, \vert\, \omega\vert_{\partial M} = 0\}, \;
	\Omega^k_t = \{\omega\in\Omega^k\, \vert\, \star\omega\vert_{\partial M} = 0\}.
\end{align}
It follows from the definitions that the Hodge star $\star: \Omega^k_n\to \Omega^{m-k}_t$ provides an isomorphism between the two. In addition, the differential $d$ preserves the normal boundary condition, whereas the codifferential $\delta$ preserves the tangential boundary condition.
The Hodge-Morrey decomposition~\cite{morrey1956variational}, which decomposes a $k$-form into an exact normal form, a coexact tangential form, and the rest, is thus orthogonal,
\begin{align}\label{eq.morreyDecomposition}
	\Omega^k = d\Omega^{k-1}_n\oplus\delta\Omega^{k+1}_t \oplus \mathcal{H}^k,
\end{align}
where $\mathcal{H}^k = \ker d\cap \ker\delta$ is the space of closed and coclosed fields. However, $\mathcal{H}^k\subset \mathcal{H}^k_{\Delta}(M)$ is infinite-dimensional instead of being determined by the topology~\cite{schwarz2006hodge}.

\noindent{\bf Finite-dimensional kernels.}
The solution is to enforce boundary conditions on the harmonic space,
\begin{align}\label{eq.hd.fundamental}
	\Omega^k &= d\Omega^{k-1}\oplus\delta\Omega^{k+1}_t\oplus \mathcal{H}^{k}_t\\
	&= d\Omega^{k-1}_n\oplus\delta\Omega^{k+1} \oplus \mathcal{H}^{k}_n.
\end{align}
Under the normal (or tangential) boundary conditions, the subspace $\mathcal{H}^k_n= \mathcal{H}^{k}\cap \Omega^{k}_n(M)$ (or $\mathcal{H}^k_t= \mathcal{H}^{k}\cap \Omega^{k}_t(M)$, resp.) is finite-dimensional. Each normal harmonic field in $\mathcal{H}^k_n$ corresponds to a unique equivalence class in the relative de Rham cohomology $H^k_{dR}(M,\partial M)$, and each tangential harmonic field in $\mathcal{H}^k_t$ corresponds to a unique equivalence class in the absolute de Rham cohomology $H^k_{dR}(M)$, i.e., $\mathcal{H}^{k}_t \cong H^k_{dR}(M)$ and $\mathcal{H}^{k}_n\cong H^k_{dR}(M, \partial M)$ ~\cite{friedrichs1955differential}. It follows that, for compact manifolds in Euclidean spaces, the dimensions of these subspaces are given by the Betti numbers $ \dim\mathcal{H}^{k}_n = \beta_{m-k}$ and $ \dim\mathcal{H}^{k}_t = \beta_k$, which are fully determined by the topology.

\noindent{\bf Complete decomposition.}
The space $\mathcal{H}^k$ can be further decomposed into three terms \cite{friedrichs1955differential}
\begin{align}
	\mathcal{H}^k = (\mathcal{H}^{k}_n + \mathcal{H}^{k}_t)\oplus (d\Omega^{k-1}\cap\delta\Omega^{k+1}).
\end{align}
The spaces $\mathcal{H}^{k}_n$ and $\mathcal{H}^{k}_t$, in general, are not orthogonal with each other. However, they are orthogonal on compact domains in Euclidean spaces~\cite{shonkwiler2009poincare}. Thus, on such domains, a 5-component Hodge decomposition is available,
\begin{align}\label{eq.hd.5subspaces}
	\Omega^k = d\Omega^{k-1}_n\oplus\delta\Omega^{k+1}_t \oplus \mathcal{H}^{k}_n \oplus \mathcal{H}^{k}_t\oplus (d\Omega^{k-1}\cap\delta\Omega^{k+1}).
\end{align}
Note that
\begin{align}\label{eq.image.d.delta}
	\im d &= d\Omega^{k-1}_n \oplus \mathcal{H}^{k}_n \oplus (d\Omega^{k-1}\cap\delta\Omega^{k+1}),\\
	\im\delta &= \delta\Omega^{k+1}_t \oplus \mathcal{H}^{k}_t\oplus (d\Omega^{k-1}\cap\delta\Omega^{k+1}).
\end{align}
Given $\omega\in\Omega^k$, by the $5$-component Hodge decomposition \eqref{eq.hd.5subspaces}, there is a unique orthogonal 5 component decomposition~\eqref{eq.l2innerprod.forms}:
\begin{align}\label{eq.hd.5components}
	\omega = d\alpha_n+\delta\beta_t + h_n + h_t + \eta,
\end{align}
where $\alpha_n\in \Omega^{k-1}_n$, $\beta_t\in\Omega^{k+1}_t$, $h_n\in\mathcal{H}^k_n$, $h_t\in\mathcal{H}^k_t$ and $\eta\in d\Omega^{k-1}\cap\delta\Omega^{k+1}$. The first two terms can be computed by first solving for the potentials $\alpha_n\in\Omega^{k-1}_n$ and $\beta_t\in\Omega^{k+1}_t$, and then applying the differential $d$ to $\alpha_n$ and codifferential $\delta$ to $\beta_t$. Note that these potentials are not uniquely determined as $d(\alpha_n+d\xi) = d(\alpha_n)$ and $\delta(\beta_t+\delta\zeta)=\delta\beta_t$ for all $\xi\in\Omega^{k-2}_n$ and $\zeta\in\Omega^{k+2}_t$.
To ensure the uniqueness of these potentials, we impose gauge conditions $\delta\alpha_n=0$ and $d\beta = 0$. We then obtain
\begin{align}\label{eq.potentials.alpha_beta}
	\delta\omega &= (\delta d + d\delta)\alpha_n = \Delta\alpha_n\\
	d\omega &= (\delta d + d\delta)\beta_t =\Delta\beta_t.
\end{align}
By resolving the rank deficiencies of Laplacians on $\Omega^{k-1}_n$ and $\Omega^{k+1}_t$, the potential $\alpha_n$ can then be solved by considering the first equation in \eqref{eq.potentials.alpha_beta} with boundary conditions $\alpha_n|_{\partial M} = 0$ and $\delta\alpha_n|_{\partial M} = 0$, while $\beta_t$ can be solved by using the second equation with boundary conditions $\star\beta_t|_{\partial M} = 0$ and $\star d\beta_t|_{\partial M} = 0$.

Since $\mathcal{H}^k_n$ and $\mathcal{H}^k_t$ are finite-dimensional, $h_n$ and $h_t$ can be calculated by projecting the form $\omega$ onto the spaces $\mathcal{H}^k_n$ and $\mathcal{H}^k_t$ in any orthonormal bases. The last term $\eta$ can then be calculated by subtracting these four terms from $\omega$.


\subsection{Relation to vector field analysis}
When $M$ is a compact domain in $\RR^m, m=2$ or $3$ with metric induced by Euclidean metric, the differential forms can be interpreted as scalar and vector fields. Specifically, in $\RR^2$, a $0$-form or a $2$-form can be identified with a scalar field and a $1$-form can be identified with a vector field, while in $\RR^3$, a $0$-form or $3$-form can be identified with a scalar field and a $1$-form or $2$-form can be regarded as a vector field. Furthermore, there are correspondences between the differential operators $d$ and codifferential operators $\delta$ on differential forms and the gradient operator $\nabla$, the curl operator $\nabla\times$ and the divergence operator $\nabla\cdot$. The identities $dd=0$ and $\delta\delta=0$, correspond to the vector field analysis identities $\nabla\cdot\nabla\times = 0$ and $\nabla\times\nabla = 0$. See \cite{desbrun2006discrete} for a comprehensive comparison.

In the following, we provide the relations of the subspaces involved in our differential form-based decomposition with the subspaces of vector fields on $M.$
Denote by $\vf$ be the infinite-dimensional space of all smooth vector fields on $M$. A natural $L^2$ inner product on $\vf$ is given as follows:
\begin{align}\label{eq.l2innerprod.vf}
	\langle V, W\rangle = \int_MV\cdot W \; d\mathrm{Vol}.
\end{align}
It follows from the musical isomorphisms that this inner product \eqref{eq.l2innerprod.vf} on vector fields corresponds to the Hodge $L^2$-inner product on differential $1$-forms~\eqref{eq.l2innerprod.forms}. More specifically given the metric tensor $g_{ij},$ $\omega_i = g_{ij} W^j$ and $W^i=g^{ij}\omega_j,$ where we assume the Einstein summation convention and $\omega=W^\flat$, or equivalently $W=\omega^\sharp.$ Furthermore, the vector field $L^2$-product corresponds to the Hodge $L^2$-inner product on $2$-forms, through $\omega =\star W^\flat.$

In \cite{cantarella2002vector}, Cantarella et al. defined 5 subspaces of $\vf$ in $\RR^3$, called grounded gradients, fluxless knots, harmonic gradients, harmonic knots, and curly gradients, denoted by $\GG, \fk, \hg, \hk$ and $\cg$, respectively. The representations of these subspaces are given as follows:
\begin{align}\label{eq.hd.vf}
	\GG &= \{\nabla\varphi;\, \varphi\in C^{\infty}(M),\, \varphi|_{\partial M} = 0\}\\
	\fk &= \{V = \nabla\times X;\ X\in\vf,\ V\cdot\mathbf{n} = 0 \}\\
	\hg &= \{V = \nabla\varphi;\, \varphi\in C^{\infty}(M),\, V\times\mathbf{n} = 0,\, \nabla\cdot V = 0 \}\\
	\hk &= \{V = \nabla\times X;\, X\in\vf,\, V\cdot\mathbf{n}=0,\, \nabla\times V = 0 \}\\
	\cg &= \{V = \nabla\varphi = \nabla\times X;\, \varphi\in C^{\infty}(M),\, X\in\vf \},
\end{align}
where $\mathbf{n}$ denotes the unit normal vector field on the boundary $\partial M$. Then there is an orthogonal decomposition of the space of vector fields with respect to the inner product~\eqref{eq.l2innerprod.vf}
\begin{align}
	\vf = \GG\oplus\fk\oplus\hg\oplus\hk\oplus\cg.
\end{align}
This decomposition \eqref{eq.hd.vf} of vector fields is a counterpart of the Hodge decomposition \eqref{eq.hd.5subspaces} of differential forms in $\Omega^1$. The five subspaces of $\vf$ correspond exactly to the five subspaces of $\Omega^1$ in \eqref{eq.hd.5subspaces} in the same order. The same decomposition also holds for vector fields in $\RR^2$.  However, we need to redefine the subspaces $\fk, \hk$ and $\cg$, as in the 2D case the curl operator is given as $J\nabla$ with $J$ being the 2d $\pi/2$ counterclockwise rotation matrix. These three redefined subspaces are given as follows
\begin{align}
	\fk &= \{V = J\nabla\psi;\ \psi\in C^{\infty}(M),\ V\cdot\mathbf{n} = 0 \}\\
	\hk &= \{V = J\nabla\psi;\, \psi\in C^{\infty}(M),\, V\cdot\mathbf{n}=0,\, \nabla\times V = 0 \}\\
	\cg &= \{V = \nabla\varphi = J\nabla\psi;\, \varphi, \psi\in C^{\infty}(M)\}.
\end{align}
With these correspondences between vector fields and differential forms in $\RR^2$ or $\RR^3$, a $5$-component decomposition of a vector field on the domain $M$ can be straightforwardly computed by applying the Hodge decomposition on the corresponding $1$-form (or $2$-form) of the vector field.
\begin{rem}
	By considering the unions of these subspaces of $\vf$, one could also obtain $2$-component, $3$-component or $4$-component orthogonal decompositions of vector fields. In particular, the Helmholtz-Hodge decomposition is given as follows
	\begin{align}
		\vf = \GG\oplus\fk\oplus H,
	\end{align}
where $H = \hg\oplus\hk\oplus\cg$ is the union of all harmonic fields. Nevertheless, even in the continuous setting, the existence and uniqueness (up to gauge conditions) of potentials rely on the preservation of the topology in $\hg$ and $\hk$.
\end{rem}

\section{Discretization of the Hodge decomposition}
\label{sec.hodgeDec.discretized}
The generalization of discrete exterior calculus (DEC)~\cite{desbrun2006discrete} from simplicial meshes to regular  grids for topology preservation is straightforward
\begin{wrapfigure}[3]{r}{0.35\columnwidth}
\vspace*{-5mm}\hspace*{-1mm}\centering
\includegraphics[width=0.37\columnwidth]{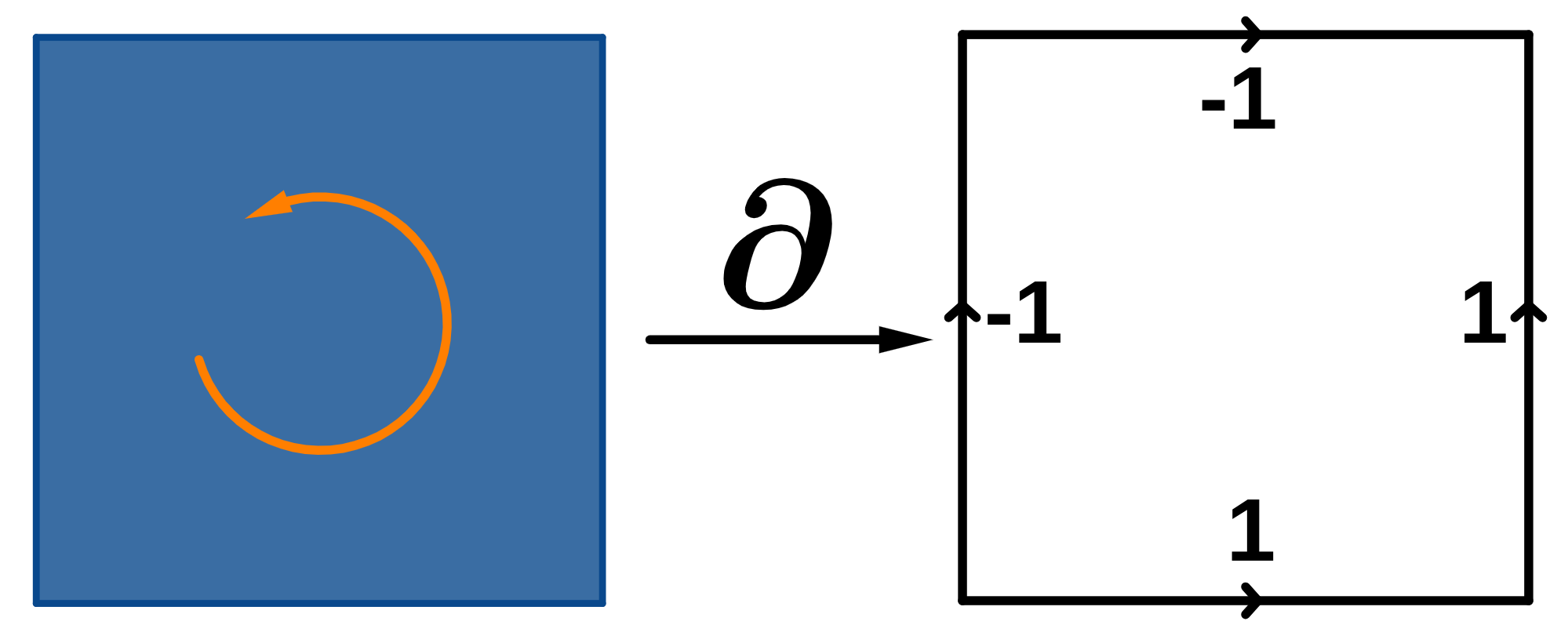}%
\end{wrapfigure}
 (see inset figure). Our focus is on keeping the topological correspondences for 2D/3D domains bounded by a level set curve/surface embedded within such grids.


\subsection{Discretization on grids}
{\bf Discrete forms.} Let $I_m$ be a rectangular $m$-dim regular Cartesian grid with $k$-cells oriented according to their alignment with the coordinate axes. With uniform grid spacing $l$ along each axis direction, each $k$-cell of the grid is a $k$-dimensional hypercube. With the grid $I_m$ treated as a cell complex tessellating a rectangular domain in $\RR^m$, a continuous $k$-form $\omega$ can be discretized by its integral value $W^i=\int_{\sigma_i} \omega$ over each $k$-cell $\sigma_i$~\cite{desbrun2006discrete}. This discretization, known as the de Rham map, linearly maps a $k$-form to a cochain, a linear functional $\int_{c} \omega=\sum_i W^i a_i$ that maps a $k$-chain $c=\sum_i a_i\sigma_i$, a formal linear combination of $k$-cells representing a $k$-dimensional subdomain, to the integral over $c$. This map is crucial since it can be seen as what induces a chain map to preserve the cohomology.

\noindent{\bf Discrete differential operator.} The discrete differential operator acting on discrete $k$-forms is represented by a sparse matrix $D^I_k.$ This matrix encodes the signed incidence between $(k\!+\!1)$-cells and $k$-cells, as in the simplicial case, i.e., $D^I_k=\partial_{k+1}^T$, the transpose of the cell boundary operator $\partial_{k+1}$ on $(k\!+\!1)$-cells. This is the direct consequence of Stokes' theorem $\int_{\sigma} d\omega= \int_{\partial \sigma} \omega$. It follows that $D^I_{k+1}D^I_k = 0,$ since the boundary of a boundary of any cell is a 0 chain (i.e., $\partial\partial = 0$).

\noindent{\bf Discrete Hodge star.} Treating the centers of $m$-cells as grid points, we construct a \emph{dual} grid that is staggered with the \emph{primal} grid $I_m$. A one-to-one correspondence between discrete $k$-forms on the primal grid and $m\!-\!k$-form on the dual grid can be established by the diagonal discrete Hodge star $\star_k$, induced by the continuous Hodge star through local averaging
\begin{align}
	\frac{1}{|\sigma_k|}\int_{\sigma_k}\omega \approx \frac{1}{|\star\sigma_k|}\int_{\star\sigma_k}\star\omega,
\end{align}
where $\star\sigma_k$ is the dual $(m\!-\!k)$-cell formed by the dual grid points associated with the $m$-cells incident to the $k$-cell $\sigma_k$. Equivalent to assuming a one-point quadrature for the integration involved in discretizing primal and dual forms, this correspondence leads to a diagonal matrix $S^I_k$ with diagonal entries given by the ratio between the volumes of the dual  $(m\!-\!k)$-cells and primal $k$-cells, $l^{m-k}/l^k=l^{m-2k}$. The associated discrete Hodge $L^2$ inner product \eqref{eq.l2innerprod.forms} of two discrete $k$-forms $V_k$ and $W_k$ on the grid $I_m$ is then
\begin{align}
	(V_k, W_k)^I = V_k^TS^I_kW_k.
\end{align}
While this treatment can be made more accurate with a nondiagonal Hodge star matrix, as we discuss below, it does not affect the topology preservation.

\noindent{\bf Discrete Laplacian.} The discrete codifferential operator $\delta$ can be assembled from the discrete differential and Hodge star operators as $(S^I_{k-1})^{-1}D^I_{k-1}S^I_k$. With the differential and codifferential in $\Delta = d\delta + \delta d$ replaced by their discrete counterparts, the resulting matrix is nonsymmetric. Therefore, the discrete Laplacian is defined as the counterpart of $\star\Delta$,
\begin{align}
	L^I_k = (D^I_k)^TS^I_{k+1}D^I_k + S^I_kD^I_{k-1}(S^I_{k-1})^{-1}(D^I_{k-1})^TS^I_k,
\end{align}
where the operator is considered null for $k<0$ or $k>m$.

\subsection{Discrete differential forms and operators on $M$}
For simplicial or polygonal meshes, identifying boundary elements is straightforward, allowing easy implementation of projection matrices to the interior or boundary of the domain $M$. However, when $M$ is defined by the volume enclosed within a level set surface, enforcing boundary conditions through projection matrices becomes challenging. Instead of tessellating the boundary cells to form new unstructured meshes, e.g., through Marching Cubes, we modify the Hodge star operators. This approach maintains consistent data structures, accommodating evolving level sets and eliminating the need for remeshing.

\begin{figure}[t]
	\includegraphics[width=0.43\columnwidth]{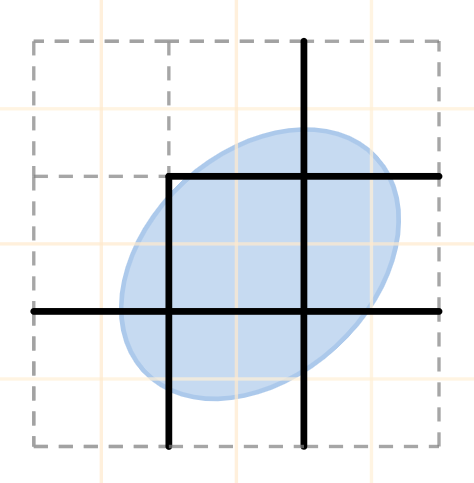}\hfil
	\includegraphics[width=0.43\columnwidth]{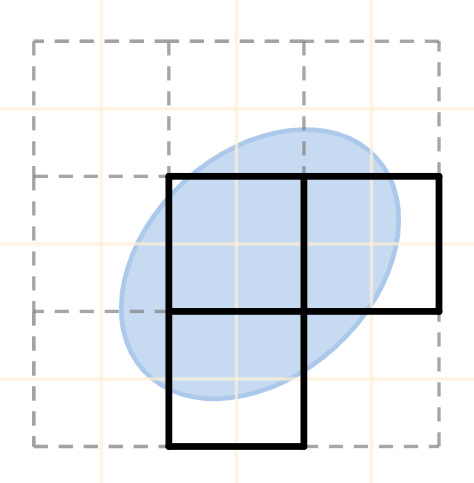}
	\caption{Distinction in supports for primal $1$-forms under normal boundary condition (left) and for those under tangential boundary condition (right).}
	\label{fig.csupports.NT}
\end{figure}

\noindent{\bf Compact supports.} While we do not cut the boundary cells, we still need to restrict the computation to the relevant cells through the inclusion or exclusion of the entire $k$-cells, similar to voxelization. However, boundary (primal or dual) $k$-cells typically intersect the boundary rather than being completely contained within it. As we aim to construct Laplacian operators with rank deficiencies corresponding to the topology of $M$, we design one rule for each of the two types of boundary conditions. For normal boundary conditions, we include every cell with at least one vertex inside $M$, while for tangential boundary conditions, we select every cell with at least one vertex of its dual cells inside $M$. The set of cells for the former is called the \emph{normal support} and the set of cells for the latter is called the \emph{tangential support}. Note that, the normal and tangential supports are typically distinct, and, unlike the mesh case, neither is necessarily a superset of the other. See Fig.~\ref{fig.csupports.NT} for one example highlighting the differences in the two supports for $1$-forms.

The 0-1 projection matrices $P_{k, n}$ and $P_{k,t}$ map $k$-chains to the normal and tangential support respectively. They can be derived from the identity matrix by eliminating the rows corresponding to $k$-cells outside the support. The relevant differential operators are
\begin{align}
	D_{k, n} = P_{k+1, n}D^I_k P_{k, n}^T, \;
	D_{k, t} = P_{k+1, t}D^I_k P_{k, t}^T.
\end{align}
The nilpotent property $D_{k+1, n}D_{k, n}= 0$ and $D_{k+1, t}D_{k, t}= 0$ remains for both, due to $D^I_{k+1}D^I_k=0$ and the following observations,
\begin{align}\label{eq:projectionRelation}
	P_{k+1,n}^T P_{k+1, n}D^I_k P_{k, n}^T = D^I_k P_{k, n}^T,\;
	P_{k+1, t} D^I_k P_{k, t}^T P_{k,t} = P_{k+1,t} D^I_k.
\end{align}

\noindent{\bf Modified Hodge stars and Laplacians.} Our discretization extends the formulations in \cite{batty2007fast,ng2009efficient,liu2015model,ribando2022graph} to accommodate the 5-component Hodge decomposition. This complete decomposition necessitates the simultaneous use of both tangential fields and normal fields. One key observation is that $d\Omega_n$ and $\delta\Omega_t$ are essential to the decompositions (Eqs.~\eqref{eq.hd.fundamental} and~\eqref{eq.hd.5subspaces}). Thus, while the forms on normal and tangential supports may be ``voxelized,'' their differentials and codifferentials, respectively, need to effectively approximate elements of $d\Omega_n$ and $d\Omega_t$.

To this end, we retain the dual cell volumes while adjusting the primal cell volumes for normal boundary conditions, and do the opposite for \begin{wrapfigure}[8]{r}{0.32\columnwidth}
\vspace*{-5mm}\hspace*{-2mm}\centering
\includegraphics[width=0.32\columnwidth]{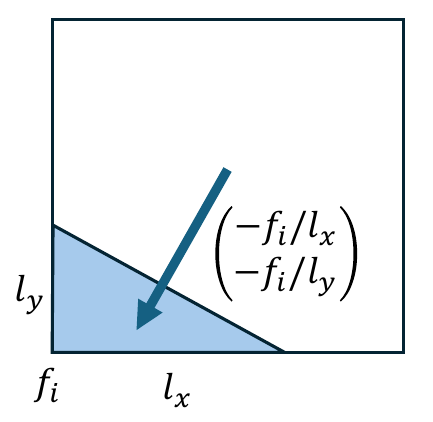}%
\end{wrapfigure}
tangential boundary conditions: keep the primal cell volumes while modifying the dual cell volumes. For instance, a gradient field $\nabla f$ of a scalar function $f$ fixed to 0 on the boundary can be represented by $D_{0,n}F$ of a normal $0$-form $F$. Its value in a boundary cell
will be forced to be orthogonal to the boundary as shown in the inset
figure. Equivalently, the codifferential of a $3$-form with tangential support also satisfies the tangential condition of the resulting $2$-form, which also corresponds to a gradient vector field orthogonal to the boundary.

Specifically, for normal (or tangential) boundary conditions, we replace the $k$-volume of each primal (resp. dual) $k$-cell. Instead of the full $k$-volume of $l^k$, where $l$ is the grid spacing, we use the $k$-volume of its intersection with $M$. For the discrete Hodge star matrix, this modified volume is placed in the denominator (or numerator) while the dual (or primal) cell volumes $l^{n-k}$ are left unchanged in the numerator (or denominator). For numerical stability, we perturb the level set function evaluated at primal/dual gridpoints to have an absolute value above $\epsilon=10^{-5}l$, which ensures fractional $k$-volumes behave well. We denote by $S_{k, n}$ and $S_{k, t}$ the resulting sparse Hodge star matrices defined on the normal and tangential supports, respectively.

The discrete $L^2$ inner products of the two types of discrete $k$-forms on the manifold $M$ under these two boundary conditions, namely,
the discrete $\Omega^{k}_n(M)$ and $\Omega^{k}_t(M)$, are thus
\begin{align}
	(\xi^k,\, \zeta^k)^n = (\xi^k)^TS_{k,n}\zeta^k 
 ,\;\;\;
	(\xi^k,\, \zeta^k)^t = (\xi^k)^TS_{k,t}\zeta^k. \label{eq.l2innerProd.discrete.t}
\end{align}

Finally, we assemble the discrete Hodge Laplacians as follows:
\begin{align}
	L_{k, n} &= D_{k, n}^TS_{k+1, n}D_{k, n} + S_{k, n}D_{k-1, n}S_{k-1, n}^{-1}D_{k-1, n}^TS_{k, n}\\
	L_{k, t} &= D_{k, t}^TS_{k+1, t}D_{k, t} + S_{k, t}D_{k-1, t}S_{k-1, t}^{-1}D_{k-1, t}^TS_{k, t}.
\end{align}
These two types of discrete Hodge Laplacians are crucial for computing potentials in our implementation of the topology-preserving Hodge decomposition in Eq.~\eqref{eq.morreyDecomposition}, not only because of the built-in boundary conditions, but also the resulting topology-determined rank deficiency.

\subsection{Vector field decomposition}
\label{sec.hhd.vecfield}
In 2D or 3D, vector fields on $M$ can be represented as discrete $1$-forms. Following the typical de Rham map $W^i=\int_{\sigma_i}\omega $, the integral of the vector field along each grid edge is the corresponding primal $1$-form (on the normal support). For instance, if the input vector field is sampled as one vector per grid point in $M$, and $0$ at grid points outside, the discrete $1$-form on each primal edge in the tangential support is the average of the edge direction component of vectors on the inside endpoints of each grid edge, multiplied by the fraction of the primal edge length within the domain. The resulting discrete $1$-form is denoted as $W_n.$

To reconstruct the field at a specific point, e.g., in the 3D case, each component of the vector can be evaluated from bilinear interpolation of the average values on the 4 edges along that direction incident to the containing grid cell. This can be seen as the Whitney map on the Cartesian grid, similar to the one that constructs a continuous field from a discrete $1$-form on a simplicial mesh.

In 3D, it is possible to use a discrete $2$-form instead. While the DoFs for $1$-forms and $2$-forms can be drastically different in the simplicial mesh case, in our Cartesian representation, it is strictly equivalent to a discrete $1$-form on a grid shifted by $(l/2,l/2,l/2)$. Thus, we limit our discussion to the $1$-form representation.

\begin{rem}
More precisely, a vector field $v$ can also be represented as a $2$-form on $M$ by averaging the normal components of vectors on inside grid points of each face times the face area inside the domain. The decomposition of the $2$-form following \eqref{eq.morreyDecomposition} provides a dual version, where the Laplacians $L_{1, n}$ and $L_{3,t}$ will be used for solving the vector and scalar potentials, respectively. One could use either of these two representations of vector fields for the Hodge decomposition since they are equivalent to one another through the duality between a Cartesian grid and its staggered dual by offsetting each grid point with $(l/2,l/2,l/2)$, so long as $M$ is at least one grid spacing away from the boundary of the grid.
\end{rem}

For the tangential representation of the 1-form, we actually follow the discretization of the normal $2$-form as described above on the dual grid, except that the average is rescaled by $l$ instead of by the inside part of the face area. This ensures that the tangential $1$-form $W_t$ samples the vector field at least at one point within the domain, as at least one of the four neighboring dual cell centers will be inside $M$. Recall that while $W_n$ and $W_t$ are merely represented in the normal and tangential supports, respectively,  neither necessarily satisfy the corresponding boundary conditions.

In the following, we first describe a naive approach for calculating the five components in the Hodge decomposition Eq.~\eqref{eq.hd.5components}, which converges to orthogonal components only in the continuous limit. Then, we describe a modified computation that achieves discrete \l-orthogonality of the 5 subspaces, mirroring the orthogonality found in mesh-based decomposition.

\subsubsection{Direct approach}

Similar to the mesh setting, to solve for the scalar potentials $A_n$ defined on vertices (i.e., $0$-cells) and $B_t$ defined on faces (i.e., $2$-cells) in the Hodge decomposition Eq.~\eqref{eq.hd.5components}, we use the discrete equivalents of Eq.~\eqref{eq.potentials.alpha_beta}:
\begin{align}\label{eq.potential.discrete}
	L_{0, n}A_n &= D_{0,n}^TS_{1,n} W_n, \quad L_{2, t}B_t = S_{2,t}D_{1,t} W_t.
\end{align}

The normal Laplacian $L_{0,n}$ is full rank, but the kernel size of $L_{2,t}$ in 3D is $\beta_2$, as determined by the topology of (voxelized) $M$. To address the rank deficiency, we add a small positive value to $\beta_2$ selected diagonal entries to $L_{2,t}$. Following the computation of the potentials, the first two terms in Hodge decomposition \eqref{eq.hd.5components} are calculated by applying the discrete differential $D_{0,n}$ to $A_n$ and the discrete codifferential $\delta_{2,t} = S_{1,t}^{-1}D_{1,t}^TS_{2,t}$ to $B_t$.

To compute the normal harmonic component $N_h$ and the tangential harmonic component $T_h$, we simply project $W$ to the kernels of the discrete Laplacians $L_{1,n}$ and $L_{1,t}$. As with the mesh-based approach, these two Laplacians correctly capture the cohomology of the underlying manifold: as the modification to the Hodge stars does not change the cohomology, their kernels have the same dimensions as those of the voxelized $M$, which are homotopic to $M$. Thus, the space of discrete normal harmonic $1$-fields is isomorphic to the first relative cohomology group, and the space of discrete tangential harmonic $1$-fields isomorphic to the first absolute cohomology group, i.e.,
\begin{align}
	\ker L_{1,n}\cong H_{dR}^1(M, \partial M), \qquad \ker L_{1,t}\cong H_{dR}^1(M).
\end{align}
It follows that $\dim\ker L_{1,n} = \beta_{m-1}$ and $\dim\ker L_{1,t} = \beta_1$. The $\beta_{m-1}$ (resp. $\beta_1$) eigenvectors corresponding to the $0$ eigenvalues of $S_{1,n}^{-1}L_{1,n}$ (resp. $S_{1,t}^{-1}L_{1,t}$) form a basis of the normal (resp. tangential) harmonic space $\mathcal{H}^1_n$ (resp. $\mathcal{H}^1_t$). Let $\mathbb{H}_{1,n}$ (resp. $\mathbb{H}_{1,t}$) be the matrix with columns formed by the basis elements of $\mathcal{H}^1_n$ (resp.$\mathcal{H}^1_t$). The projections are
\begin{align}
	N_h = \mathbb{H}_{1,n}\mathbb{H}_{1,n}^TS_{1,n} W_n, \qquad T_h = \mathbb{H}_{1,t}\mathbb{H}_{1,t}^TS_{1,t} W_t.
\end{align}

To compute the final term $E$ in \eqref{eq.hd.5components}, it is important to account for the different dimensions between the discrete normal forms $D_{0,n} A_n, N_h$ and tangential forms $S_{1,t}^{-1}D_{1,t}^TS_{2,t} B_t, T_h$. This difference in dimensions arrives due to the relation between the normal and tangential support. Consequently, one has to resolve the consistency issue, e.g., by converting all representations to tangential support.
To avoid excessive averaging, we construct the linear conversion operator $C_{n\to t}$ as follows: if an edge is in both supports, we simply rescale the normal 1-form value by $l^{2-m} S_{1,n},$ otherwise, the vectors at the incident cell centers are reconstructed from the normal 1-form (i.e., Whitney map), then the tangential 1-form discretization procedure is carried out on that edge.
With this conversion operator, we obtain
\begin{align}
	E = W_t &- C_{n\to t} D_{0,n} A_n
	 - \delta_{2,t} B_t - C_{n\to t} N_h - T_h.
\end{align}
However, the resulting 5 components are only $L^2$-orthogonal in the continuous limit.

\subsubsection{Discrete orthogonal decomposition}

To achieve a discrete L2-orthogonal decomposition, we first have to choose a consistent L2-inner-product and a consistent functional space. Given that the Hodge duality on regular Cartesian grids is straightforwardly established by shifting the grid by half a grid spacing, we use $S_{1,t}$ and the tangential support for the following discussion without loss of generality.

\noindent{\bf Tangential decomposition.} According to Eq.~\eqref{eq.hd.fundamental}, a 3-component discrete orthogonal decomposition containing two of the components we need is
\begin{align}\label{eq.hd.discrete.fundamental.1t}
	W_t = D_{0,t} A_t + \delta_{2,t} B_t + T_h,
\end{align}
where $A_t$ and $B_t$ can be solved from
\begin{align}
	L_{0, t} A_t &= D_{0,t}^TS_{1,t}W_t,\quad L_{2,t}B_t = S_{2,t}D_{1,t}W_t,
\end{align}
where the rank deficiency of $L_{0,t}$ is $\beta_0$, the number of connected components, which can be fixed in the same way as $L_{2,t}$.  The third term $T_h$ can be obtained either from the projection procedure Eq.~\eqref{eq:projectionRelation}, or by
\begin{align}
	T_h = W_t - D_{0,t}A_t - \delta_{2,t}B_t.
\end{align}

\noindent{\bf Gradient field decomposition.}
We further decompose $D_{0,t}A_t$ into three components according to Eq.~\eqref{eq.image.d.delta},
\begin{align}\label{eq.hd.discrete.fundamental.2n}
	D_{0,t} A_t = D_{0,t}\tilde{A}_n + \tilde{N}_h + \tilde{E} =D_{0,t}(\tilde{A}_n+\tilde{A}_h+\tilde{A}_E),
\end{align}
where $\tilde{A}_n$ is a normal 0-form on the tangential support, approximating a function that vanishes at the boundary, similar to $A_n$ except not defined on the normal support; $\tilde{N}_h=D_{0,t}\tilde{A}_h$ is similar to $N_h$ but on tangential support with the potential $\tilde{A}_h$; and $\tilde{E}=D_{0,t}\tilde{A}_E$ is similar to $E$ but defined on tangential support with its scalar potential $\tilde{A}_E$.

To implement this, we build two nested linear subspaces of $\im D_{0,t}.$ The first subspace, representing the discrete version of ${\im D_{0,t}\cap \mathcal{H}^1}$, is the space onto which an \l-projection will produce $\tilde{A}_h +\tilde{A}_E$. The second space is a subspace of the first, denoted as $\widetilde{H_{1,n}}$, corresponding to an \l-projection providing $\tilde{A}_h$. The final component can then be computed through $\tilde{A}_n=A_t-(\tilde{A}_h+\tilde{A}_E).$ All these potentials are unique, up to one constant shift per connected component (kernel of $L_{0,t}$).

As the first linear subspace is the space of harmonic exact 1-forms, we seek to express it through discrete harmonic 0-form potentials, which are determined by its restriction to the boundary. However, in contrast to the mesh-based case, not all boundary grid points are present in the tangential support of $0$-forms. Therefore, we establish an extended support for the harmonic potential $A$, which includes all grid points incident to any grid edge in either tangential or normal support. With the projection $P_{E\to T}$ representing the projection from the extended support to the tangential support, $\tilde{A}_E+\tilde{A}_h= P_{E\to T} A$ is the harmonic potential whose differential leads to $\tilde{E} + \tilde{N}_h$.  The linear space is thus defined as $\{A\;|\; L_{0,E} A = 0\},$ where $L_{0,E}$ is the graph Laplacian for 0-forms in the extended support evaluated on all grid points in the normal support. The graph Laplacian is necessary and sufficient here since we are enforcing neither the tangential nor the normal boundary condition for this component.

The projection to the first linear subspace can be obtained from the linear system resulting from the constrained minimization,
\begin{align}\label{eq.constrained.minimization}
	\min_{A} \|D_{0,t} P_{E\to T} A - D_{0,t} \tilde{A}_n\|_t^2
 + \Lambda^T L_{0,E} A,
\end{align}
where $\|W\|_t = \sqrt{W^T S_{1,t} W}$ is the L2-norm for a 1-form $W$ in tangential support, and $\Lambda^T$ is the Lagrange multiplier to enforce the constraint $L_{0,E} A = 0$. To eliminate rank deficiency in the linear system, $A$ is set to 0 at one vertex per connected component within tangential support.

The second subspace is constructed as a subspace of normal harmonic forms, which is a subspace of the space of exact harmonic forms. As the dimension is low ($\beta_2$), it is efficient to construct its basis. We first find, for the $i$-th basis 1-form $N_{h,i}$ for the harmonic normal 1-form space, the closest $1$-form $D_{0,t}\bar{A}_{h,i}$ within the tangential support can be obtained by first solving
\begin{align}
L_{0,t} \bar{A}_{h,i} = D_{0,t}^T S_{1,t} C_{n\to t} N_{h,i},
\end{align}
We then project $D_{0,t}\bar{A}_{h,i}$ to the first subspace using Eq.~\ref{eq.constrained.minimization}, and perform the Gram-Schmidt procedure to form an orthonormal basis for $\widetilde{H_{1,n}}.$ The resulting basis 1-forms can be assembled as $\tilde{\mathbb{H}}_{1,n}.$ The projection is thus
\begin{align}
	\tilde{N}_h = \tilde{\mathbb{H}}_{1,n}\tilde{\mathbb{H}}_{1,n}^TS_{1,t} D_{0,t} P_{E\to T} A.
\end{align}



\subsection{Kernel dimensions and $L_2$-orthogonality}
In this section, we elaborate on the aforementioned cohomology preservation and the related \l-orthogonality among the spaces.
\noindent{\bf Voxelized cell complex.} In contrast to the Lagrangian case in \cite{zhao20193d}, the discrete version of de Rham's theorem on the isomorphism between homology and de Rham cohomology is not immediately apparent on the normal/tangential supports of Cartesian grids. However, while the ``voxelized'' supports do not follow the actual boundary surface, the cohomology ($\ker D/ \im D$) still exists as $D_t D_t=0$ ($D_n D_n = 0$) still holds, and it depends only on $D_t$ (or $D_n$) but not on $S_t$ (or resp. $S_n$). In fact, as the tangential support contains each primal $k$-cell that has a dual $(n\!-\!k)$-cell with at least one internal dual grid point, it forms a voxelized cell complex mesh. This is due to that any $k$-cell in the tangential support has each of its $(k\!-\!1)$-faces also in the tangential support, because the dual $(n\!-\!k\!+\!1)$-cell has at least one internal dual grid point. Thus, $\ker D_{k,t}/ \im D_{k-1,t}$ is isomorphic to $H^k(M)$, since the voxelized mesh is homeomorphic to $M$. The cohomology on the normal support $\ker D_{k,n}/ \im D_{k-1,n}$ likewise corresponds to a voxelized dual cell complex, and is thus isomorphic to $H_{n-k}(M)\cong H^k(M,\partial M)$.

As $S_t$ and $S_n$ are nonsingular by construction, it then follows that $\ker L_{k,t}\cong H^k(M)$ and $\ker L_{k,n} \cong (M,\partial M)$, as $\ker L_k$ contains one unique representative from every equivalence class $[W+ D F]$ in $\ker D_k/\im D_{k-1,t}$. For instance, if both $W$ and $W+D_{0,t}F$ belong to $\ker L_{1,t}$, then $\delta_{1,t} D_{0,t} F=0;$ so $F$ belongs to the kernel of $L_{0,t}$ (constant functions on each connected component), thus $D_{0,t}F=0.$

\noindent{\bf Orthogonality.}
The naive approach does not have strict discrete orthogonality among the 5 subspaces as confirmed in our experiments. In the following, we show why our projection-based approach establishes the $S_{1,t}$-orthogonality in the decomposition on the tangential support.

\begin{figure*}[t]
	\centering
	\includegraphics[trim=3cm 3cm 3cm 4cm, clip, width=0.12\columnwidth]{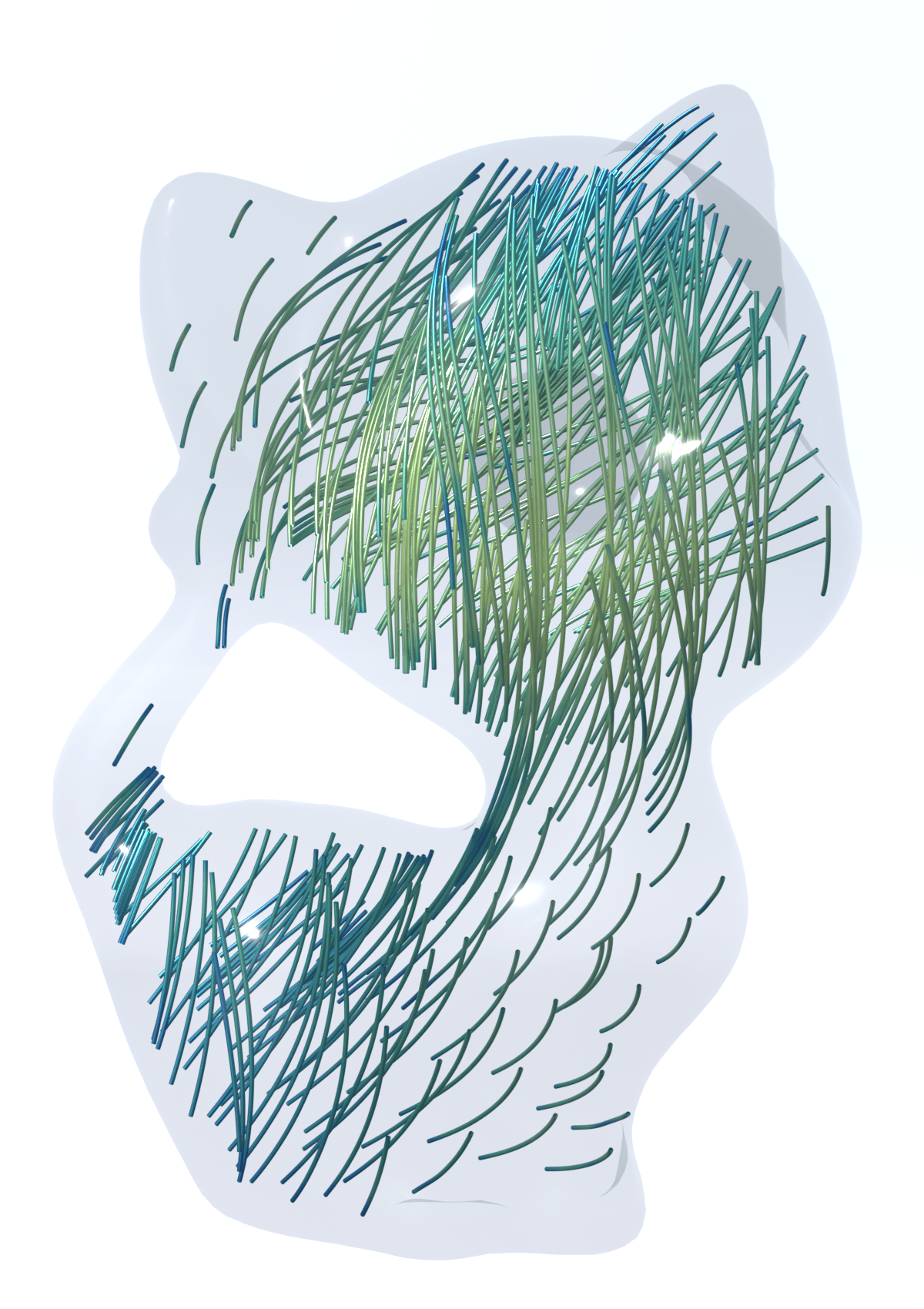}
 \raisebox{1cm}{\Large$=$}
	\includegraphics[trim=3cm 3cm 3cm 4cm, clip, width=0.12\columnwidth]{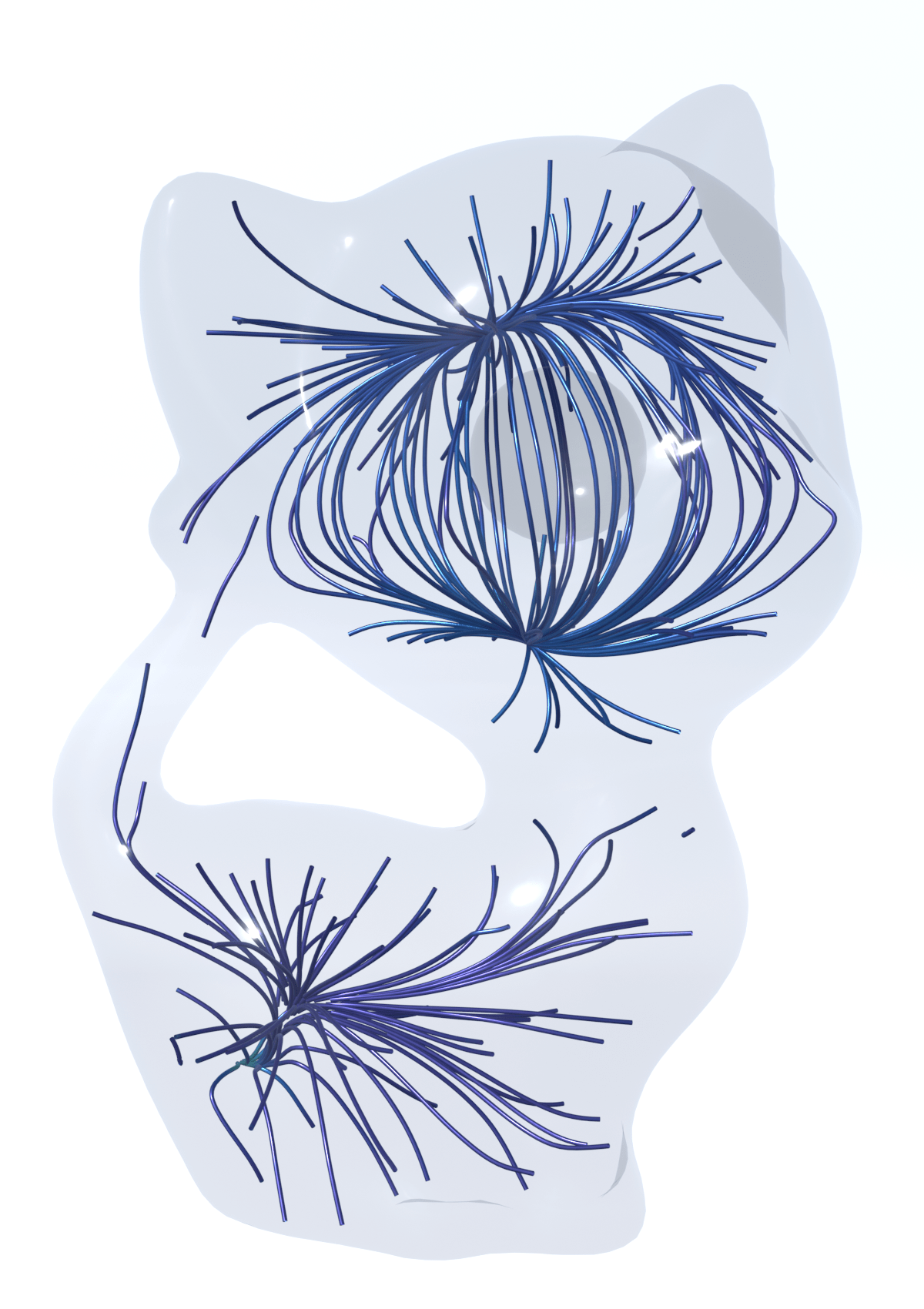}
  \raisebox{1cm}{\Large$+$}
	\includegraphics[trim=3cm 3cm 3cm 4cm, clip, width=0.12\columnwidth]{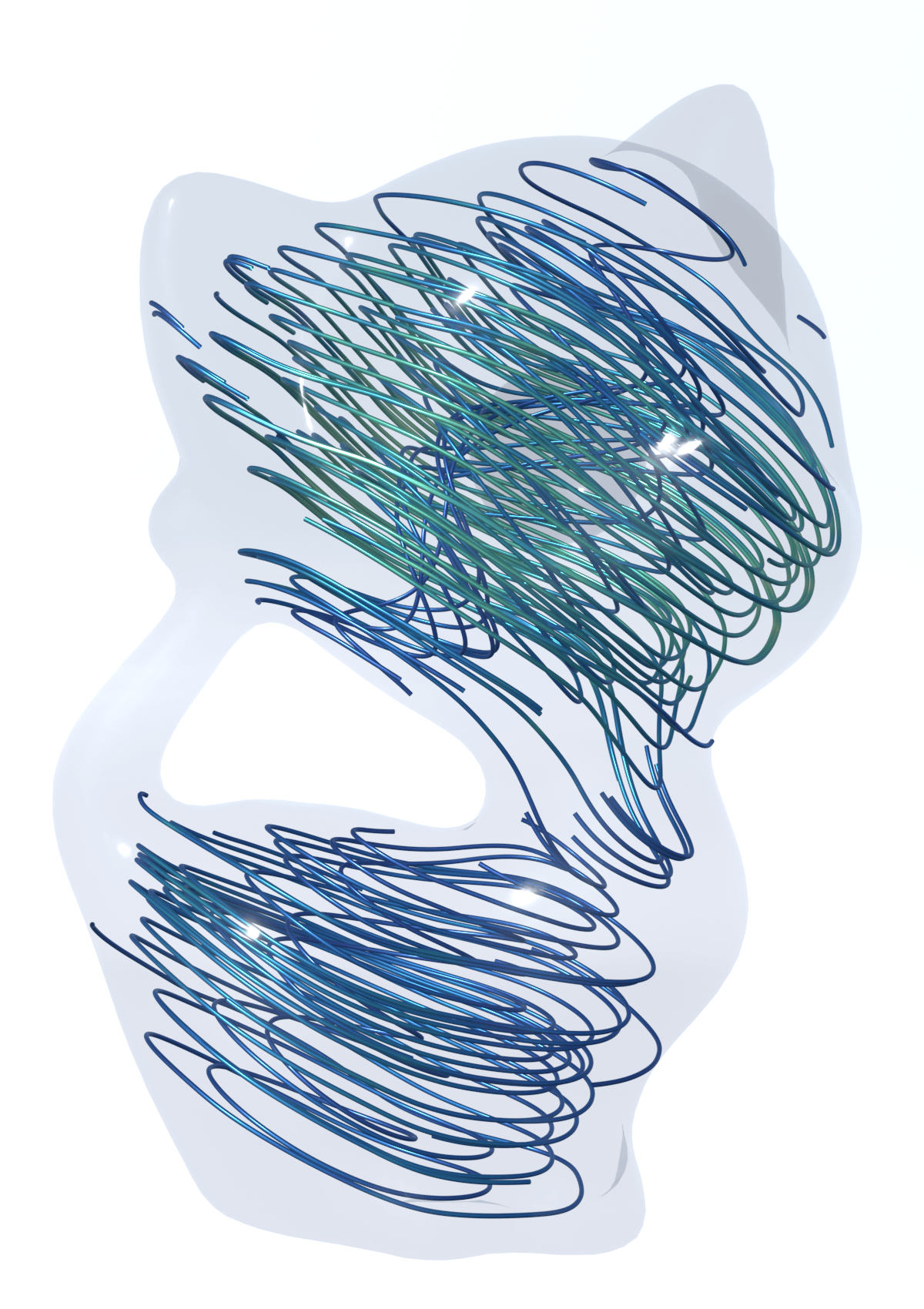}
  \raisebox{1cm}{\Large$+$}
	\includegraphics[trim=3cm 3cm 3cm 4cm, clip, width=0.12\columnwidth]{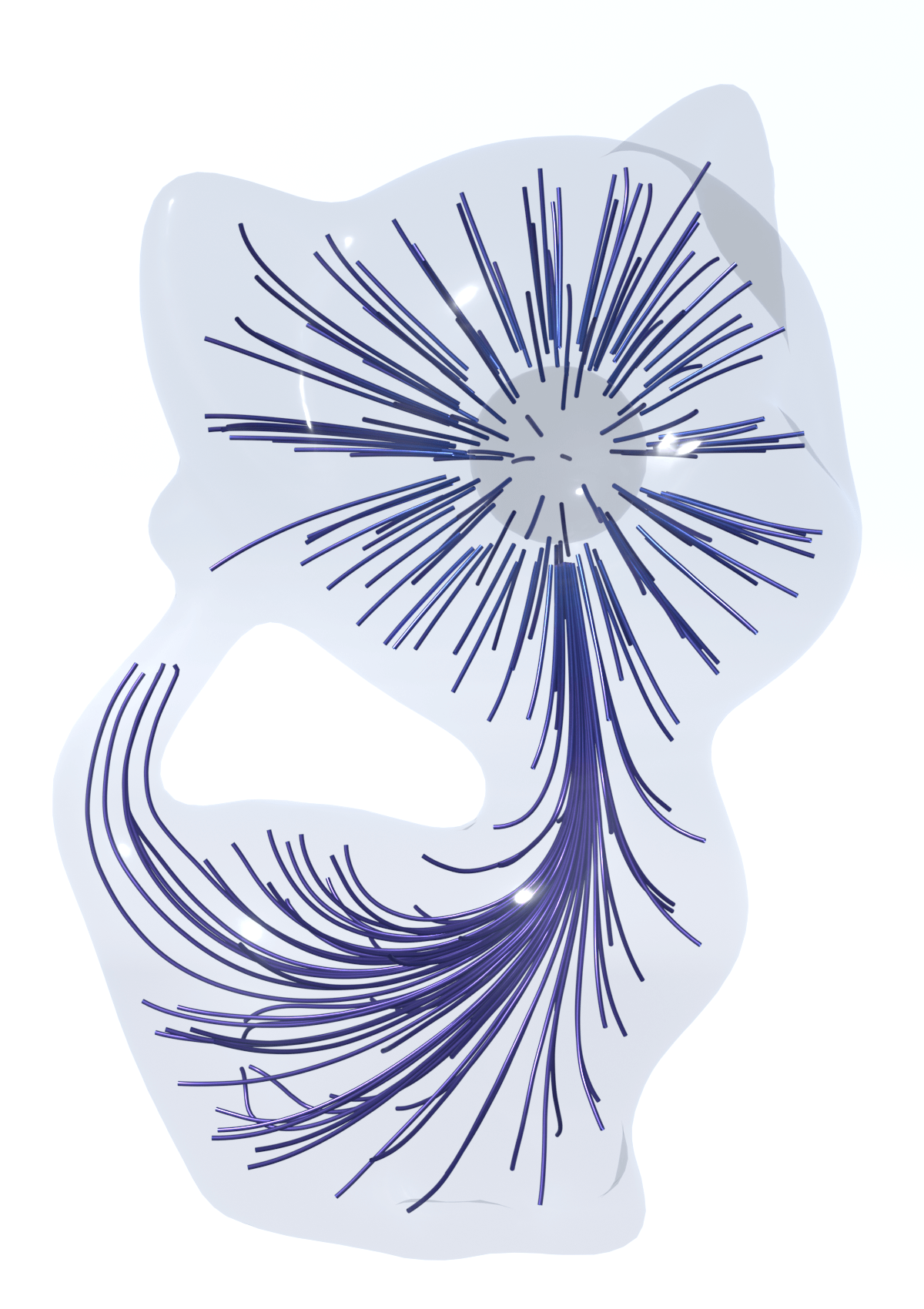}
  \raisebox{1cm}{\Large$+$}
	\includegraphics[trim=3cm 3cm 3cm 4cm, clip, width=0.12\columnwidth]{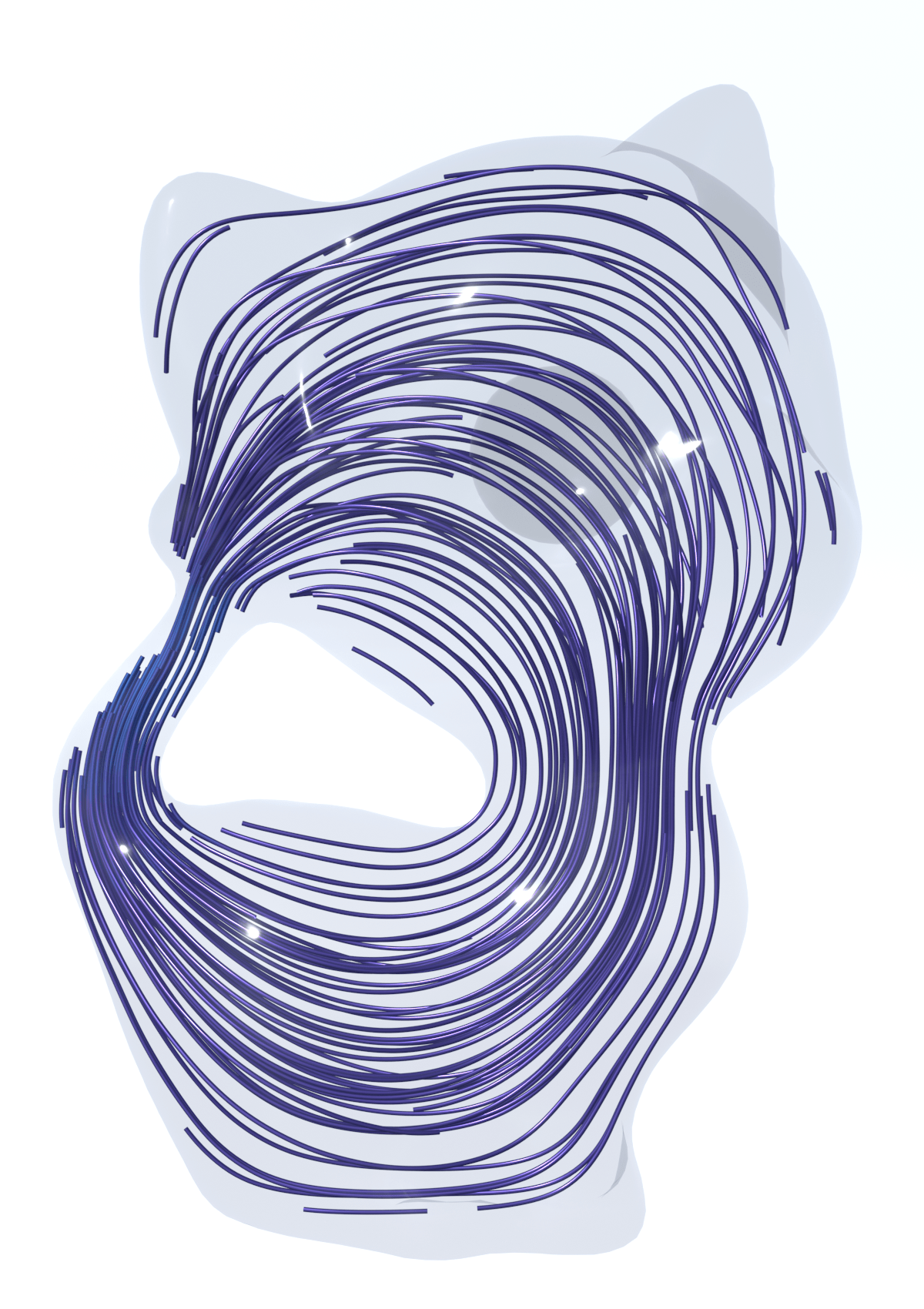}
  \raisebox{1cm}{\Large$+$}
	\includegraphics[trim=3cm 3cm 3cm 4cm, clip, width=0.12\columnwidth]{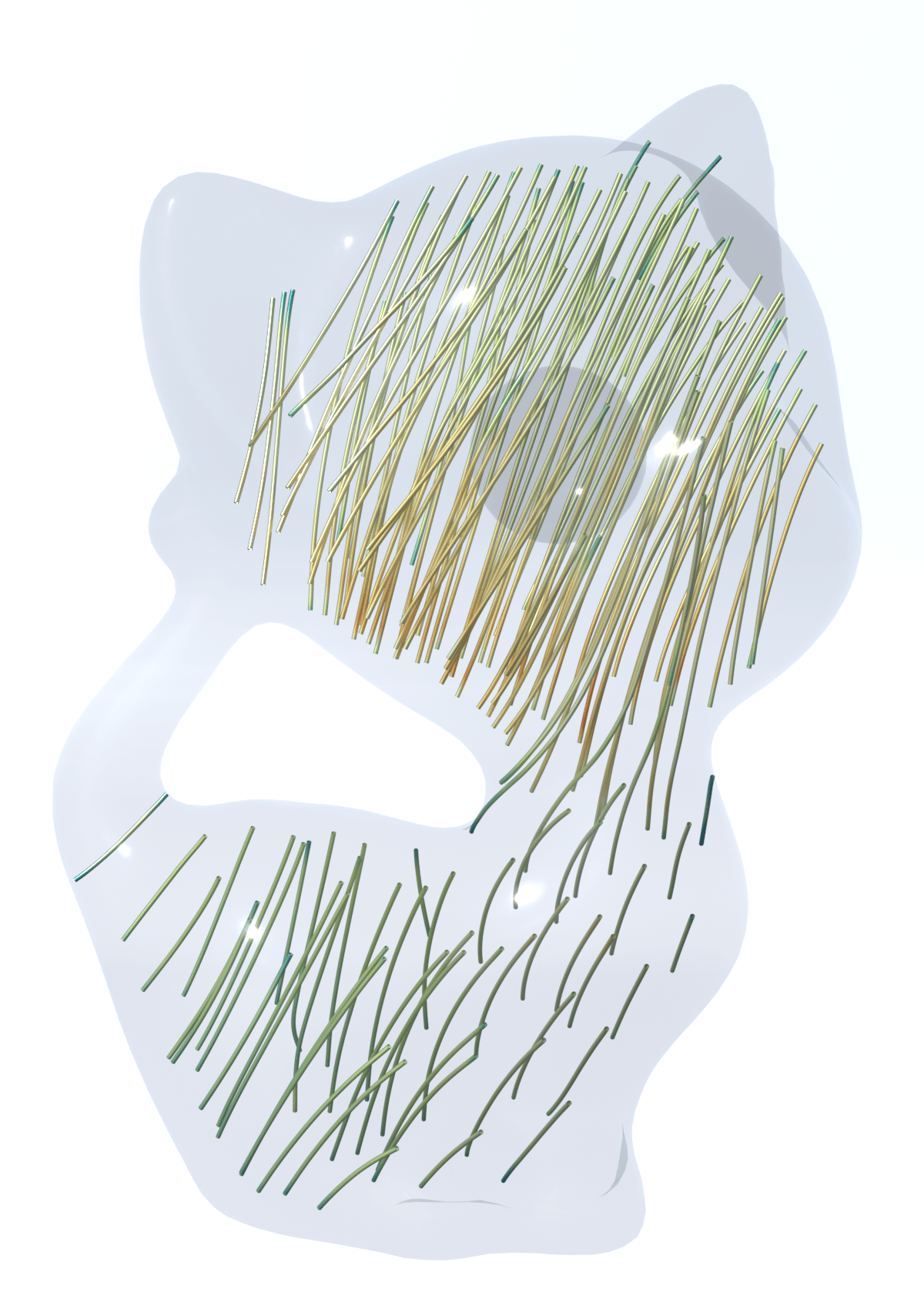}
	\caption{3D Hodge decomposition on kitty model. From left to right: the original vector field, the normal gradient field, the tangential curl field, the normal harmonic field, the tangential harmonic field, and the curly gradient field.}
	\label{fig.hd.3d.kitty}
\end{figure*}

By mimicking the Eq.~\eqref{eq.hd.duality} on tangential representations, we show the orthogonality among $D_{0,t}A$, $\delta_{2,t}B_t$ and $T_h$ as follows
\begin{align}
	(D_{0,t}A)^T S_{1,t} \delta_{2,t} B_t &= A^T S_{0,t} (\delta_{1,t} \delta_{2,t}) B_t =0,\\
    (D_{0,t}A)^T S_{1,t} T_h &= A^T S_{0,t} (\delta_{1,t} T_h) =0,\\
    (\delta_{2,t} B)^T S_{1,t} T_h &= B^T S_{2,t} (D_{1,t} T_h) =0.
\end{align}
The first orthogonality results from the nilpotent property of $\delta,$ and the second and third from the closed and coclosed property of $T_h.$
The orthogonality of the 5-component decomposition then follows from the fact that we constructed 3 mutually orthogonal subspaces of $\im D_{0,t}$.

\section{Numerical experiments}

Our algorithm was implemented in MATLAB and tested on a laptop with 16GB memory. The images were rendered in Blender. To demonstrate the effectiveness, we provide both 2D and 3D examples of Hodge decomposition of vector fields defined on compact domains. These domains are represented by signed distance functions (SDF) on regular Cartesian grids. Denoting by $\rho$ the level set function, the compact domain is
\begin{align}\label{eq.manifold}
	M = \{x\;|\; \rho(x)\leq 0 \},
\end{align}
with the boundary given by $\partial M = \{x\;|\; \rho(x)= 0\}$.

\begin{figure}[ht]
	\centering
	\includegraphics[trim=0cm 4cm 0cm 4cm, clip, width=0.12\textwidth]{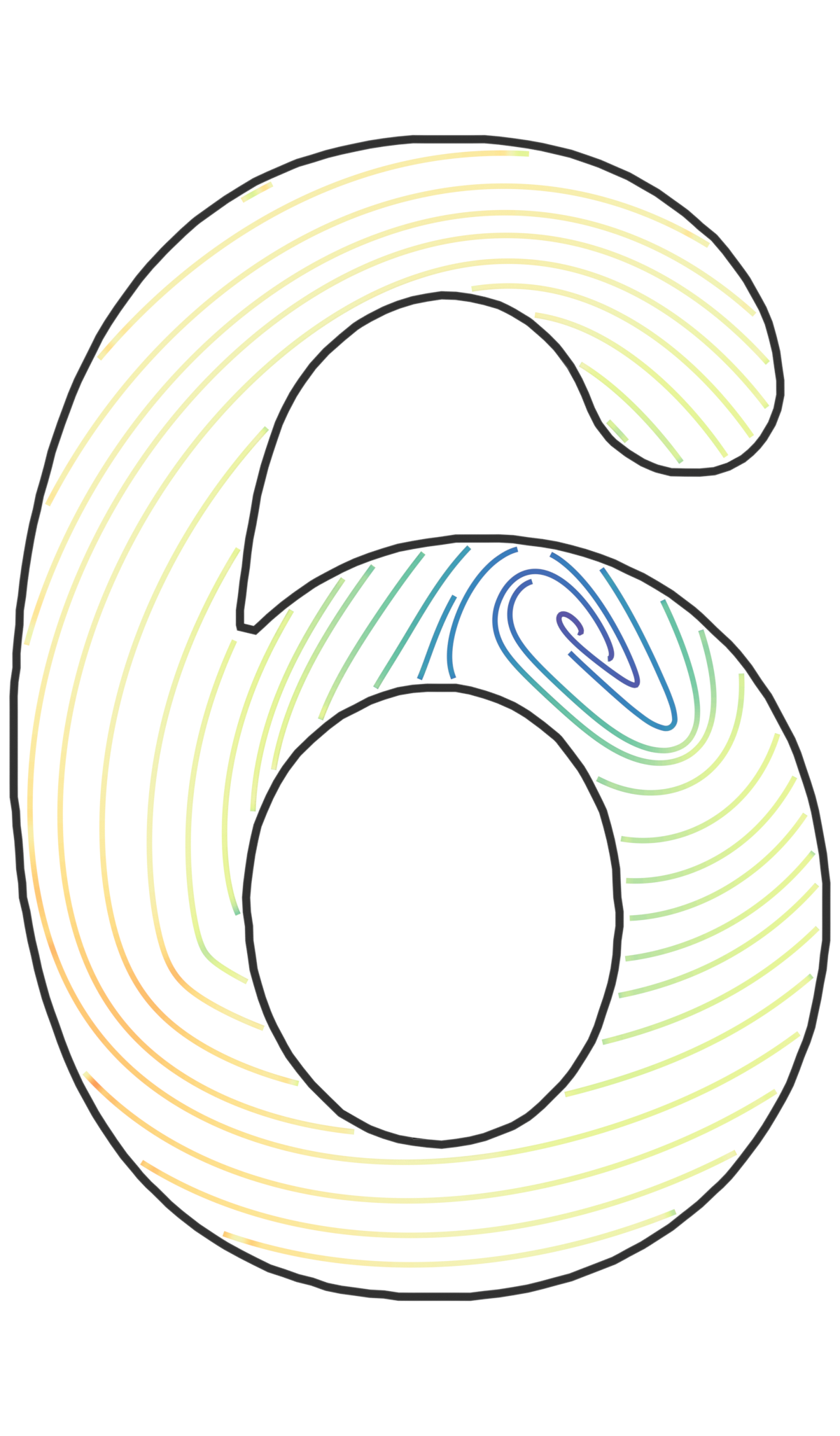}
 \raisebox{1cm}{\Large$=$}
	\includegraphics[trim=0cm 4cm 0cm 4cm, clip,width=0.12\textwidth]{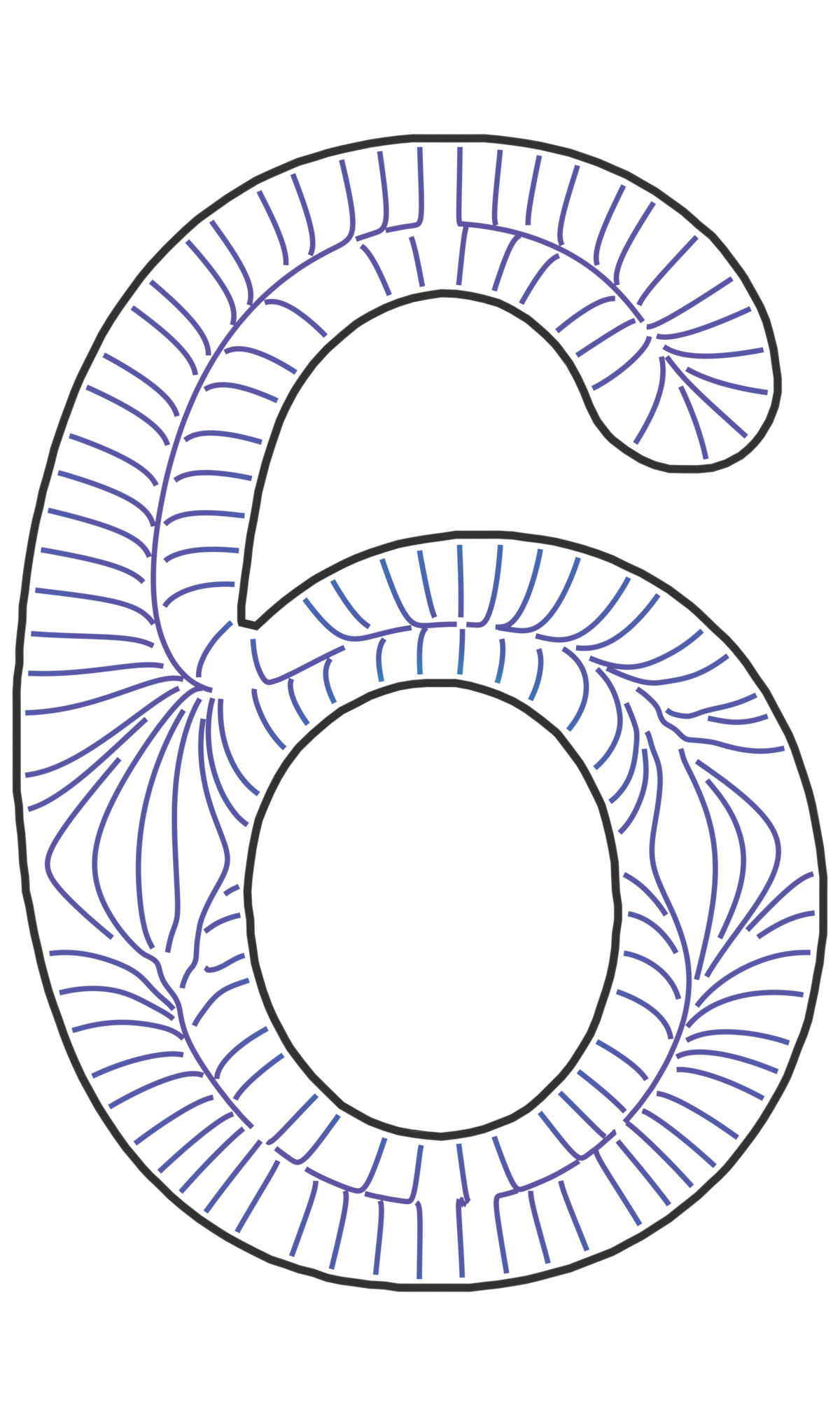}
  \raisebox{1cm}{\Large$+$}
	\includegraphics[trim=0cm 4cm 0cm 4cm, clip,width=0.12\textwidth]{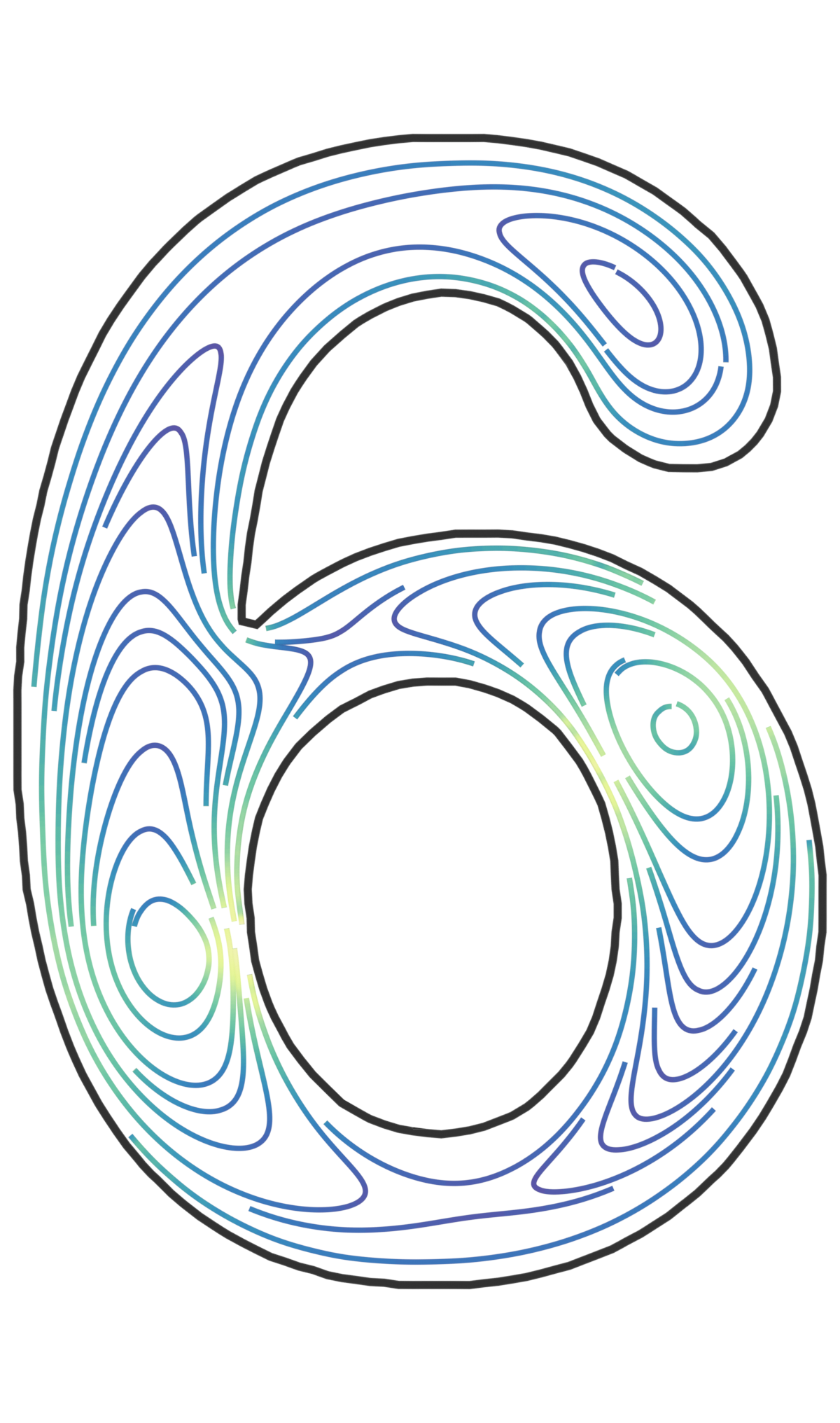}
  \raisebox{1cm}{\Large$+$}
	\includegraphics[trim=0cm 4cm 0cm 4cm, clip,width=0.12\textwidth]{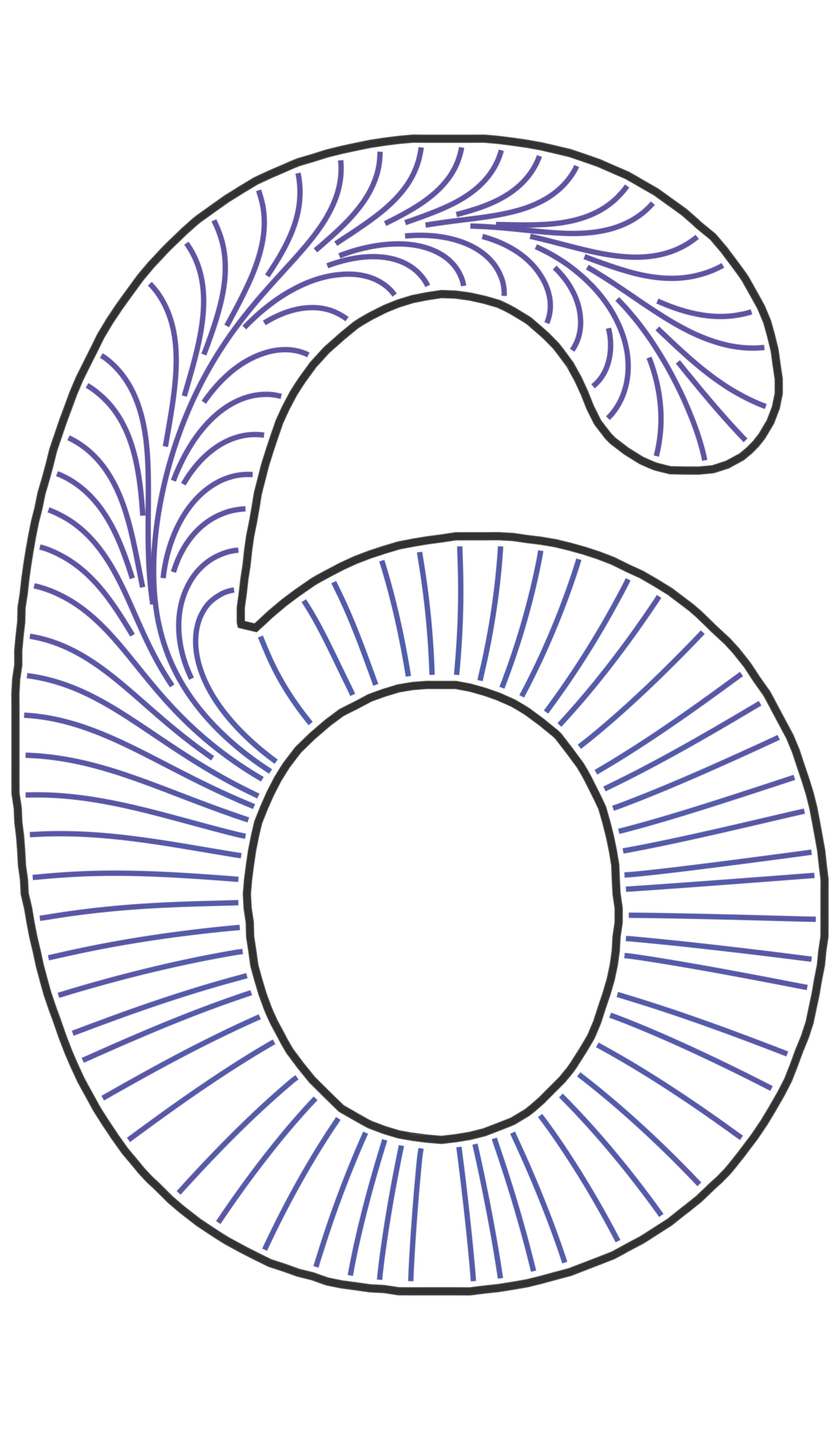}
  \raisebox{1cm}{\Large$+$}
	\includegraphics[trim=0cm 4cm 0cm 4cm, clip,width=0.12\textwidth]{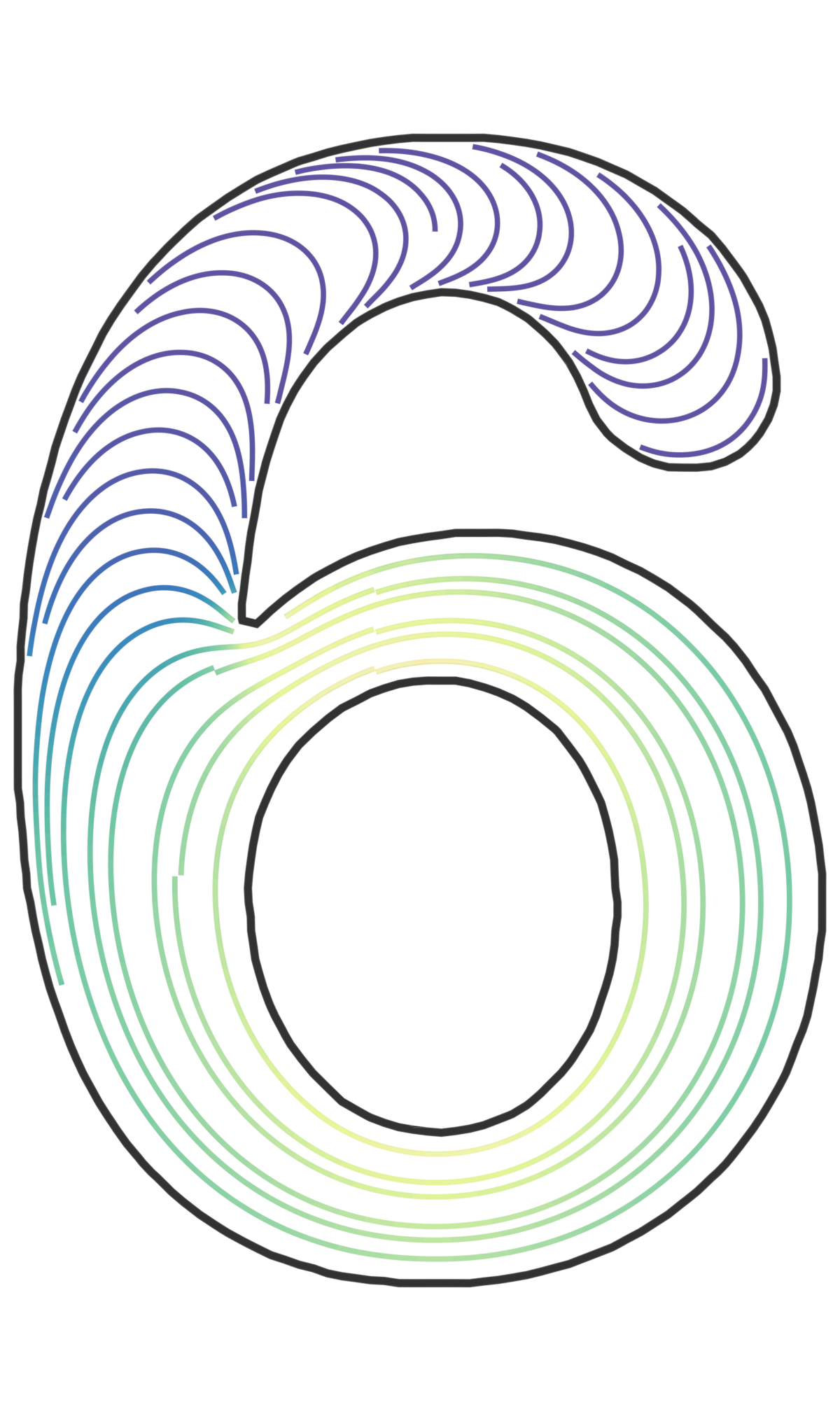}
  \raisebox{1cm}{\Large$+$}
	\includegraphics[trim=0cm 4cm 0cm 4cm, clip,width=0.12\textwidth]{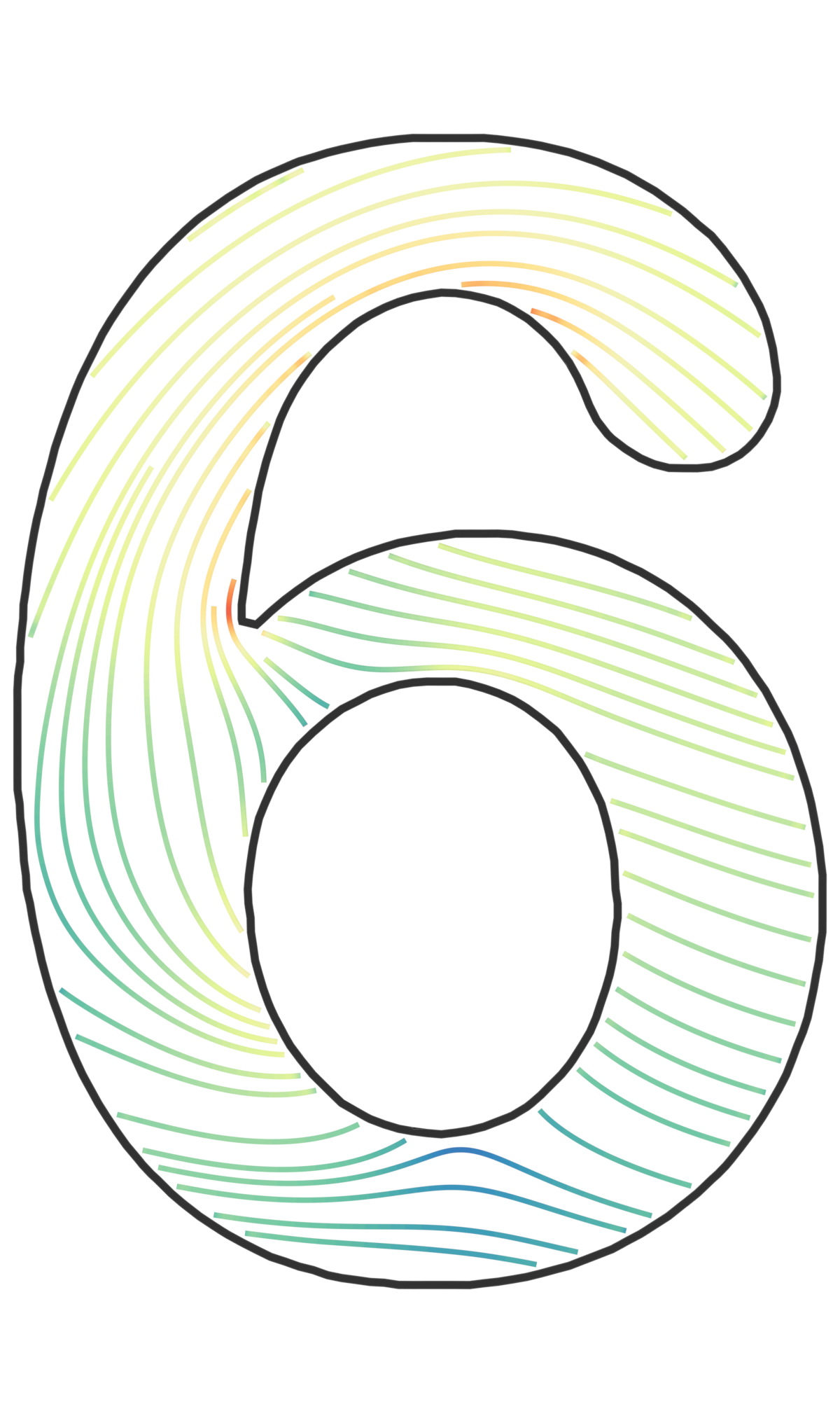}
	\caption{5-component Hodge decomposition on a 2D figure 6 domain. From left to right: the original vector field, the normal gradient field, the tangential curl field, the normal harmonic field, the tangential harmonic field, and the curly gradient field.}
	\label{fig.hd.2d.6.stepwise}
\end{figure}

\begin{figure}[ht]
	\centering
	\includegraphics[trim=7cm 4cm 7cm 5cm, clip, width=0.12\textwidth]{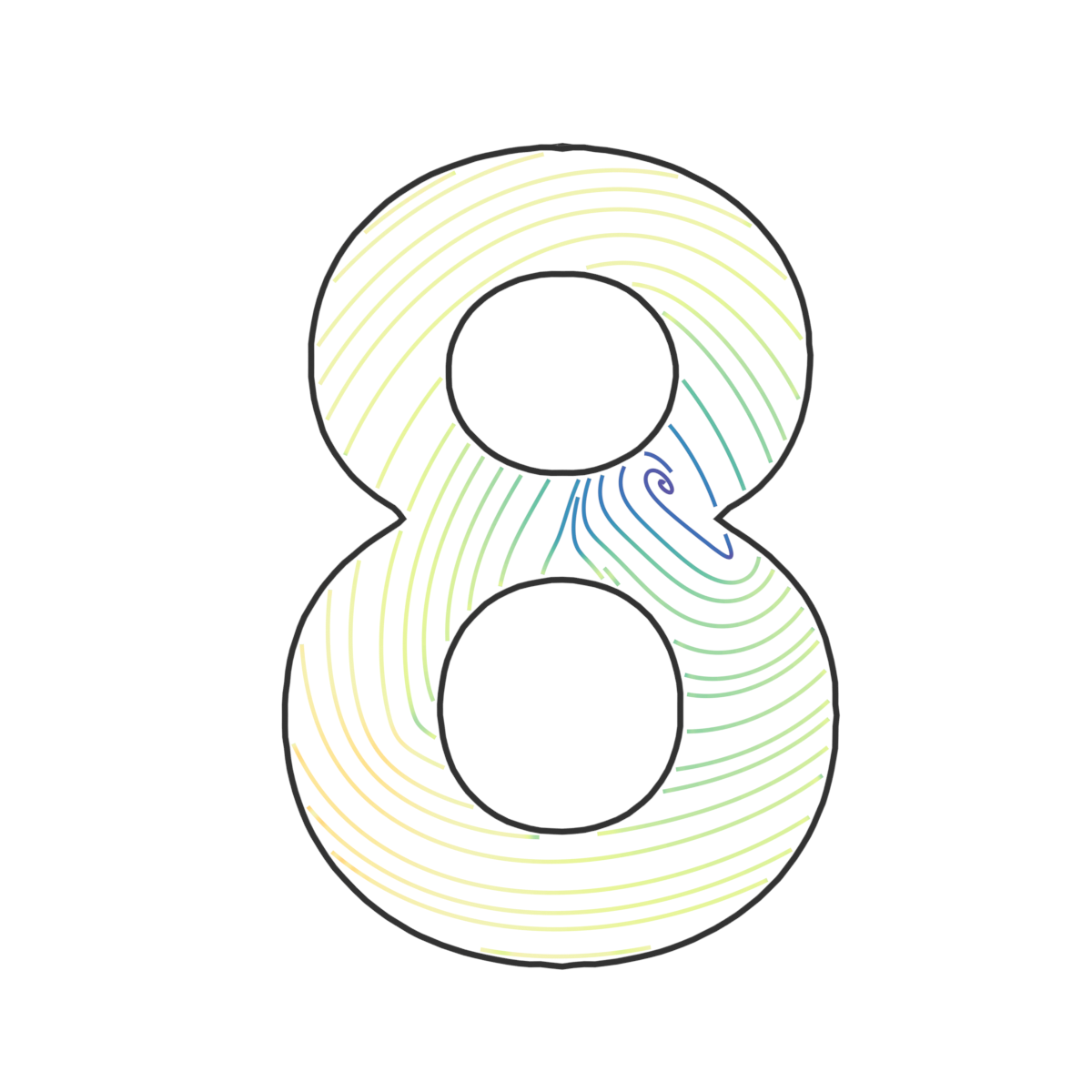}
 \raisebox{1cm}{\Large$=$}
	\includegraphics[trim=7cm 4cm 7cm 5cm, clip,width=0.12\textwidth]{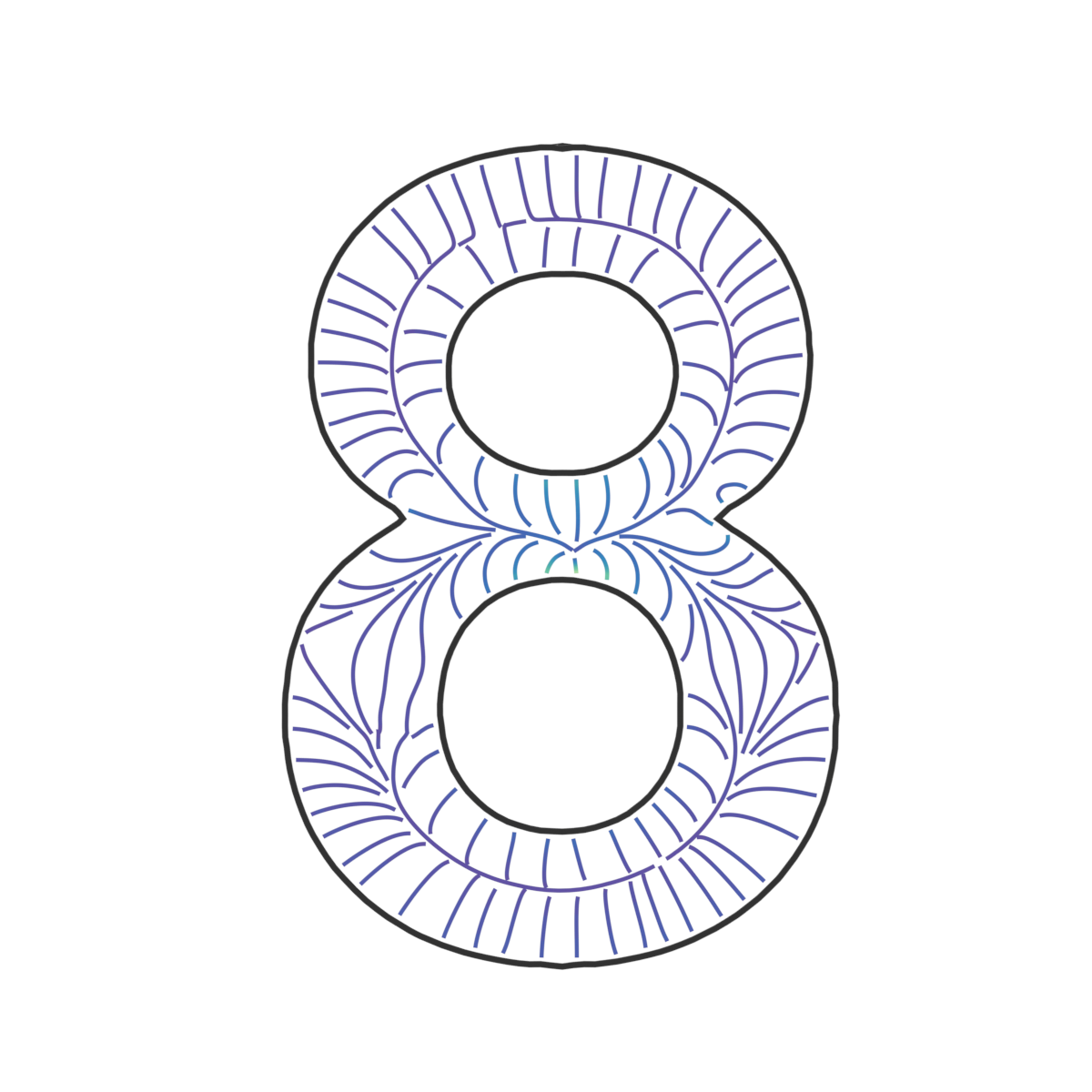}
  \raisebox{1cm}{\Large$+$}
	\includegraphics[trim=7cm 4cm 7cm 5cm, clip,width=0.12\textwidth]{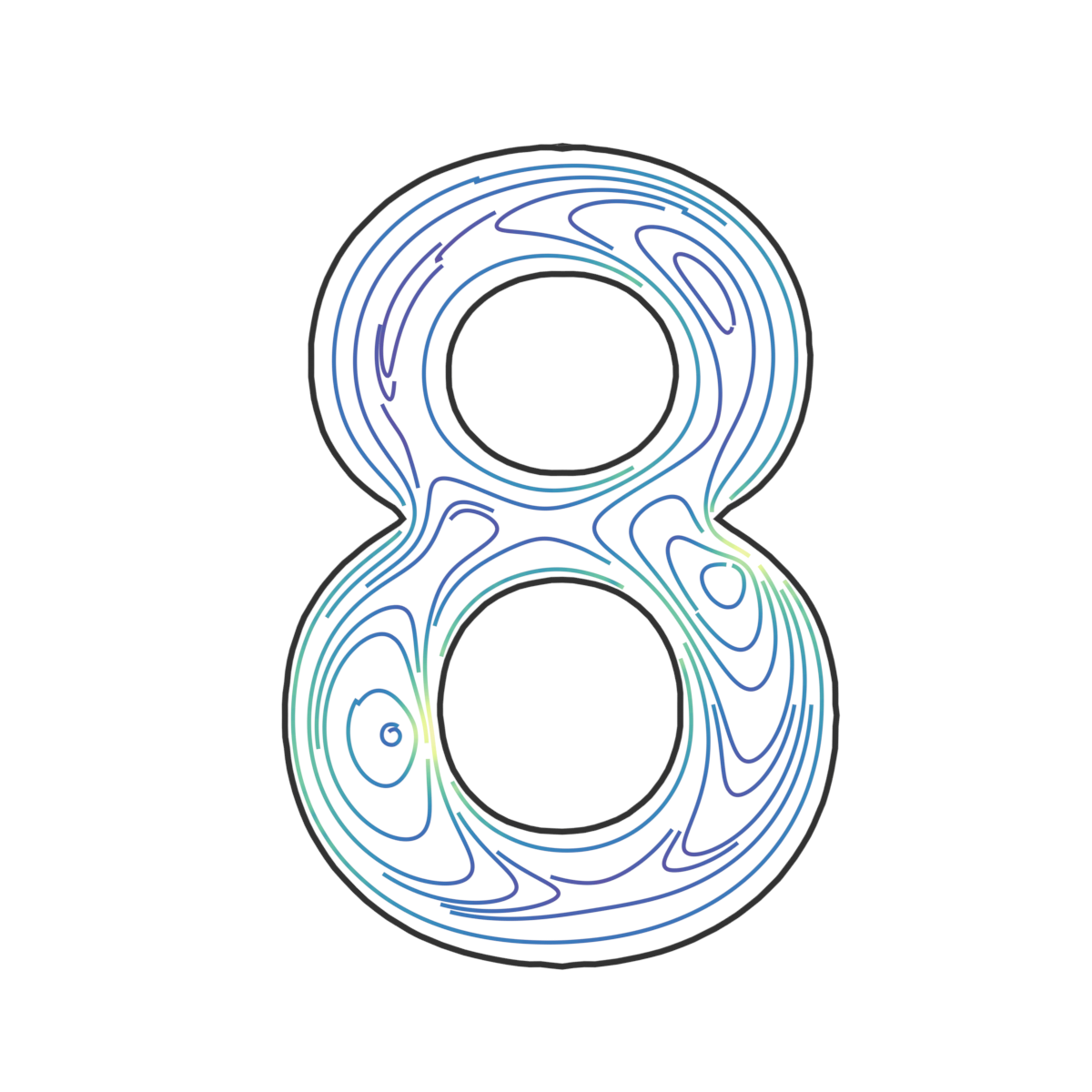}
  \raisebox{1cm}{\Large$+$}
	\includegraphics[trim=7cm 4cm 7cm 5cm, clip,width=0.12\textwidth]{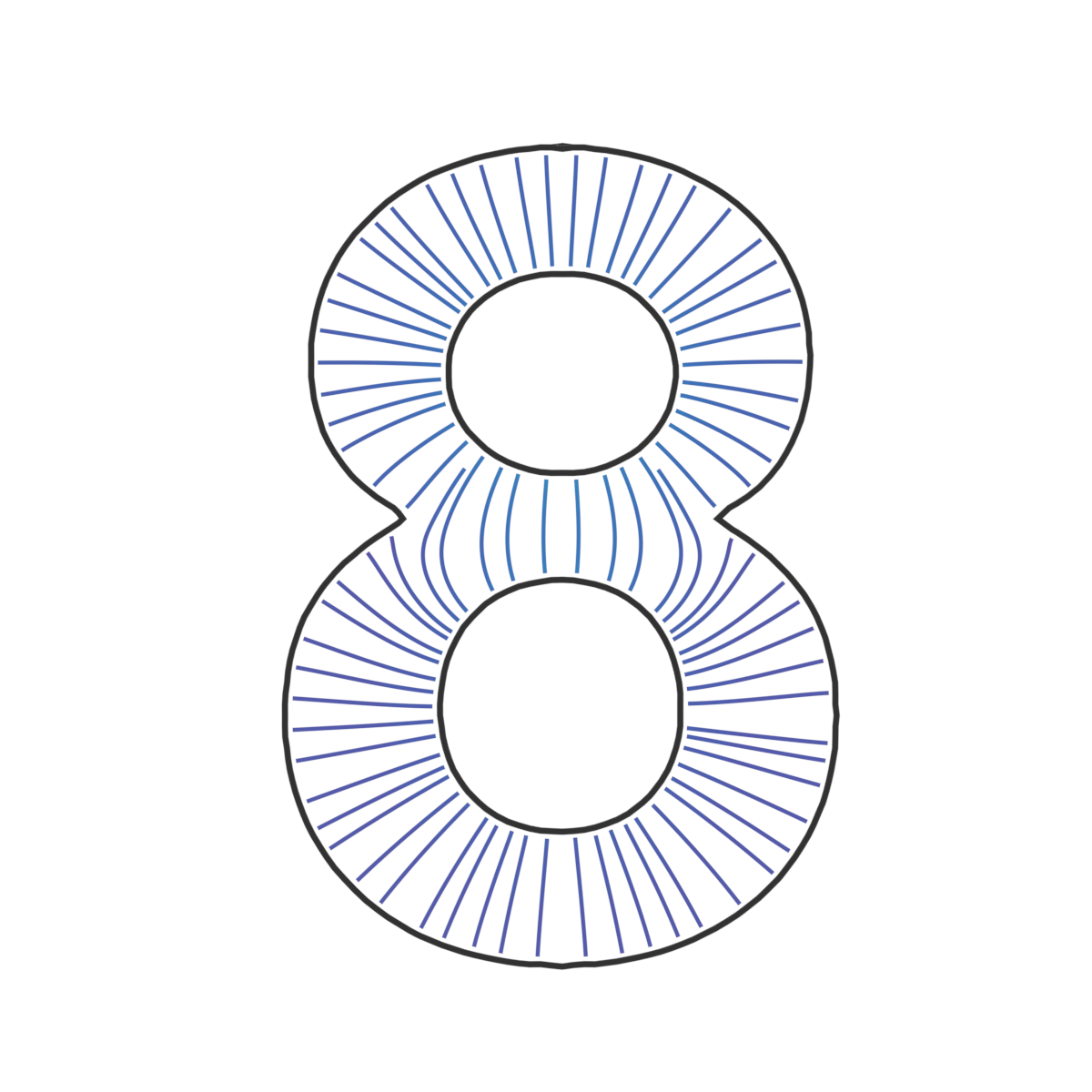}
  \raisebox{1cm}{\Large$+$}
	\includegraphics[trim=7cm 4cm 7cm 5cm, clip,width=0.12\textwidth]{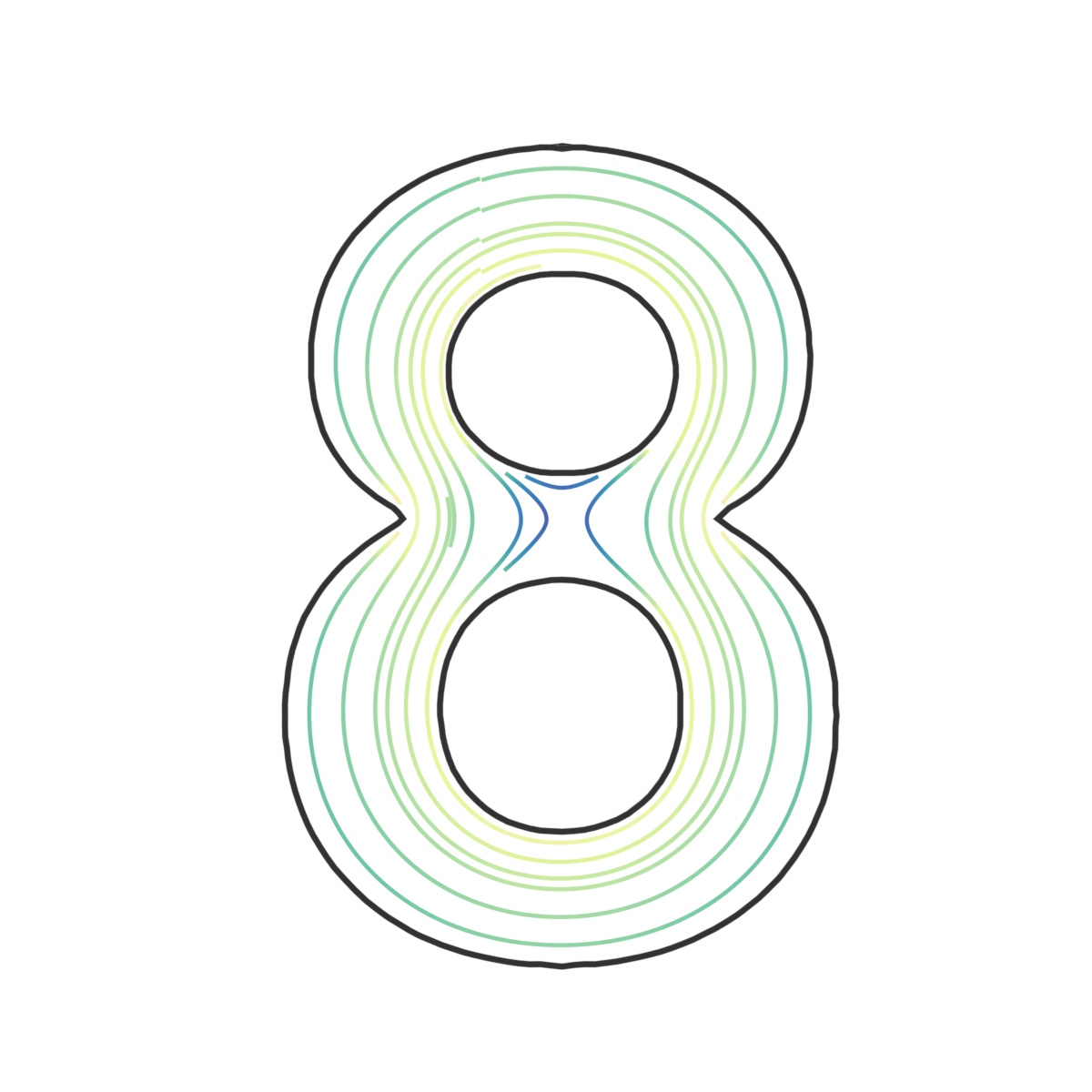}
  \raisebox{1cm}{\Large$+$}
	\includegraphics[trim=7cm 4cm 7cm 5cm, clip,width=0.12\textwidth]{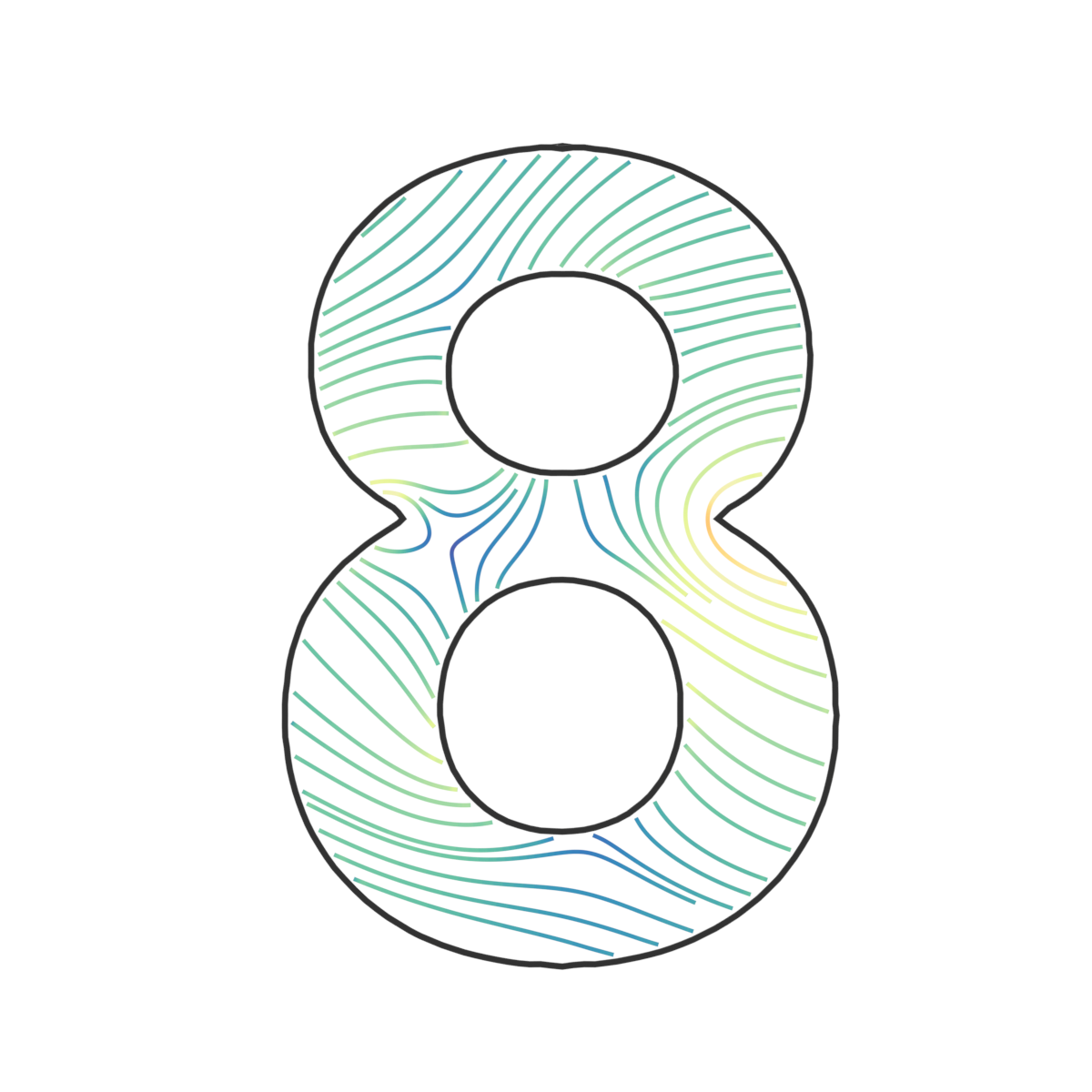}
	\caption{5-component Hodge decomposition on a 2D figure 8 domain. From left to right: the original vector field, the normal gradient field, the tangential curl field, the normal harmonic field, the tangential harmonic field, and the curly gradient field.}
	\label{fig.hd.2d.8.stepwise}
\end{figure}

In Fig.~\ref{fig.hd.2d.bunny.stepwise}, the 5-component Hodge decomposition is computed in our approach on a 2D bunny-shaped domain with one hole. Both the normal harmonic space $\mathcal{H}_n$ and the tangential harmonic space $\mathcal{H}_t$ are one-dimensional, since $\beta_1=1$ due to the annulus topology. The 5 components are pairwise orthogonal, with an \l-inner product of 0 up to the precision of the linear solver used. Fig.~\ref{fig.hd.2d.6.stepwise} and Fig.~\ref{fig.hd.2d.8.stepwise} present two additional examples of the $5$-component Hodge decomposition in the 2D case, with vector fields defined on a 2D figure-6 model and a 2D figure-8 model, respectively. The dimension of $\mathcal{H}_n$ and $\mathcal{H}_t$ is determined by the topology of each domain, given as one for the figure-6 model, and two for the figure-8 model, respectively.

\begin{figure*}[t]
	\centering
	\includegraphics[trim=2cm 2cm 2cm 1cm, clip, width=0.12\textwidth]{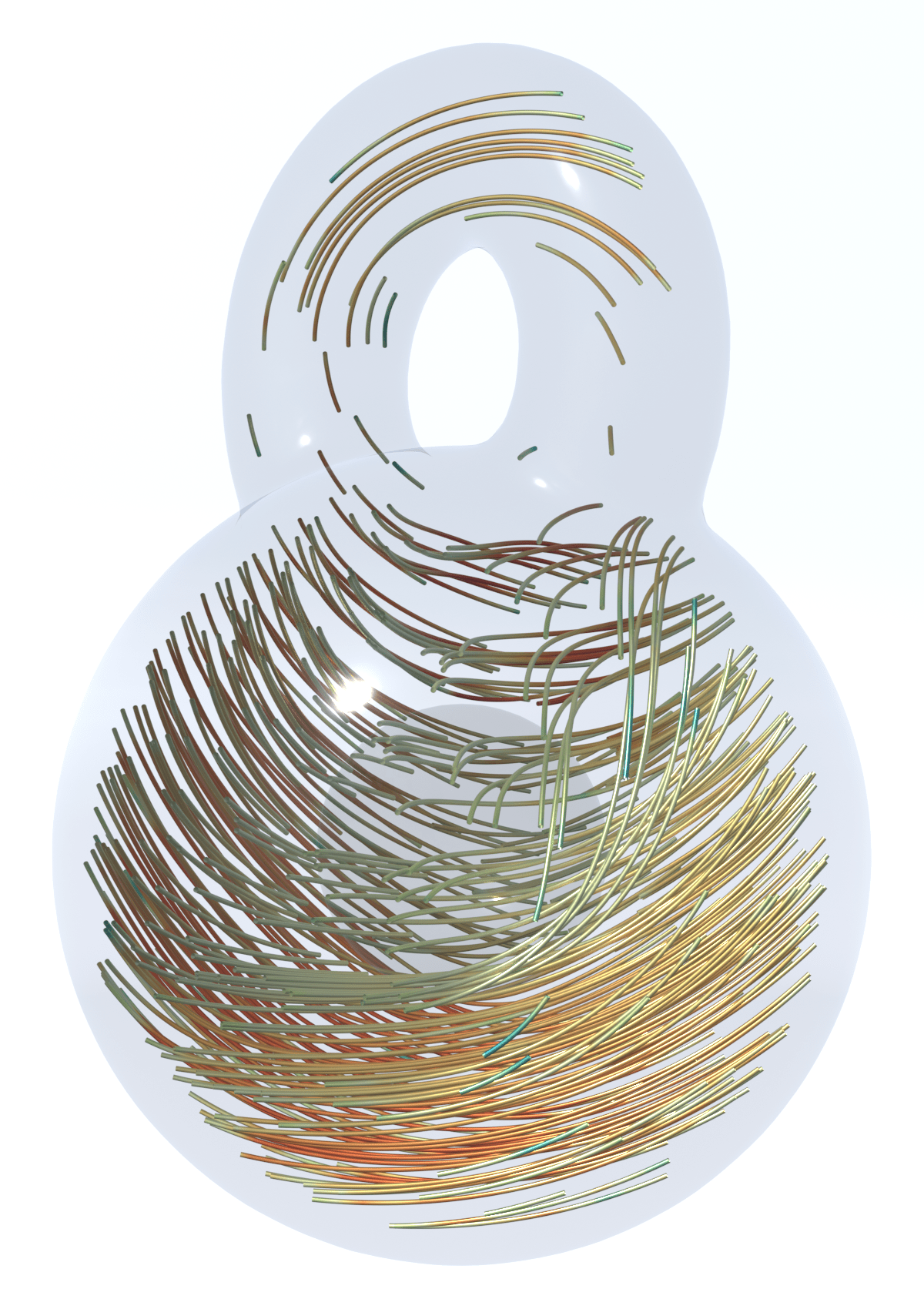}
 \raisebox{1cm}{\Large$=$}
	\includegraphics[trim=2cm 2cm 2cm 1cm, clip, width=0.12\textwidth]{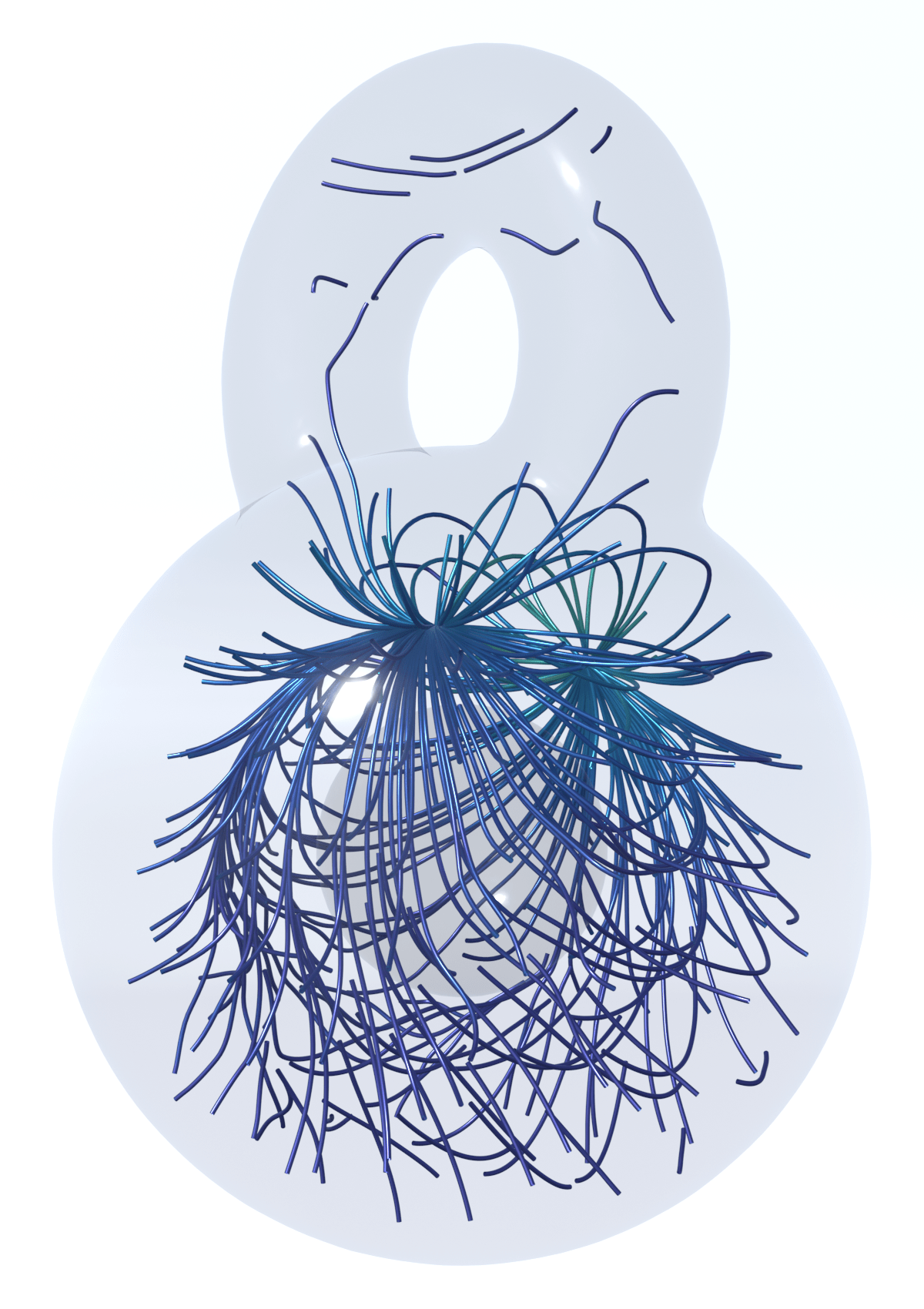}
  \raisebox{1cm}{\Large$+$}
	\includegraphics[trim=2cm 2cm 2cm 1cm, clip, width=0.12\textwidth]{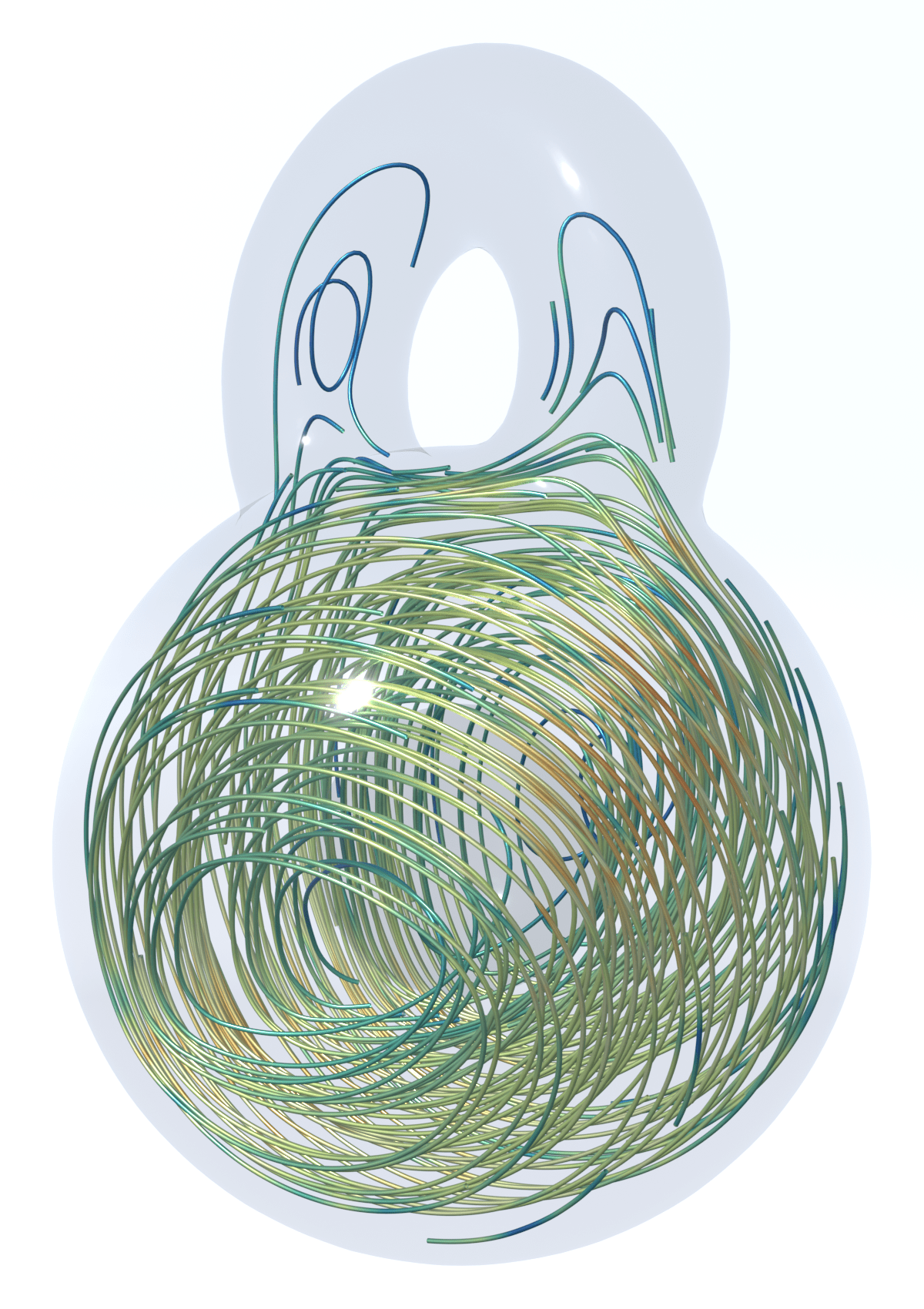}
  \raisebox{1cm}{\Large$+$}
	\includegraphics[trim=2cm 2cm 2cm 1cm, clip, width=0.12\textwidth]{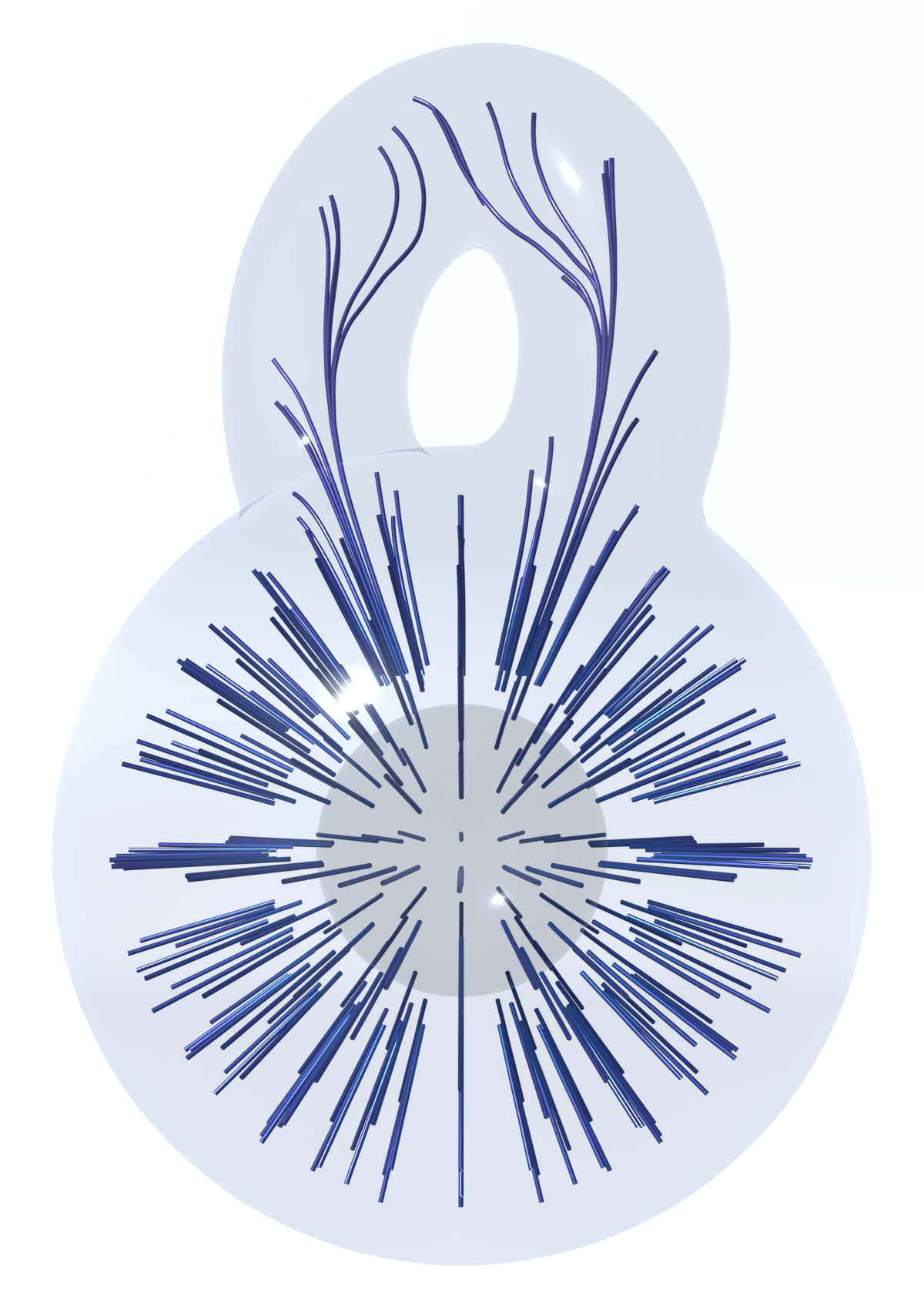}
  \raisebox{1cm}{\Large$+$}
	\includegraphics[trim=2cm 2cm 2cm 1cm, clip, width=0.12\textwidth]{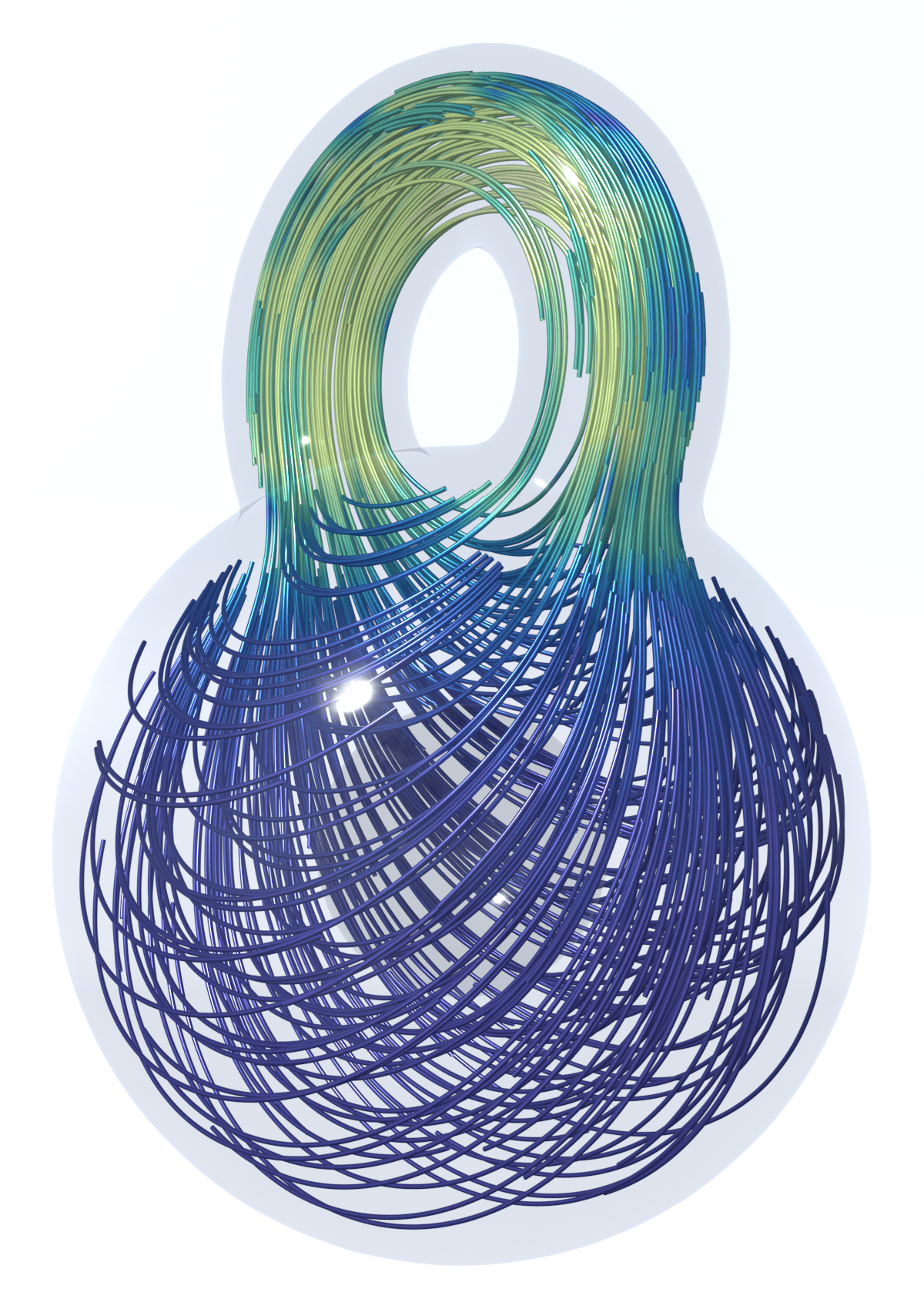}
  \raisebox{1cm}{\Large$+$}
	\includegraphics[trim=2cm 2cm 2cm 1cm, clip, width=0.12\textwidth]{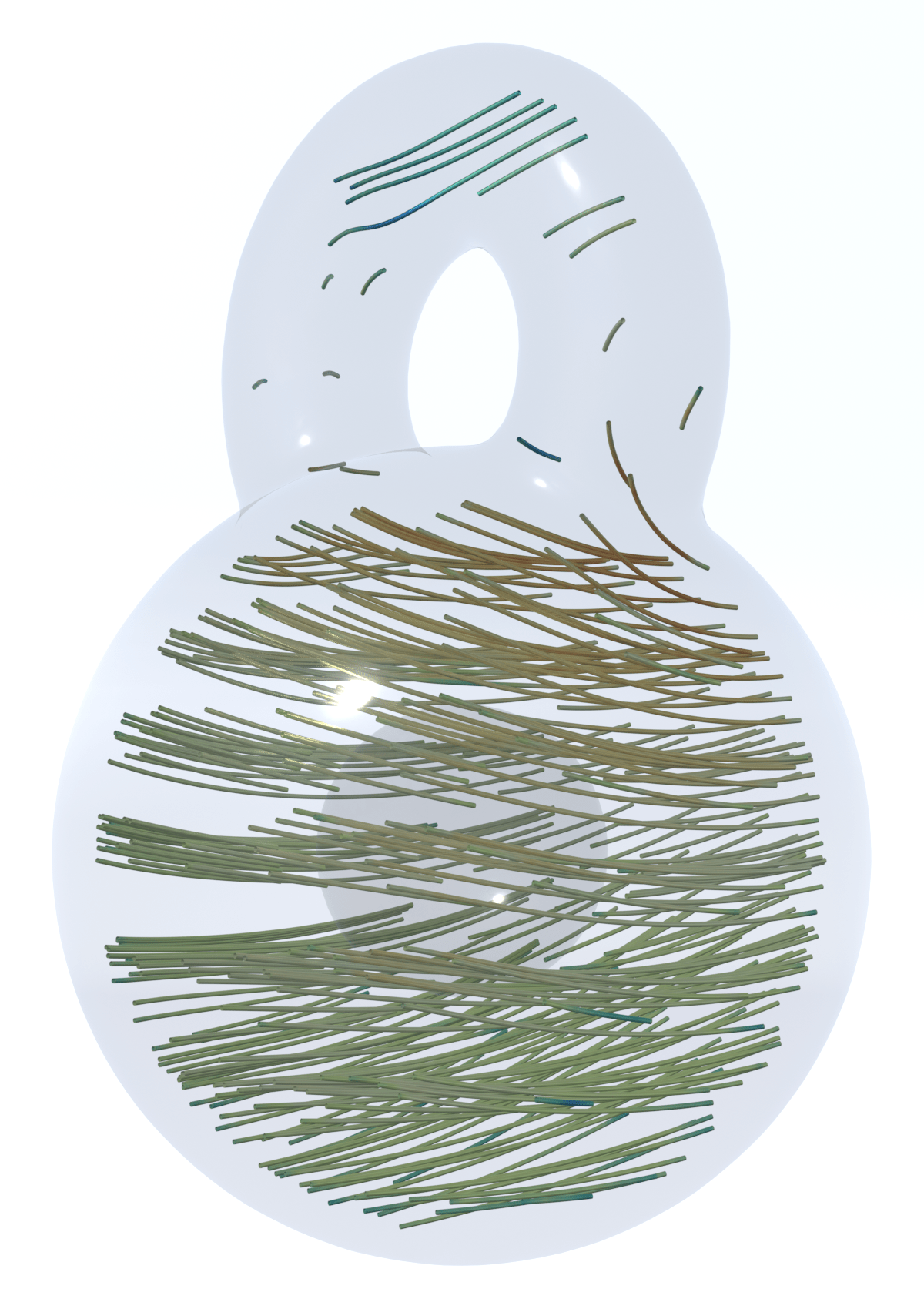}
	\caption{3D Hodge decomposition on kettlebell shape. From left to right: the original vector field, the normal gradient field, the tangential curl field, the normal harmonic field, the tangential harmonic field, and the curly gradient field.}
    \label{fig.hd.3d.kettlebell}
\end{figure*}

\begin{figure*}[t]
	\centering
	\includegraphics[trim=2cm 2cm 2cm 1cm, clip, width=0.12\textwidth]{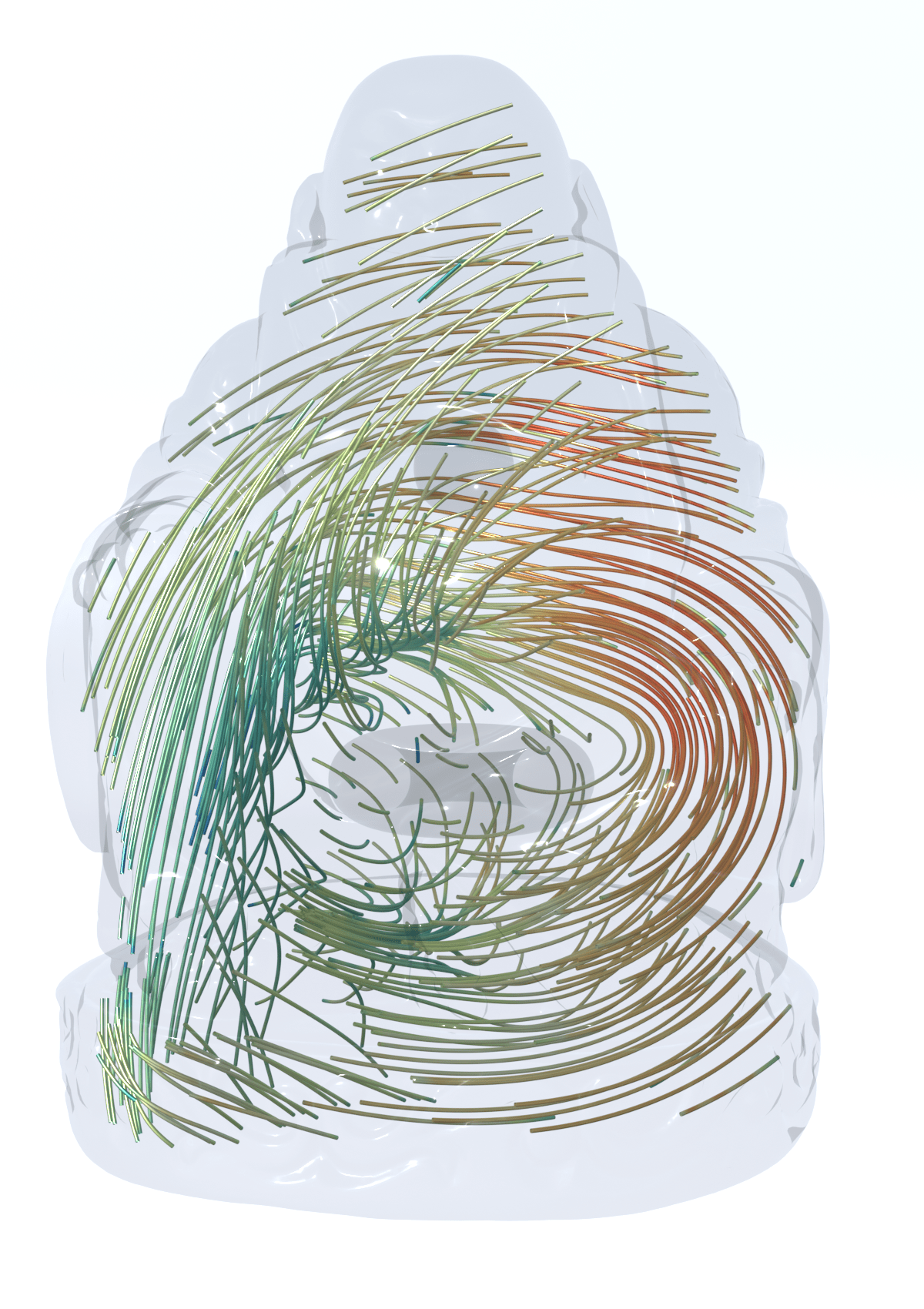}
 \raisebox{1cm}{\Large$=$}
	\includegraphics[trim=2cm 2cm 2cm 1cm, clip, width=0.12\textwidth]{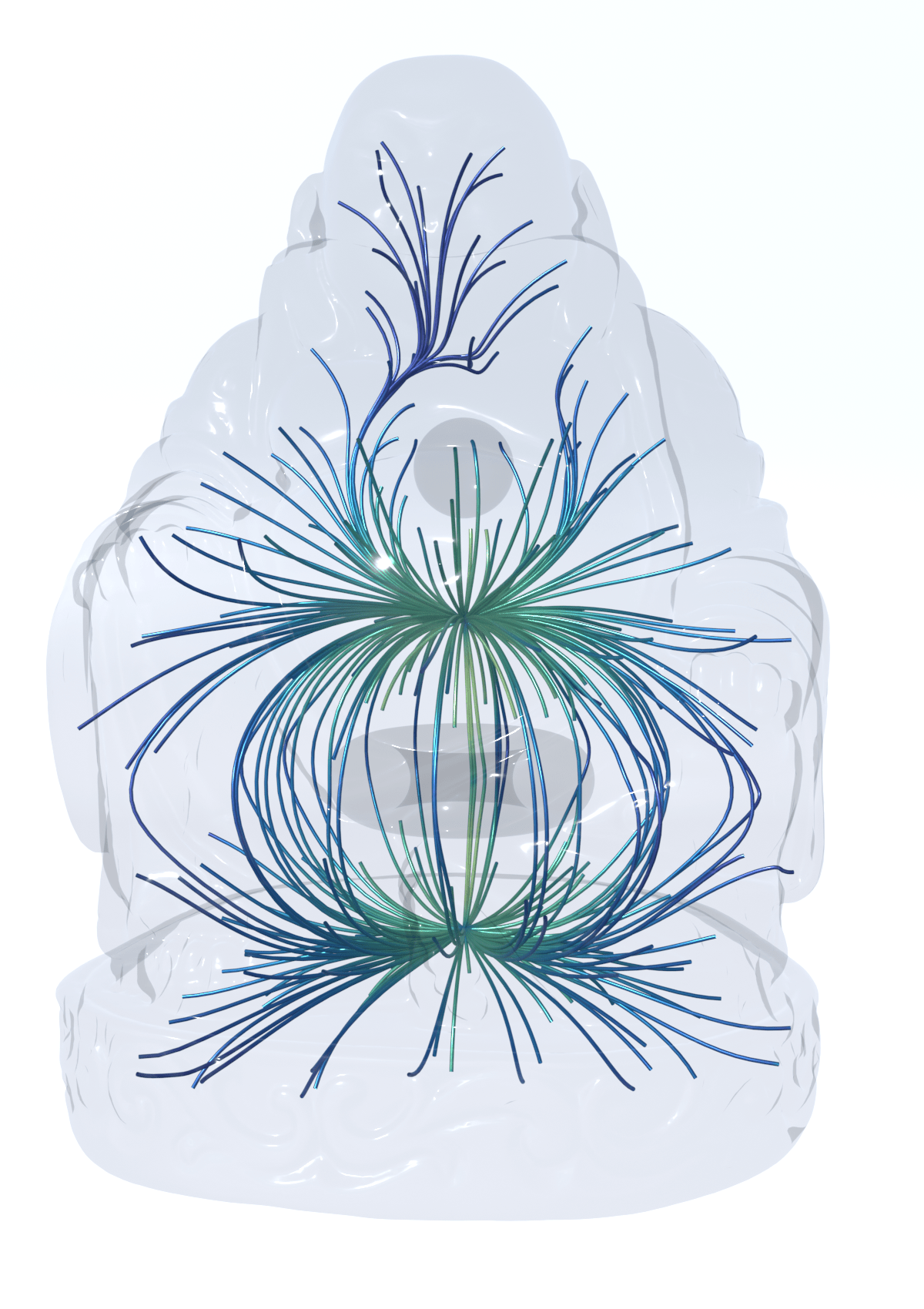}
  \raisebox{1cm}{\Large$+$}
	\includegraphics[trim=2cm 2cm 2cm 1cm, clip, width=0.12\textwidth]{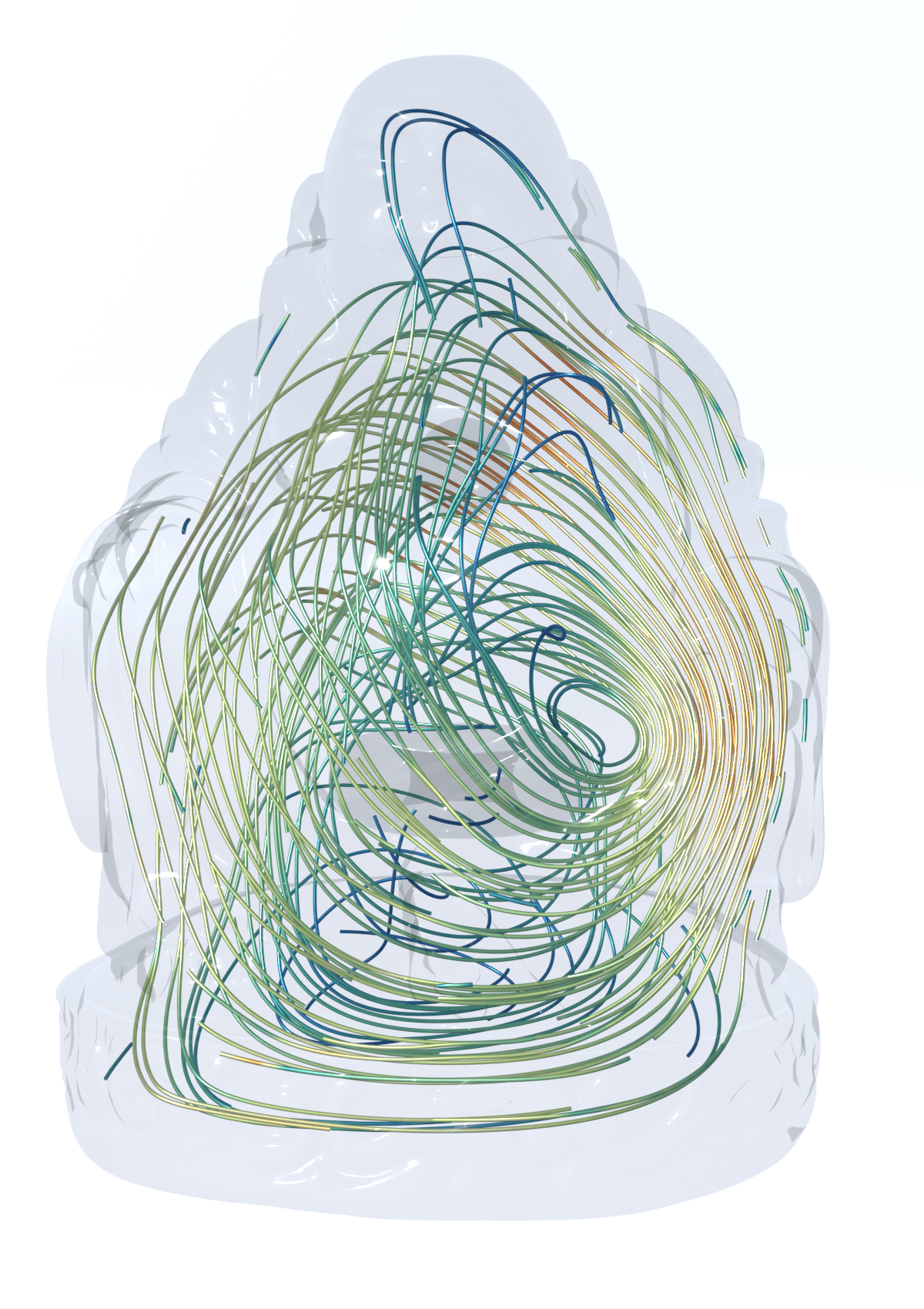}
  \raisebox{1cm}{\Large$+$}
	\includegraphics[trim=2cm 2cm 2cm 1cm, clip, width=0.12\textwidth]{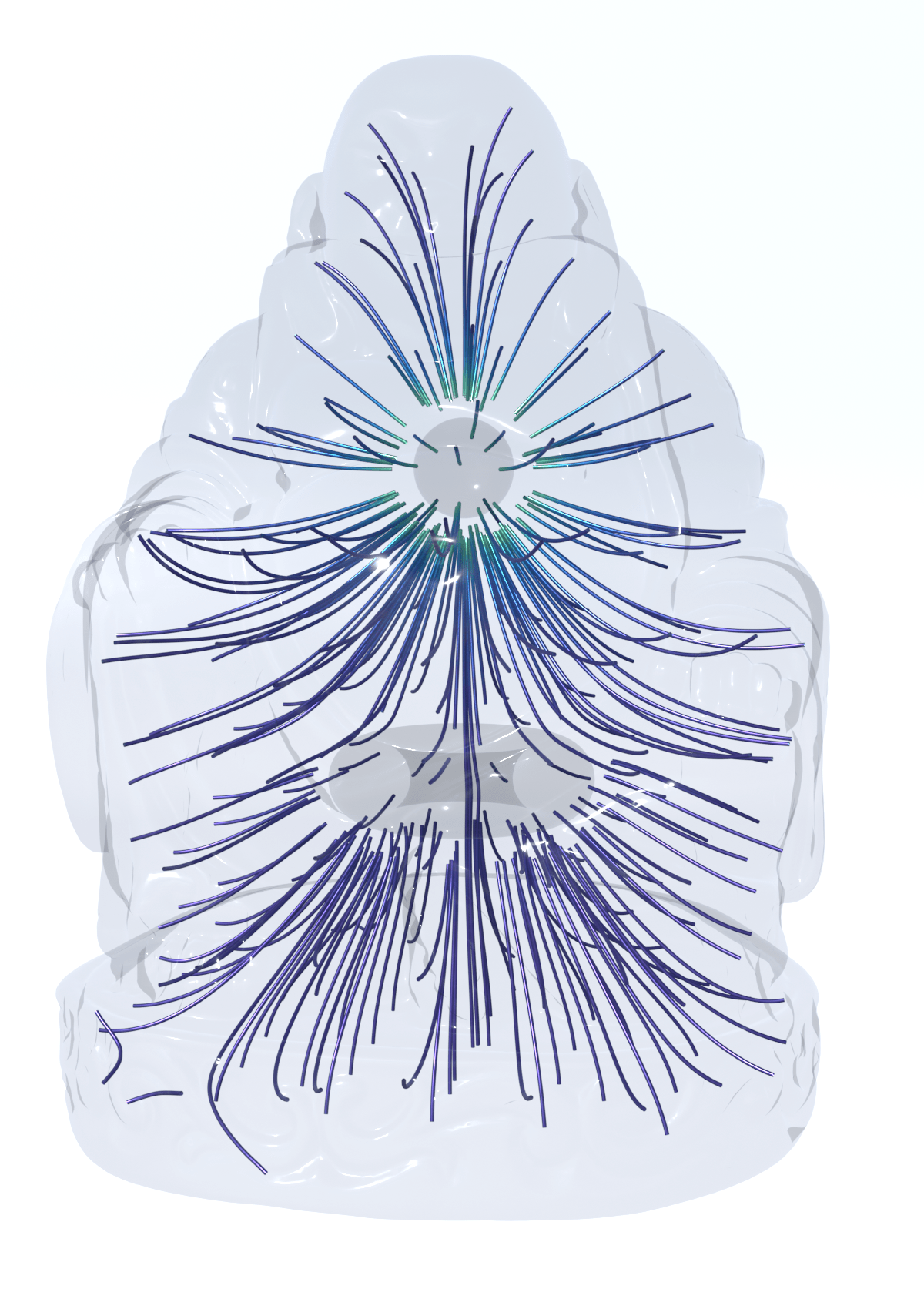}
  \raisebox{1cm}{\Large$+$}
	\includegraphics[trim=2cm 2cm 2cm 1cm, clip, width=0.12\textwidth]{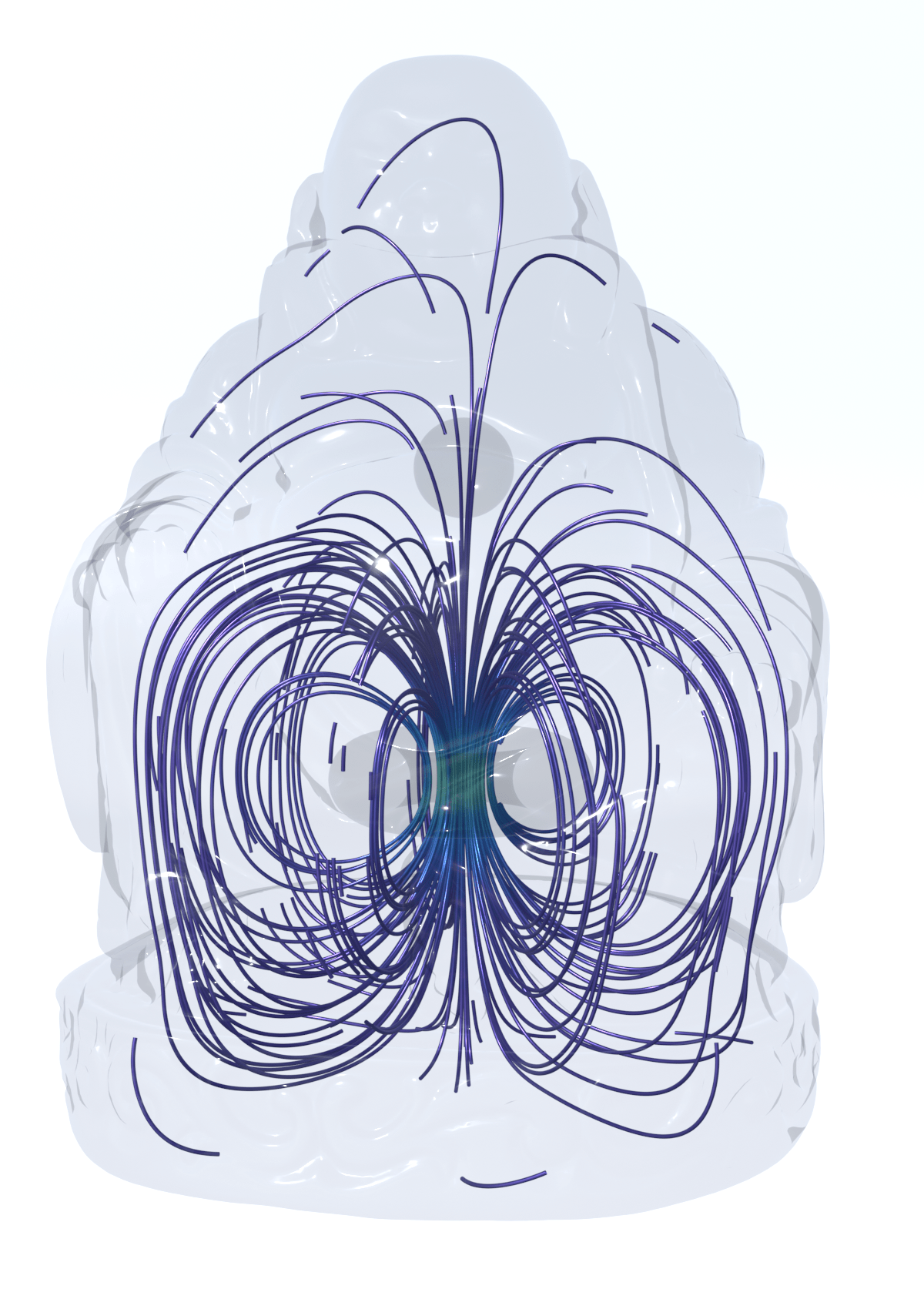}
  \raisebox{1cm}{\Large$+$}
	\includegraphics[trim=2cm 2cm 2cm 1cm, clip, width=0.12\textwidth]{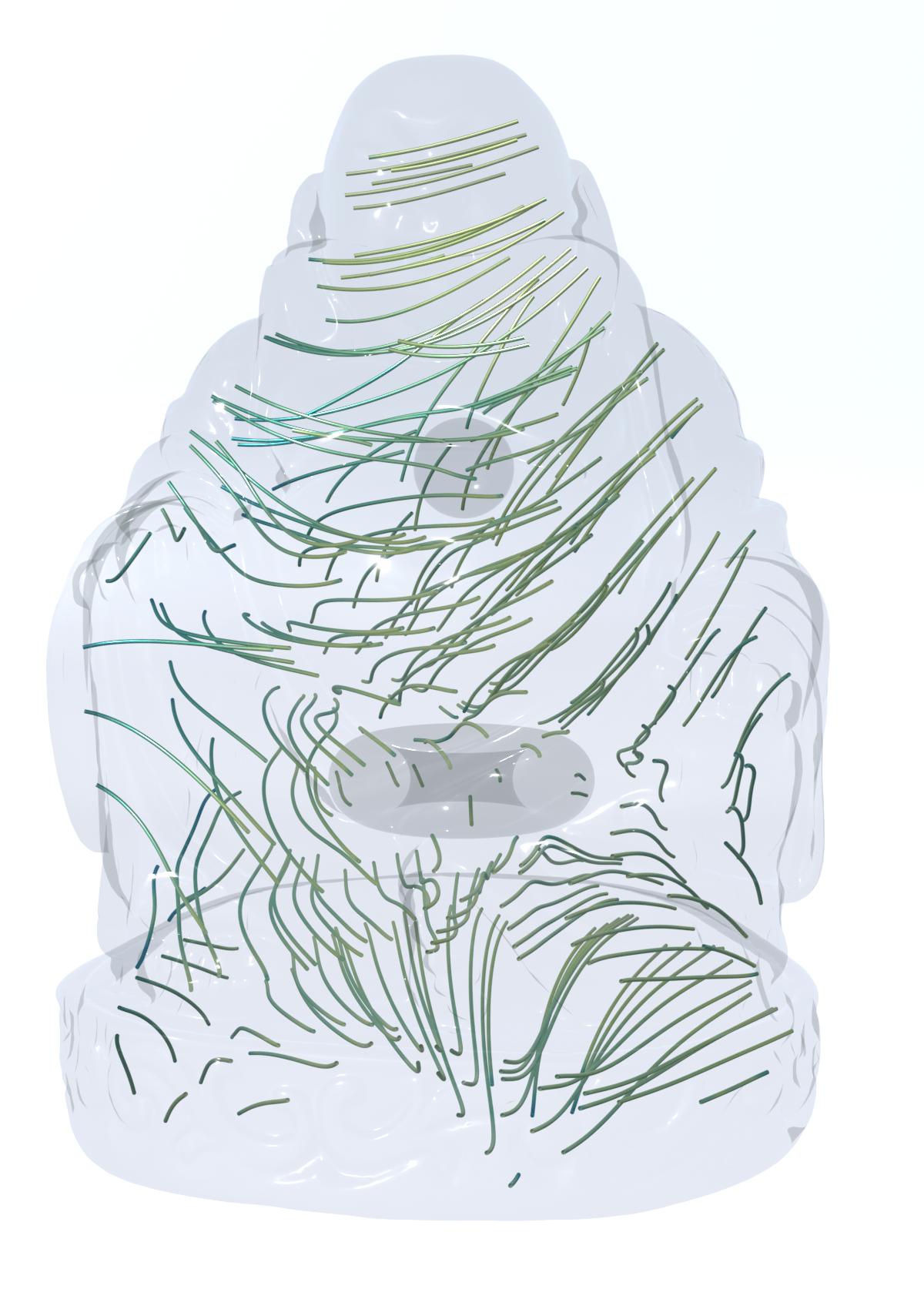}
	\caption{3D Hodge decomposition on Buddha model. To create variations in the topology, a ball and a torus are cut from the inside.}
    \label{fig.buddha}
\end{figure*}

\begin{figure*}[t]
	\centering
	\includegraphics[trim=2cm 1cm 2cm 1cm, clip, width=0.15\textwidth]{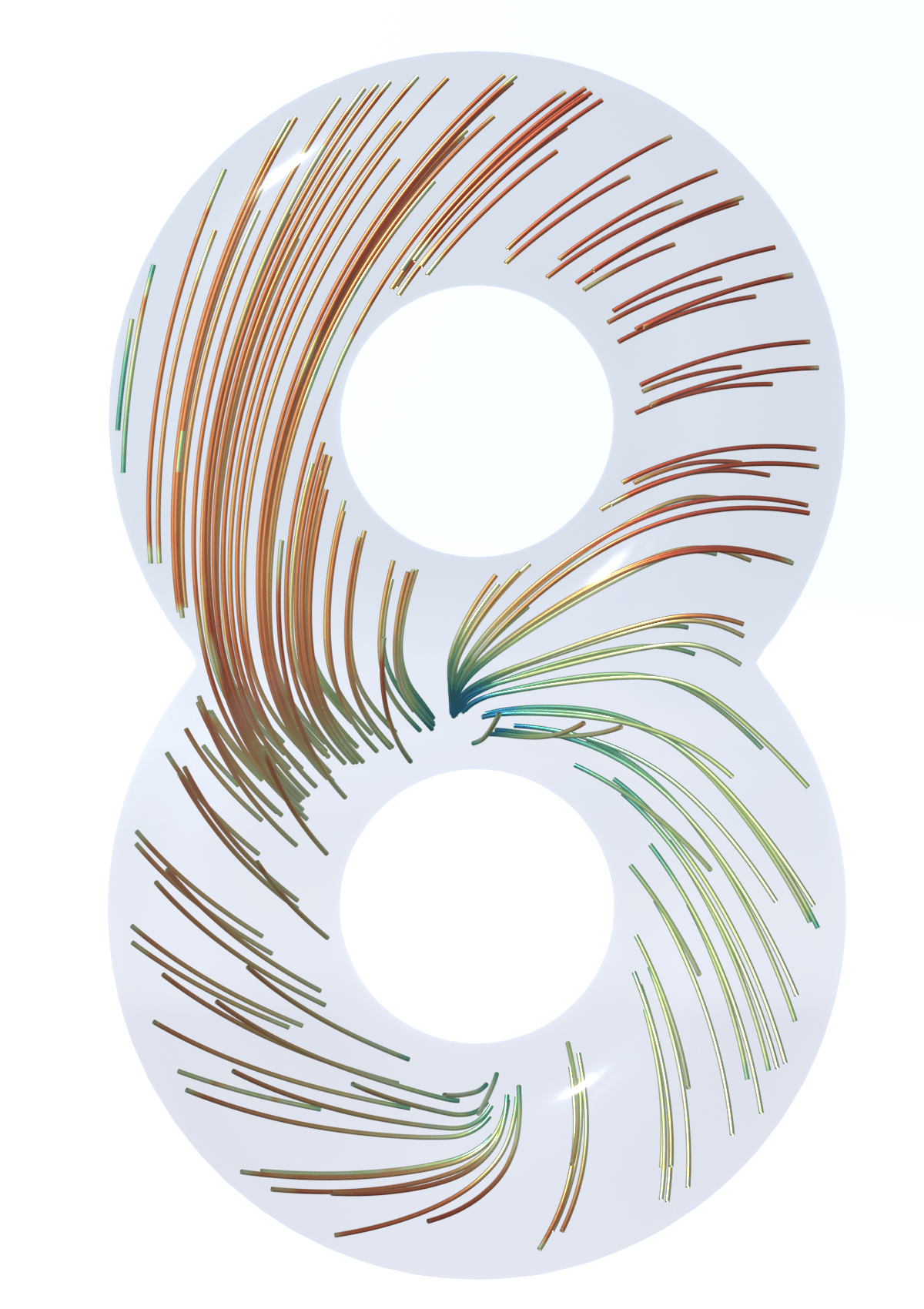}
 \raisebox{1cm}{\Large$=$}
	\includegraphics[trim=2cm 1cm 2cm 1cm, clip, width=0.15\textwidth]{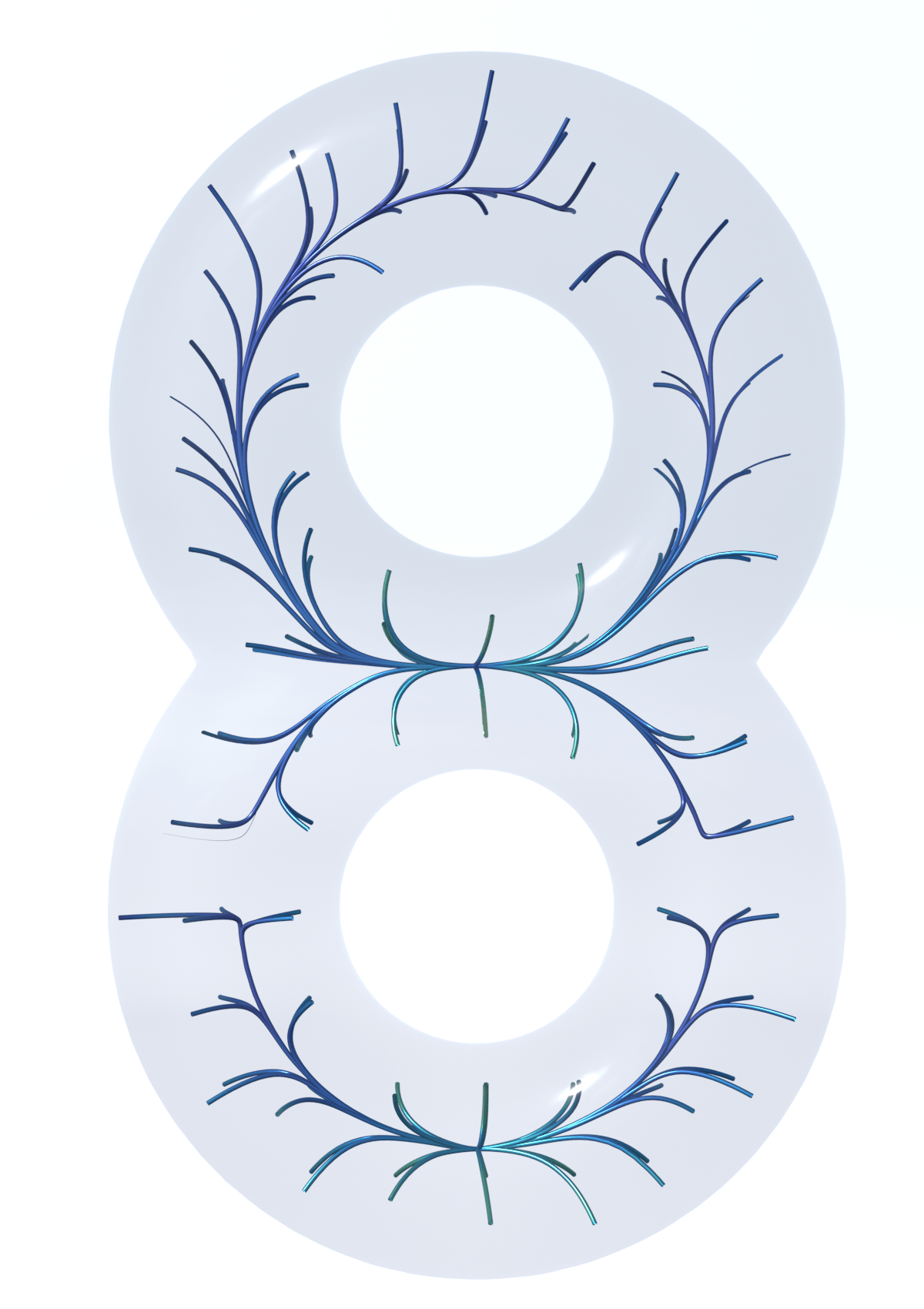}
  \raisebox{1cm}{\Large$+$}
	\includegraphics[trim=2cm 1cm 2cm 1cm, clip, width=0.15\textwidth]{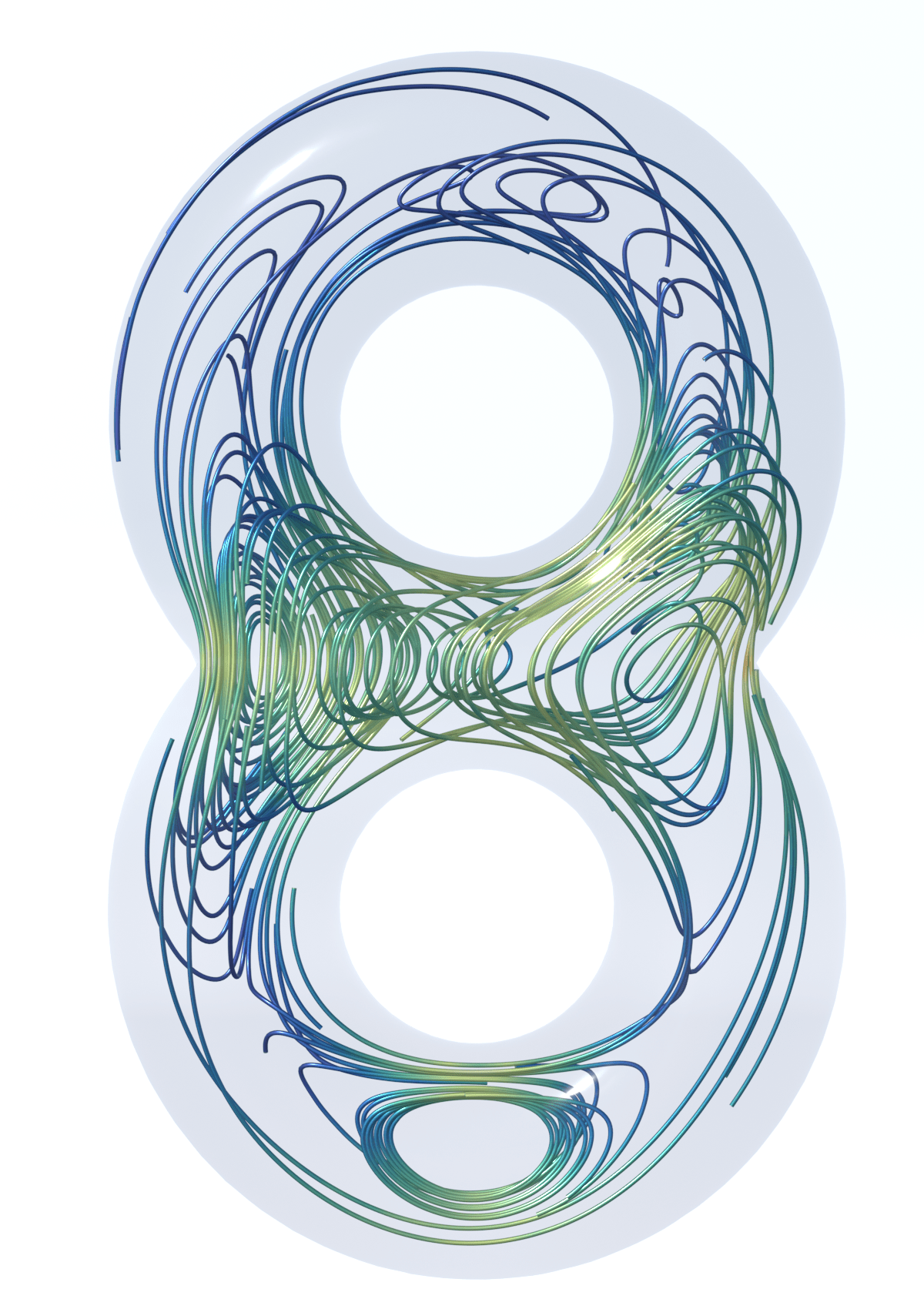}
  \raisebox{1cm}{\Large$+$}
	\includegraphics[trim=2cm 1cm 2cm 1cm, clip, width=0.15\textwidth]{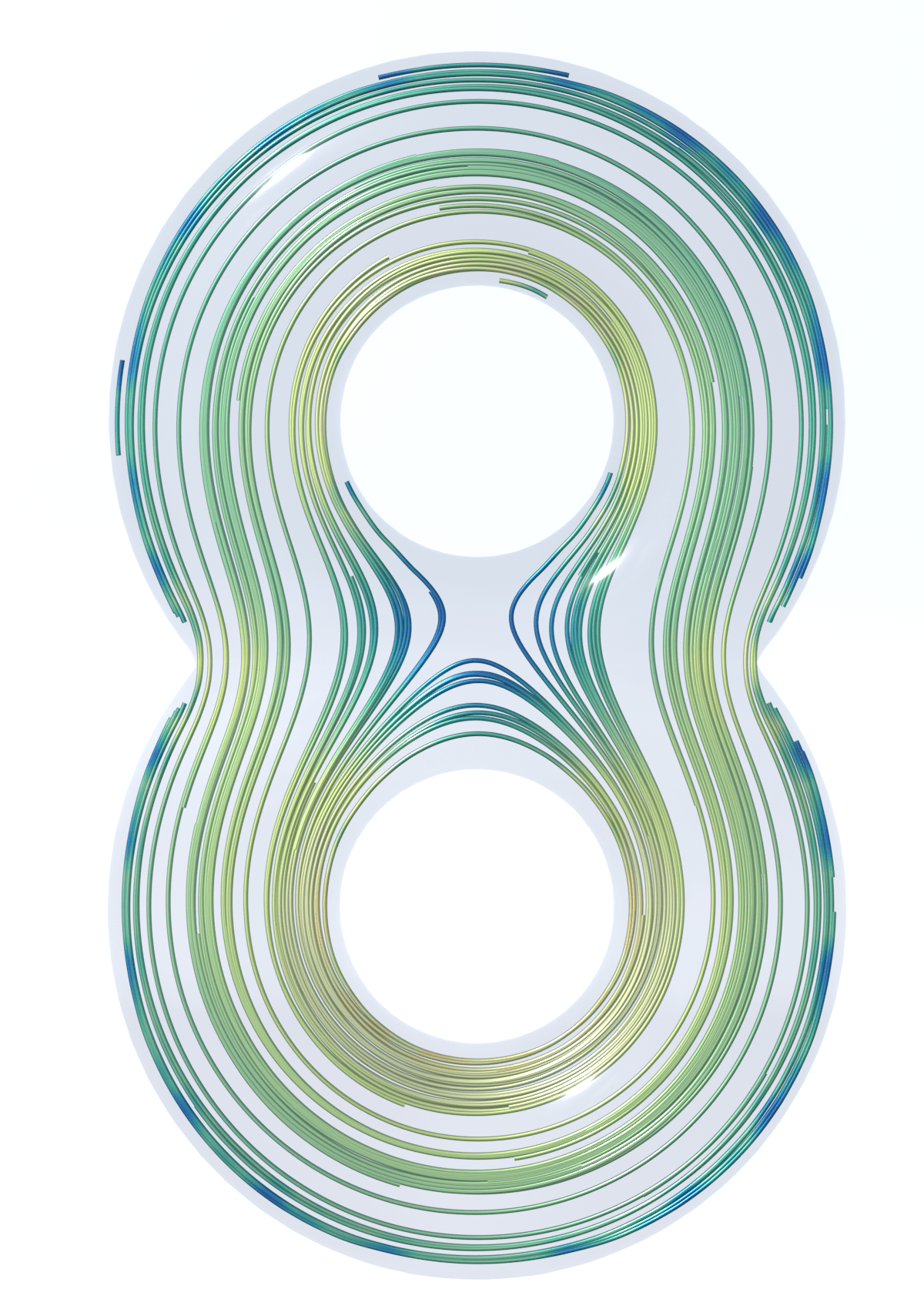}
  \raisebox{1cm}{\Large$+$}
	\includegraphics[trim=2cm 1cm 2cm 1cm, clip, width=0.15\textwidth]{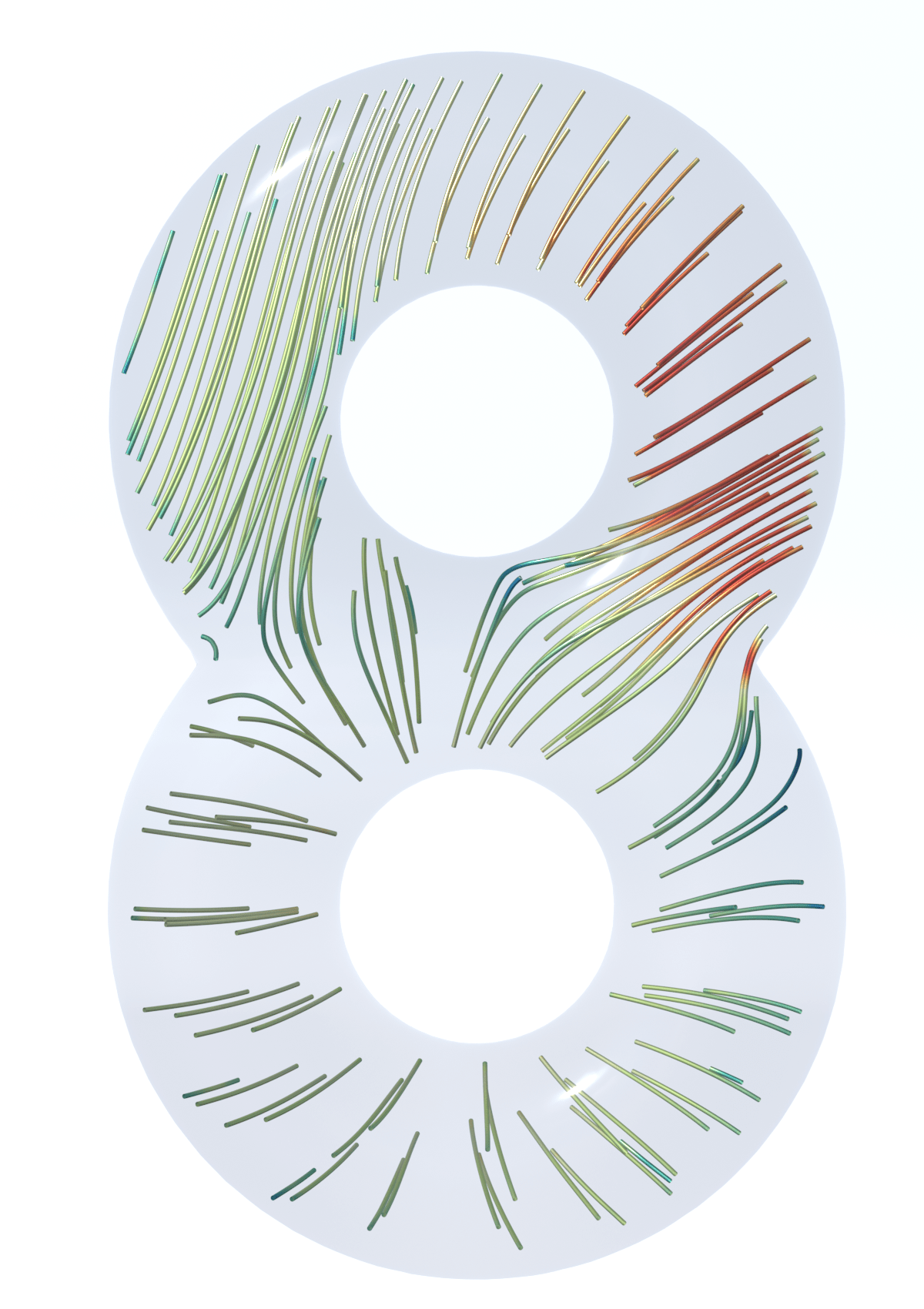}
	\caption{4-component Hodge decomposition on figure 8. From left to right: the original vector field, the normal gradient field, the tangential curl field, the tangential harmonic field, and the curly gradient field. There is no normal harmonic field, as the boundary has a single connected component ($\beta_2=0$).}
	\label{fig.hd.3d.8}
\end{figure*}

In Fig.~\ref{fig.hd.3d.kitty}, a vector field on a 3D kitty-shaped domain with a spherical cavity and one handle is decomposed into five components. With $\beta_1 = 1$ and $\beta_2 = 1$, all $5$ components can be nonzero. Figs.~\ref{fig.hd.3d.kettlebell} and~\ref{fig.buddha} show the decomposition on topologically modified kettlebell and buddha models, both with nonzero $\beta_1$ and $\beta_2$. For the 3D figure-8 model in Fig.~\ref{fig.hd.3d.8}, $\beta_1=2$ but $\beta_2=0$, thus the space of normal harmonic fields $h_t$ vanishes, leaving only four nonzero components. 

\subsection{Empirical convergence analysis}

\noindent{\bf Eigenvalue error analysis on simple shapes.} 
The convergence of the Hodge Laplacians for compact domains in the Cartesian grids has been discussed in \cite{ribando2022graph}. Here we experimented with additional examples to show the homological preservation of DEC methods and the eigenvalue error analysis on simple shapes with known eigenvalues. In the 2D case, we consider the example as in Arnold's book \cite{arnold2018finite}, where the compact domain is given by $M = (0,3)\times(0,3)\backslash ([2/3, 2]\times [3/4,2])$. Arnold showed that for this example, the standard finite element approach is insufficient and leads to entirely wrong eigenvalues for the vector Laplacian under normal (or tangential) boundary conditions, which in our case corresponds to $L_{1,n}$ (or $L_{1,t}$). In contrast, the finite element exterior calculus (FEEC) provides correct results with the first eigenvalue being 0. Note that the kernel dimension of the Laplacian is determined by the first Betti number, i.e., the number of holes inside the domain. FEEC accurately captures the topology of the underlying domain. Here we show that our extension to DEC can also provide accurate eigenvalue results and preserve the homology of the underlying domain. In Table~\ref{table.conv.arnold}, we present, as in \cite{arnold2018finite}, the first two eigenvalues computed using DEC with different grid sizes for domain $M$ with grid boundary $[-0.1, 3.1]^2$. One can see that the Laplacian has an exactly one-dimensional kernel and also its second eigenvalue converges to the correct value, given as 0.617.
\renewcommand{\arraystretch}{1.5}
\begin{table}
	\caption{First and second eigenvalues for the vector Laplacian under tangential boundary condition computed using DEC in Cartesian grids, showing the convergence to the exact solutions of 0 and 0.617.}
	\label{table.conv.arnold}
	\begin{tabular}{cccccccc}
		\hline\hline
		& grid size & $128$ & $256$ & $512$ & $1024$ & $2048$ & $ 4096 $\\
		\hline
		& $\lambda_1$ & 0  & 0 & 0 & 0 &  0 & 0 \\
		&$\lambda_2$  & 0.610  & 0.615 & 0.615 & 0.617 &  0.617 & 0.617 \\
		\hline\hline
	\end{tabular}
\end{table}

\begin{figure}[h]
	\centering
	\includegraphics[width=\textwidth]{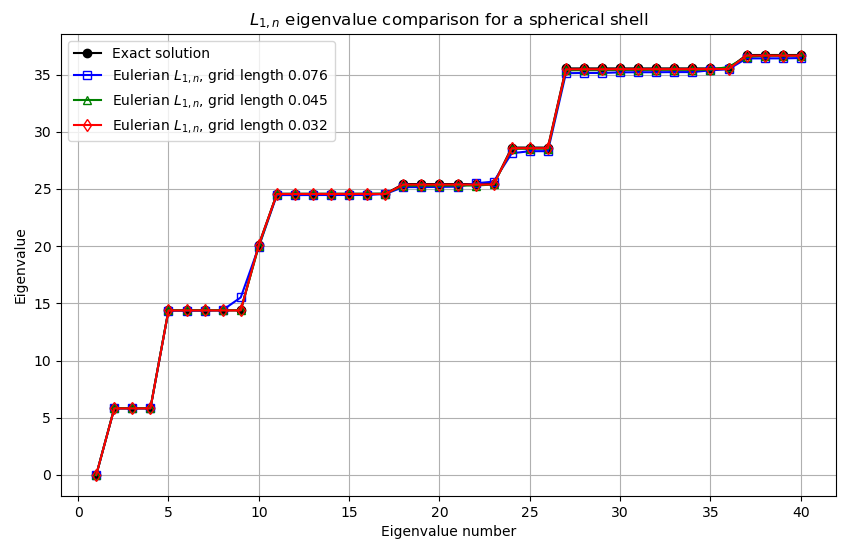}
	\caption{The first 40 eigenvalues for $L_{1, n}$ for the spherical shell with outer radius 1 and inner radius 0.3.}
	\label{fig.eigvalue.L1n.shell}
\end{figure}

To show the structure-preserving property of our framework in the 3D case, we consider the vector field Laplacian (i.e., $L_{1,n}$) under the normal boundary condition for the spherical shell centered at the origin with outer radius 1 and inner radius 0.3, as the exact spectra of the continuous Laplacian can be calculated using spherical Bessel functions. We computed the first $40$ eigenvalues of the discrete Hodge Laplacian $L_{1,n}$ and compared them with the predicted spectra. One can see in Fig.~\ref{fig.eigvalue.L1n.shell}, our numerical results closely match their corresponding exact solutions, and in particular, the computation accurately produces the correct kernel dimensional of the Laplacian, given by the second Betti number, i.e., the number of enclosed cavities. This demonstrates the accuracy of our proposed method, and its ability to preserve the manifold topology. 

\noindent{\bf Decomposition error analysis on simple shapes.}

For the 2d case, we present one example of Hodge decomposition for a flat annulus centered at the origin, with outer radius $4$ and inner radius $1$. The formulas of the vector fields are given as follows:
\begin{align}
     \omega &= \left(2(x-y)+1, 2(x+y)+1\right)\\
     d\alpha_n &= \left(2x-\frac{15}{\log(4)}\frac{x}{x^2+y^2}, 2y-\frac{15}{\log(4)}\frac{y}{x^2+y^2}\right)\\
     \delta\beta_t &= \left(-(2y-\frac{15}{\log(4)}\frac{y}{x^2+y^2}), 2x-\frac{15}{\log(4)}\frac{x}{x^2+y^2}\right)\\
     h_n &= \frac{15}{\log(4)}\left(\frac{x}{x^2+y^2}, \frac{y}{x^2+y^2}\right)\\
     h_t &= \frac{15}{\log(4)}\left(\frac{-y}{x^2+y^2}, \frac{x}{x^2+y^2}\right)\\
     \eta &= (1, 1)
\end{align}
Here $\omega = d\alpha_n + \delta\beta_t + h_n + h_t + \eta$ with $d\alpha_n, \delta\beta_t, h_n , h_t$ and $\eta$ being the components of $\omega$ in subspaces of the Hodge decomposition \eqref{eq.hd.5subspaces} following the same order. We calculated the $L^2$ norms of the decomposed components of $\omega$ using both these formulas and our framework with respect to the Hodge $L^2$ product \eqref{eq.l2innerProd.discrete.t} on the tangent support, where the annulus is modeled as a sublevel set of the signed distance function to its boundary in the grid box $[-4.1, 4.1]^2$ with different grid sizes. The results, along with the relative error between the decomposed components and their counterparts from formulas, are shown in Table~\ref{table.l2error.hd.annulus}, which demonstrates the accuracy of our framework for the Hodge decomposition in the Cartesian case.

\begin{table}
	\caption{Comparisons for the 2D annulus model between components obtained using their formula and numerically by our proposed framework measured by the Hodge $L^2$-norm induced by \eqref{eq.l2innerProd.discrete.t} on the tangent support.}
	\label{table.l2error.hd.annulus}
	\begin{tabular}{lccccc}
		\hline\hline
		Grid size  & Component & exact $L^2 norm$  & $L^2 norm$ & Relative Error\\
		\hline
		\multirow{5}{*}{100} 
            & $d\alpha_n$ & 24.178155 & 24.069865 &  0.094824\\
		& $\delta\beta_t$ & 24.145964 & 24.145907 & 0.001984\\
		& $h_n$ & 31.959631 & 31.944683 & 0.008224\\
            & $h_t$ & 31.946590 & 31.921890 & 0.001339\\
            & $\eta$ & 9.710827 & 9.972311 & 0.230518\\
		\hline
		\multirow{5}{*}{200} 
            & $d\alpha_n$ & 24.123723 & 24.078006 & 0.062341\\
		& $\delta\beta_t$ & 24.137028 & 24.137023 & 0.000589\\
		& $h_n$ & 31.922784 & 31.934750 & 0.006385\\
            & $h_t$ & 31.937387 & 31.930619 & 0.000384\\
            & $\eta$ & 9.707744 & 9.838218 & 0.158971\\
            \hline
		\multirow{5}{*}{300} 
            & $d\alpha_n$ & 24.122660 & 24.099589 & 0.043752\\
		& $\delta\beta_t$ & 24.135332 & 24.135331 & 0.000301\\
		& $h_n$ & 31.927916 & 31.930511 & 0.005049\\
            & $h_t$ & 31.935581 & 31.932422 & 0.000193\\
            & $\eta$ & 9.707372 & 9.763358 & 0.107245\\
		\hline\hline
	\end{tabular}
\end{table}

For the 3D case, we consider vector fields following \cite{poelke2017hodge} with the following formulas:
\begin{align}
	d\alpha_n &= (x, y, z)\\
	\delta\beta_t &= (y, -x, 0)\\
	h_n &= \frac{1}{\left(x^2+y^2+z^2\right)^{3/2}}(x, y, z)\\
	h_t &= \frac{1}{x^2+y^2}(y, -x, 0)\\
	\eta & = -\frac{1}{2}(1,1,1)
\end{align}
and three types of domains where these vector fields are defined, i.e., a solid ball, a spherical shell, and a torus. In our experiments, we choose the following three domains: a solid unit ball centered at the origin, a spherical shell centered at the origin with outer radius $1$ and inner radius $0.5$, and a torus centered at the origin with major radius $1$ and minor radius $0.5$. Those domains are modeled as sublevel sets of the signed distance functions to their boundaries, with grid box $[-1.1,1.1]^3$ used for the solid ball and the spherical shell, and $[-1.65,1.65]^3$ for the solid torus. We consider the vector field $\omega = d\alpha_n + \delta\beta_t + \eta$ for the solid ball model as these three fields are well defined and represent elements in their corresponding subspaces of the Hodge decomposition on the unit solid ball. The comparison results between the $L^2$ norms of the decomposed components with those obtained using their formulas on the tangential support are presented in Table~\ref{table.l2error.hd.3dball}. Note that $\omega = h_n$ is normal to spheres centered at the origin and $\omega = h_t$ are tangential to any tori with a major circle parallel to the $xy$-plane, and both are harmonic except at the origin. We test $\omega = h_n$ and $\omega = h_t$ on the spherical shell domain and the solid torus domain, respectively. The results are listed in Table~\ref{table.l2error.hd.3dshell} and Table~\ref{table.l2error.hd.3dtorus}. Since the specific information of the domains and the $L^2$ errors between the smooth and the decomposed components were not reported in \cite{poelke2017hodge}, we only report our results without comparing them to theirs. As shown in the tables, the results for the components on the tangential support obtained from our framework closely match those obtained from formulas. However, we still obtain errors due to the numerical errors from the discretizations, in particular, when the support of the opposite boundary condition is used. One can see this from the results in Table~\ref{table.l2error.hd.3dshell}, the algorithm leads to a large relative error of $0.381$ for the normal harmonic form on the tangential support when grid size is $50$, while it provides $0.045$ and $0.093$ when grid size is $30$ and $70$, respectively.

\begin{figure*}[t]
	\centering
	\includegraphics[trim=2cm 1cm 2cm 1cm, clip, width=0.15\textwidth]{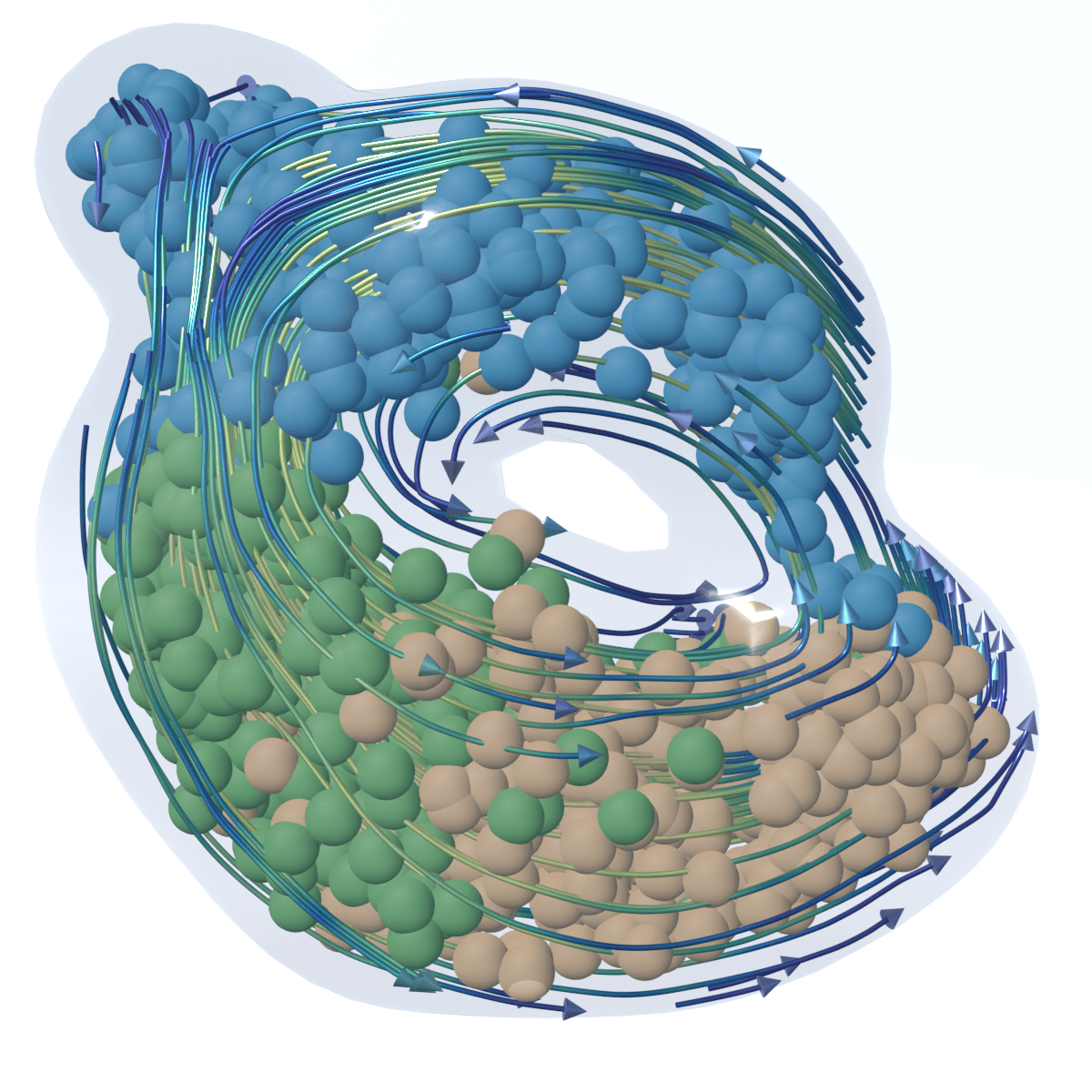}
 \raisebox{1cm}{\Large$=$}
	\includegraphics[trim=2cm 1cm 2cm 1cm, clip, width=0.15\textwidth]{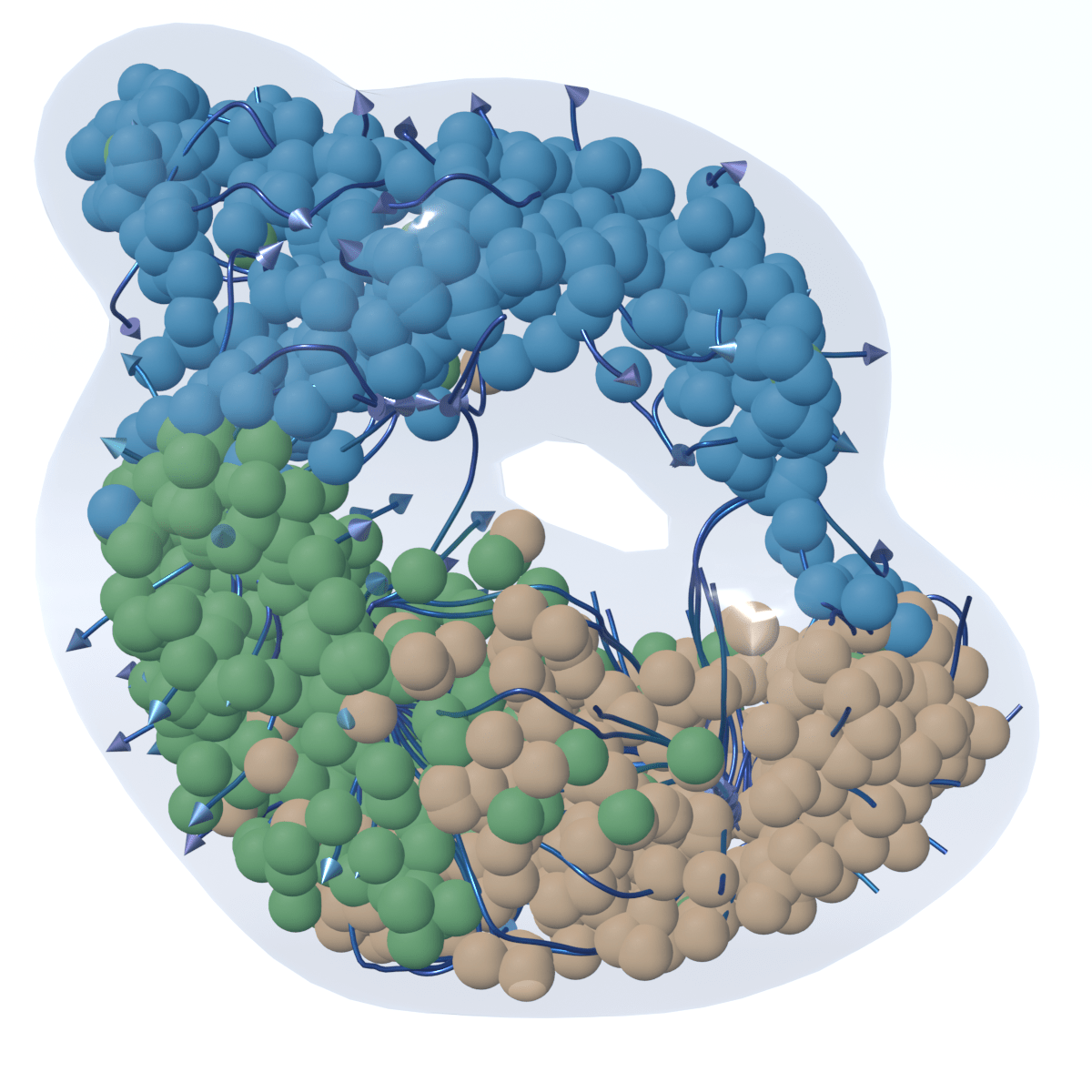}
  \raisebox{1cm}{\Large$+$}
	\includegraphics[trim=2cm 1cm 2cm 1cm, clip, width=0.15\textwidth]{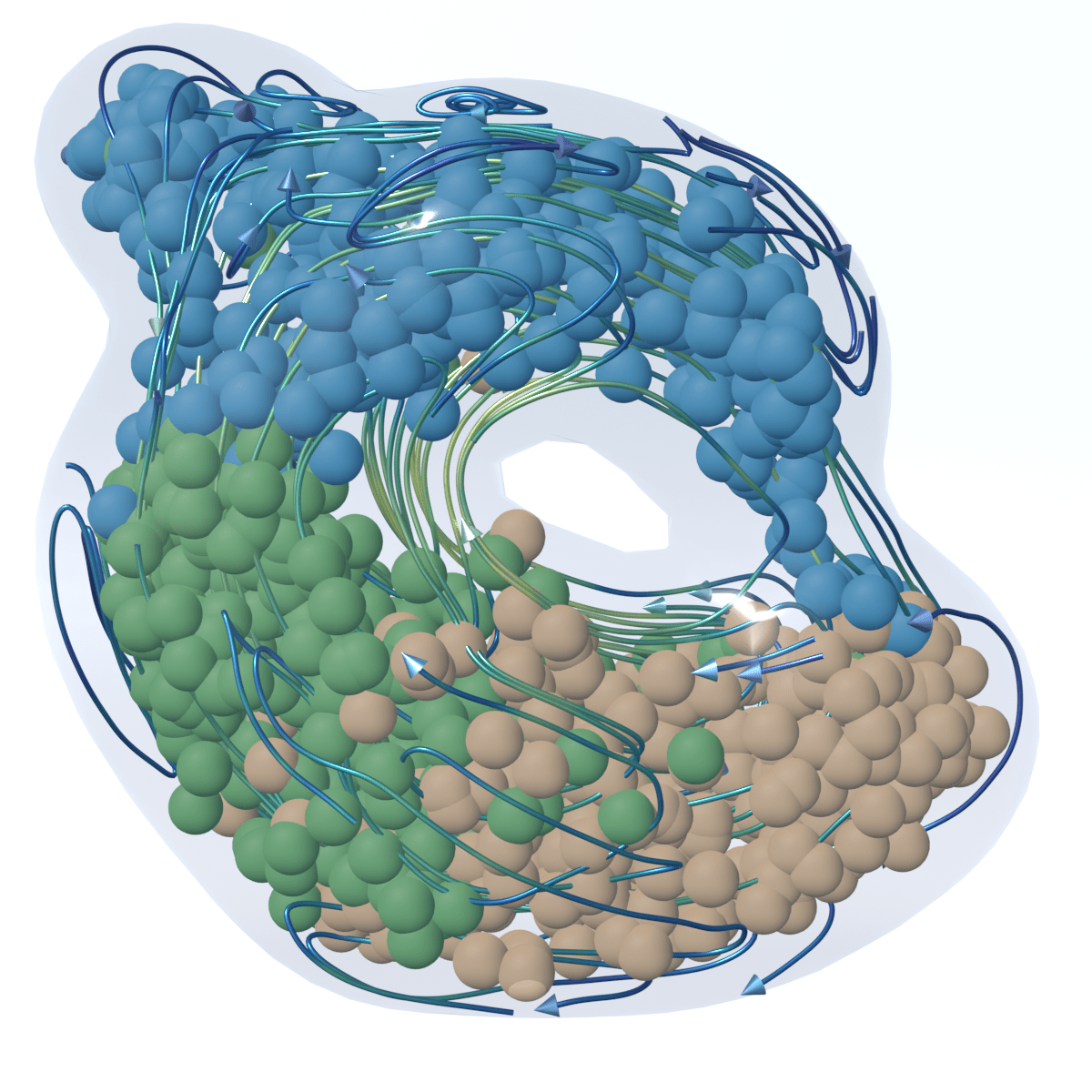}
  \raisebox{1cm}{\Large$+$}
	\includegraphics[trim=2cm 1cm 2cm 1cm, clip, width=0.15\textwidth]{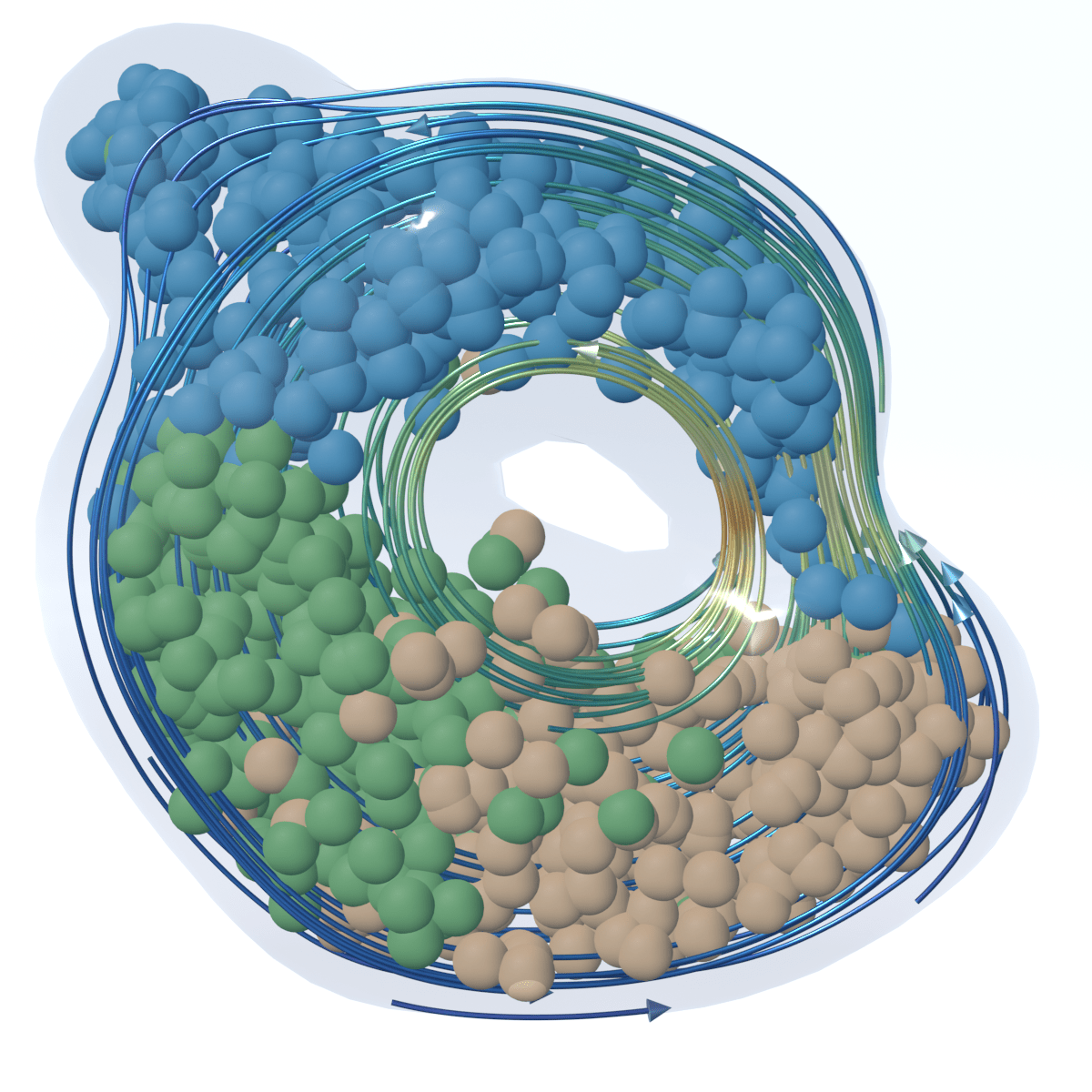}
  \raisebox{1cm}{\Large$+$}
	\includegraphics[trim=2cm 1cm 2cm 1cm, clip, width=0.15\textwidth]{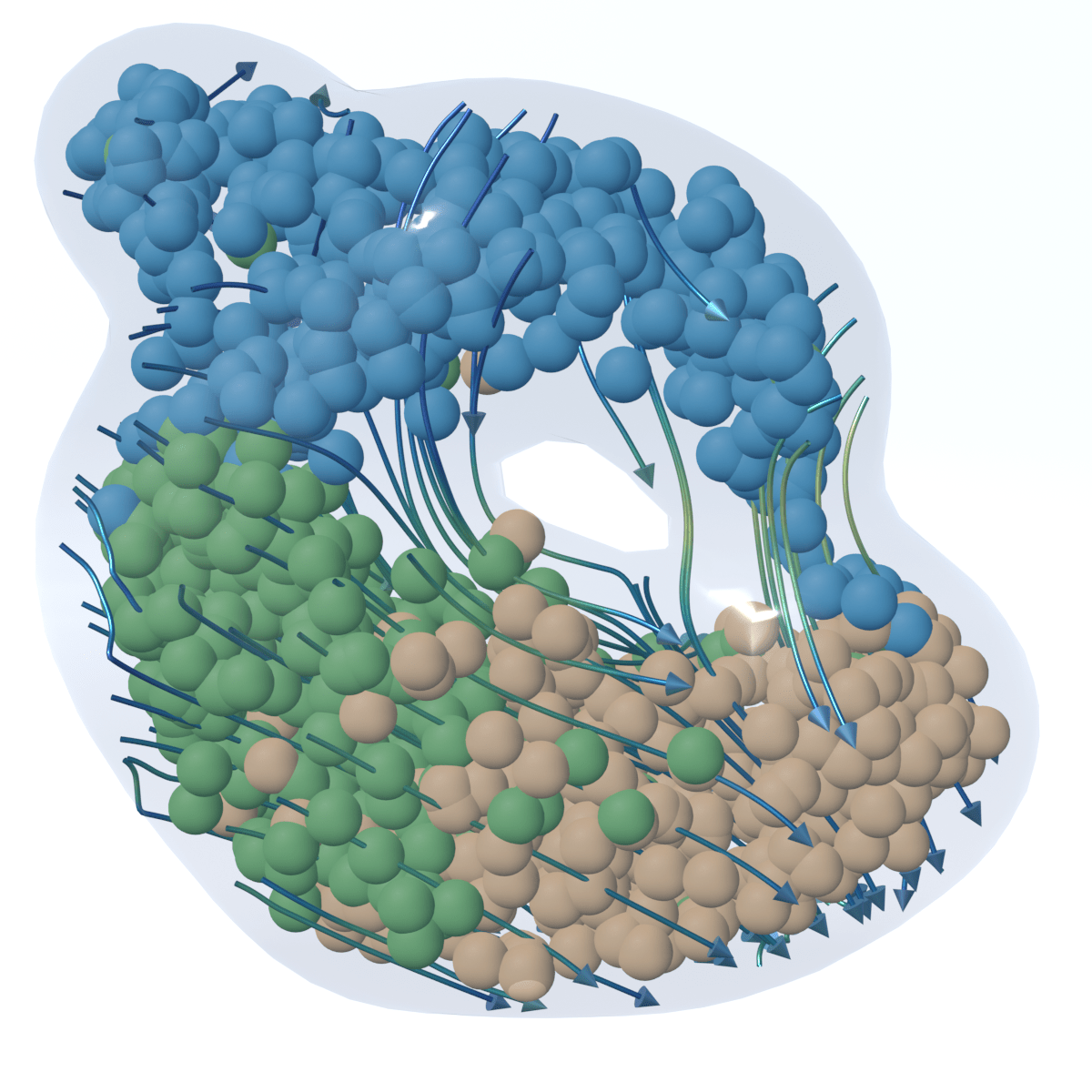}
	\caption{4-component Hodge decomposition on single-cell RNA velocity for the cell cycle dataset \cite{mahdessian2021spatiotemporal}. There is no normal harmonic field, as the boundary has a single connected component ($\beta_2=0$).}
	\label{fig.hd.RNA.cellcycle}
\end{figure*}

\begin{figure*}[t]
	\centering
	\includegraphics[trim=2cm 1cm 2cm 1cm, clip, width=0.18\textwidth]{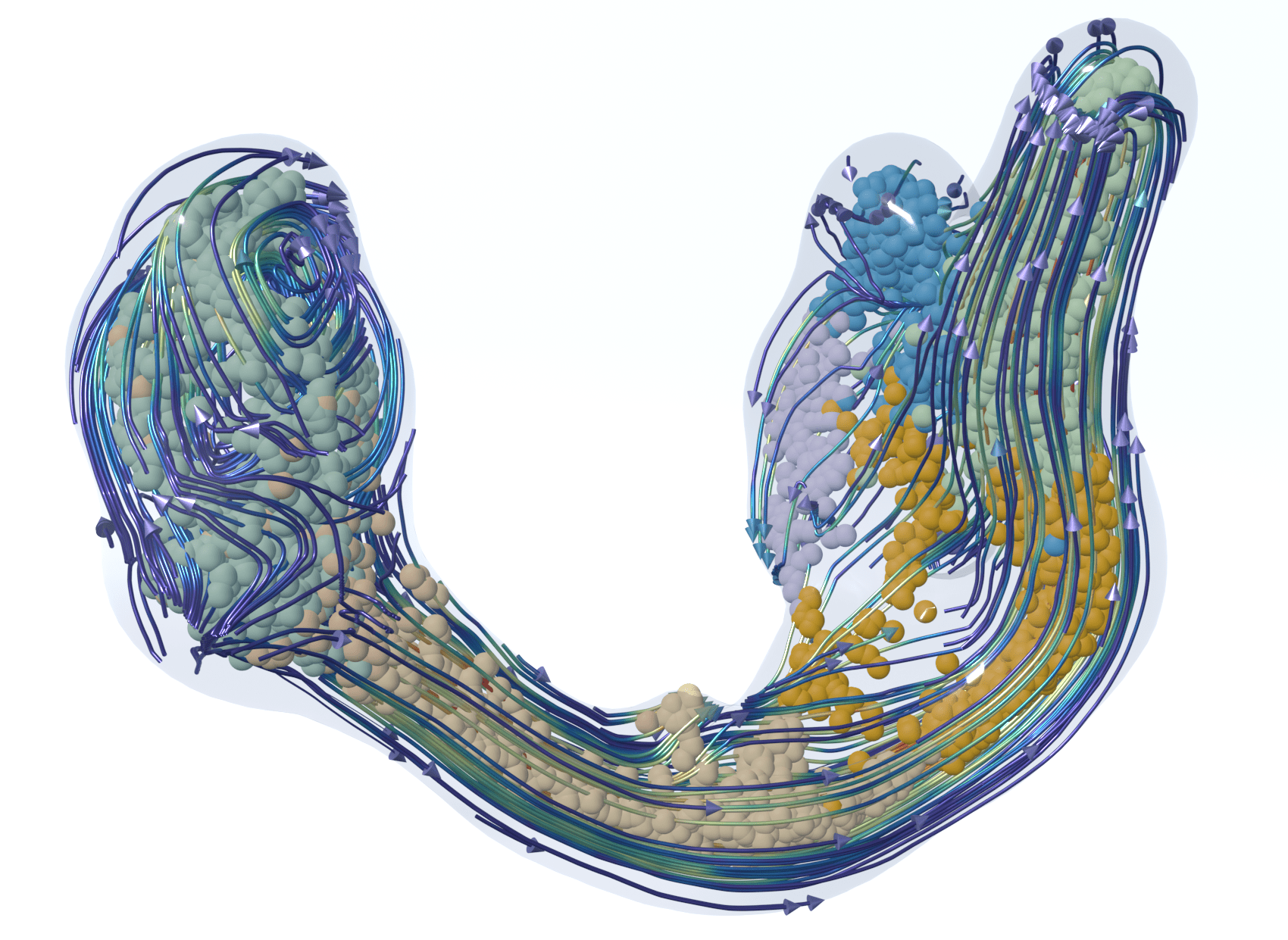}
 \raisebox{1cm}{\Large$=$}
	\includegraphics[trim=2cm 1cm 2cm 1cm, clip, width=0.18\textwidth]{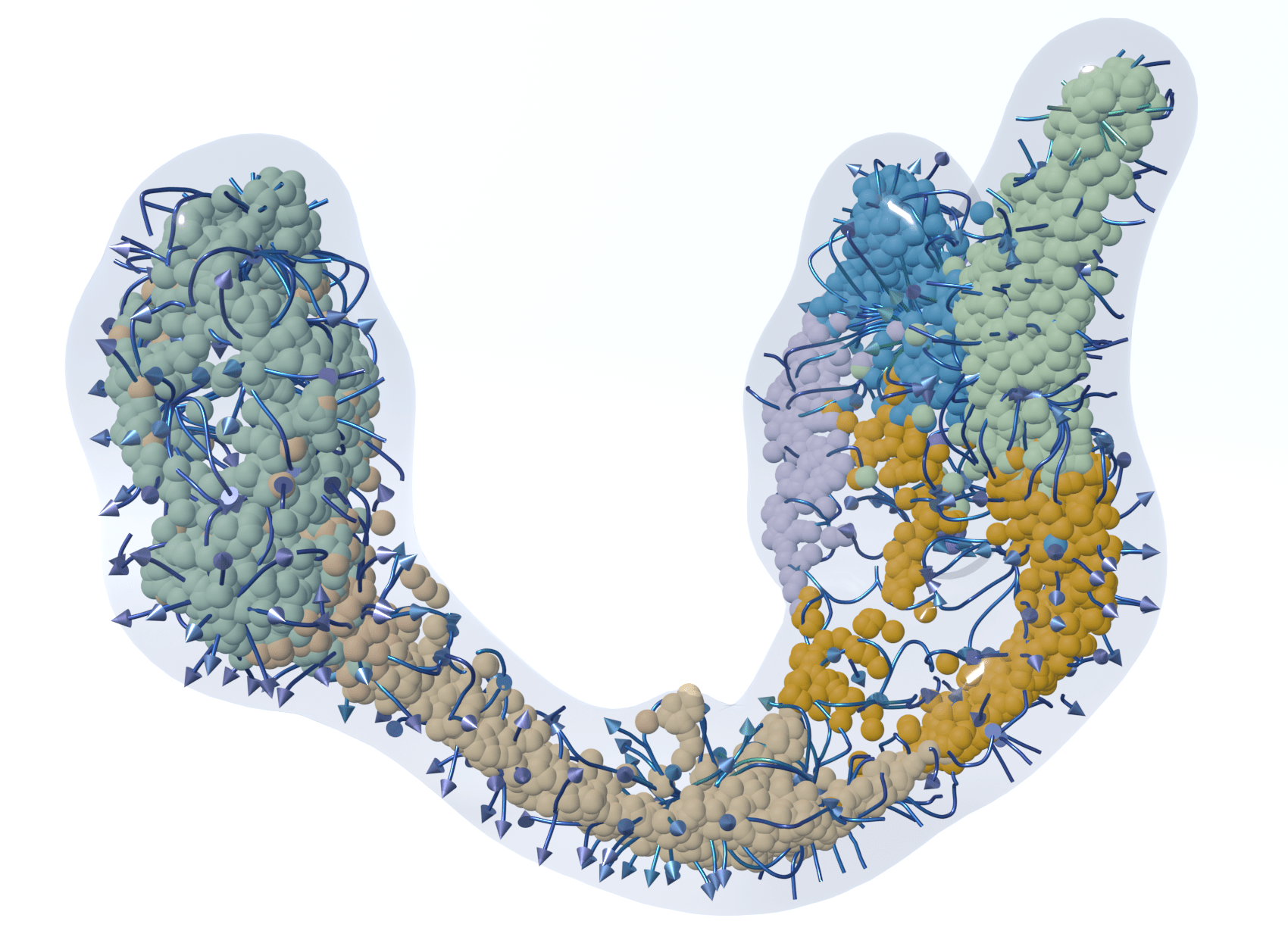}
  \raisebox{1cm}{\Large$+$}
	\includegraphics[trim=2cm 1cm 2cm 1cm, clip, width=0.18\textwidth]{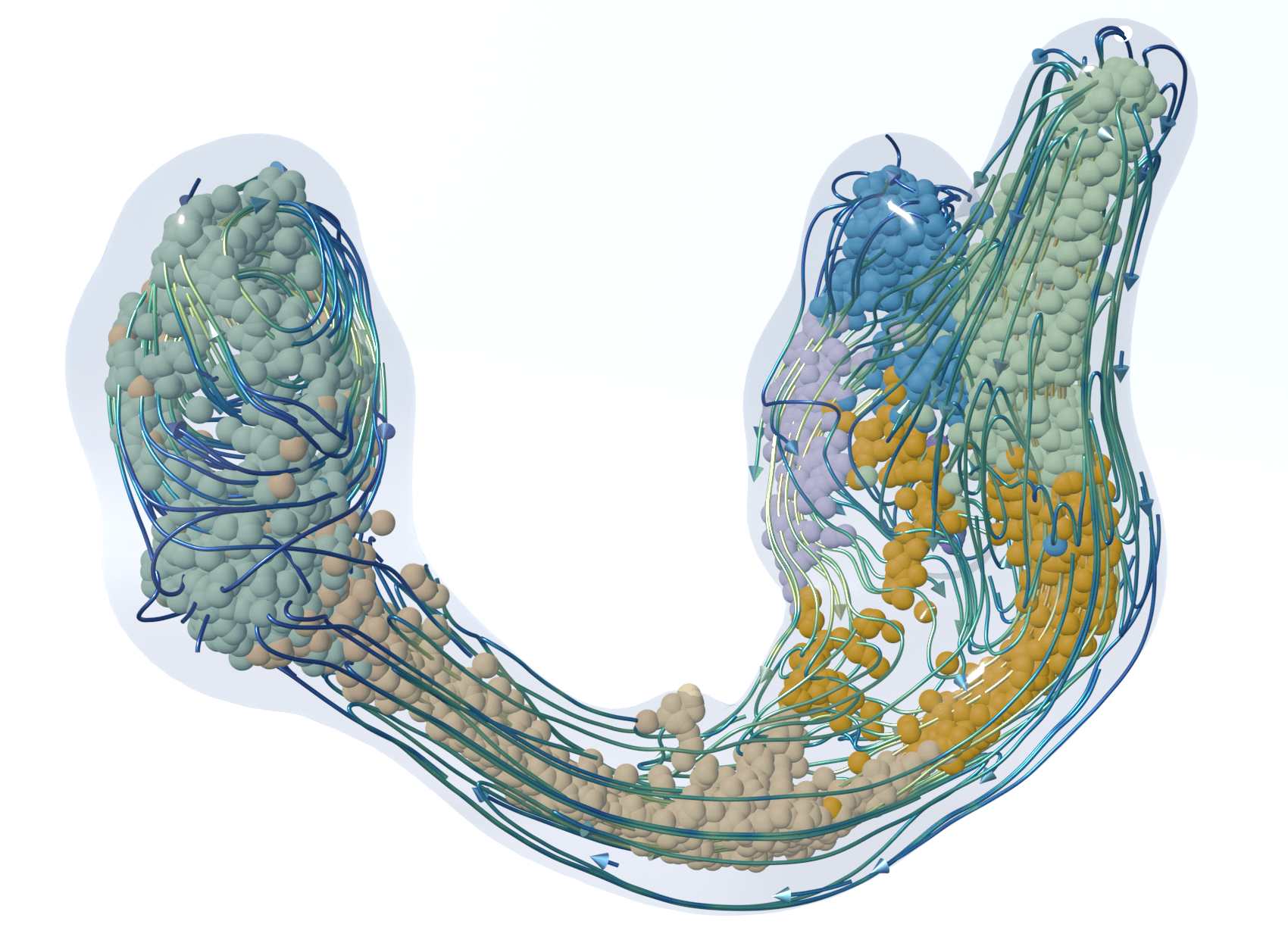}
  \raisebox{1cm}{\Large$+$}
	\includegraphics[trim=2cm 1cm 2cm 1cm, clip, width=0.18\textwidth]{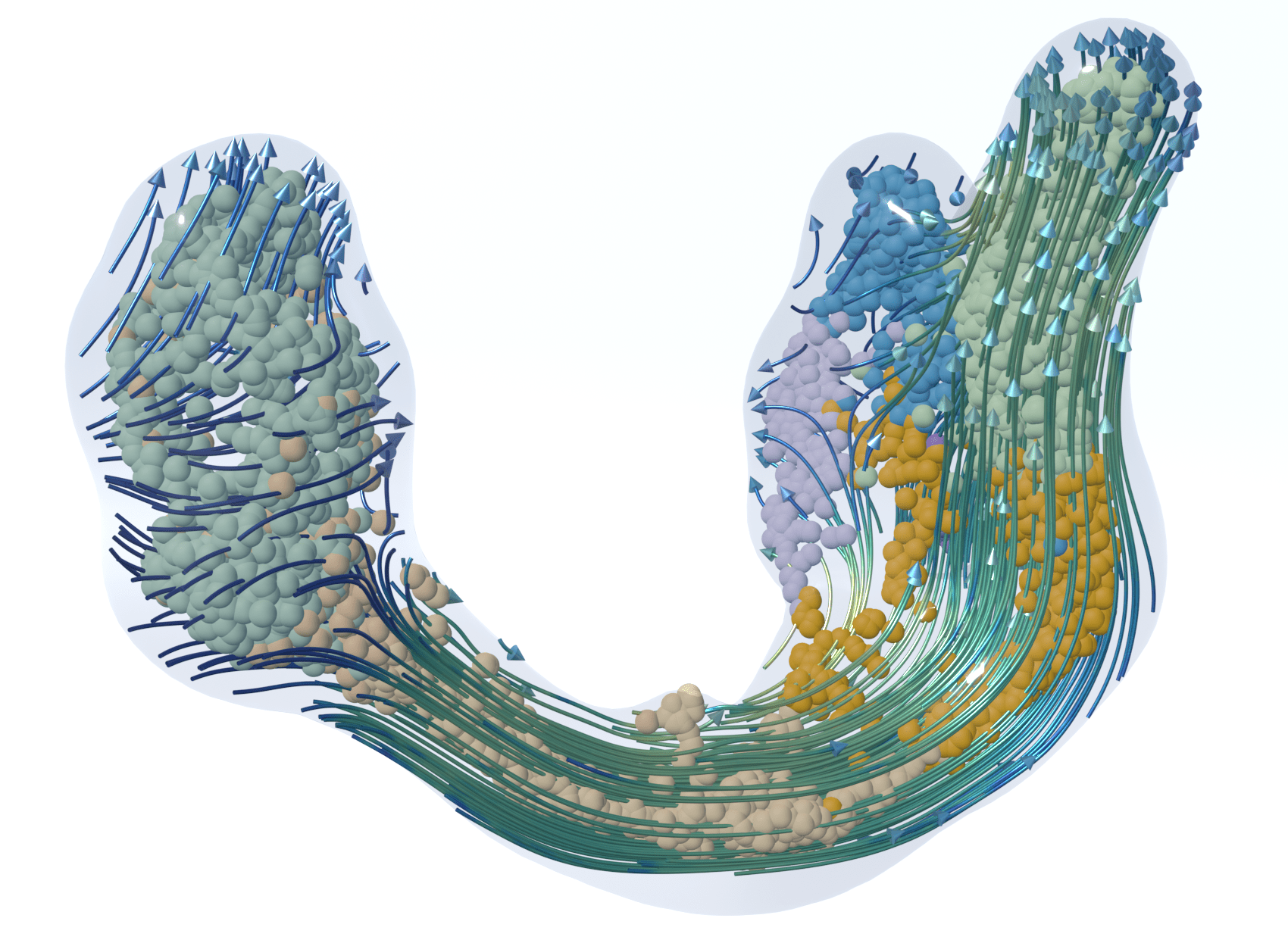}
	\caption{4-component Hodge decomposition on single-cell RNA velocity for the pancreas dataset from Scvelo \cite{bergen2020generalizing}. There is no normal and tangential harmonic field, as the manifold has no cavities and tunnels.}
	\label{fig.hd.RNA.pancreas}
\end{figure*}

\begin{figure*}[t]
	\centering
	\includegraphics[trim=2cm 1cm 2cm 1cm, clip, width=0.18\textwidth]{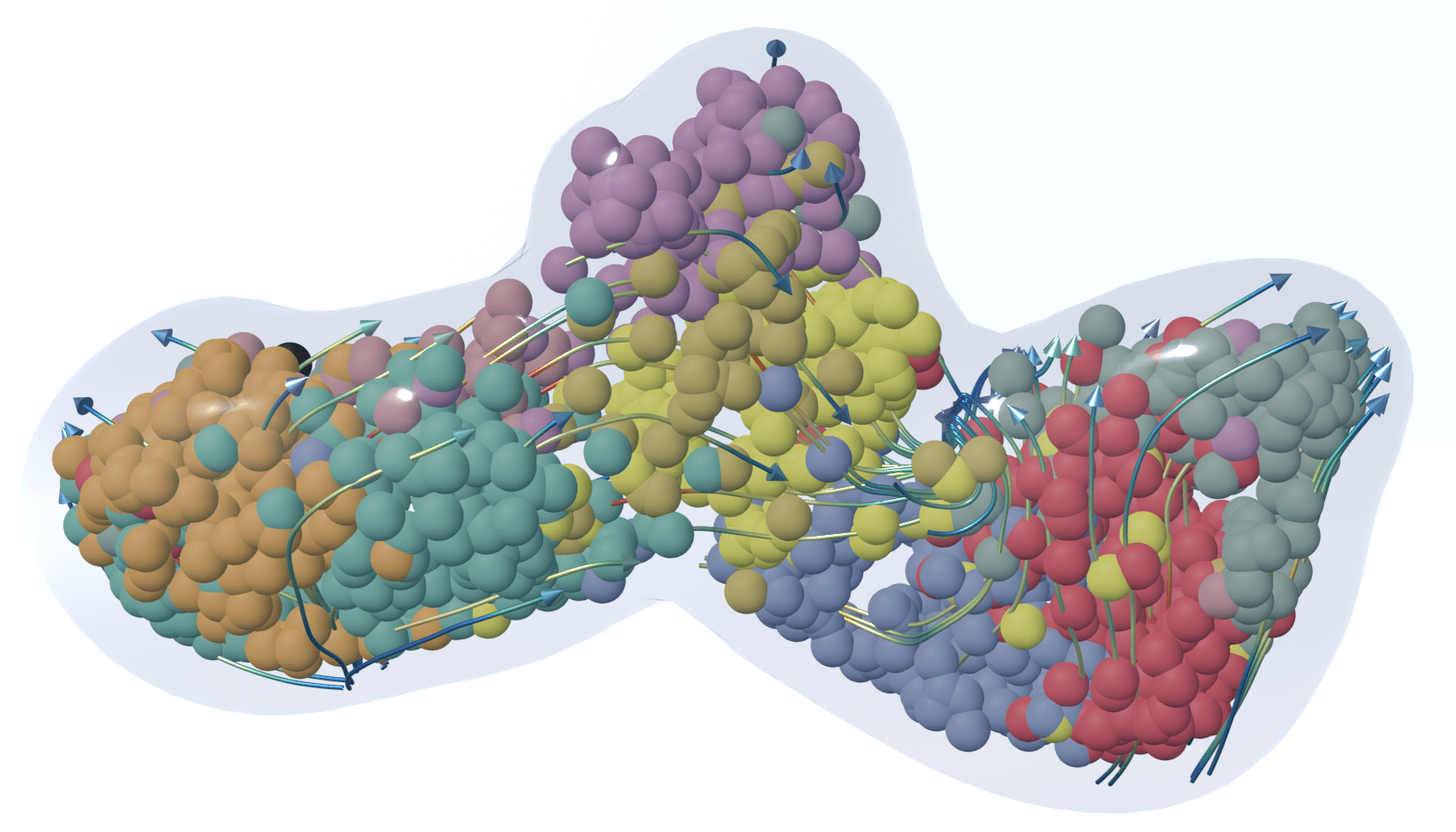}
 \raisebox{1cm}{\Large$=$}
	\includegraphics[trim=2cm 1cm 2cm 1cm, clip, width=0.18\textwidth]{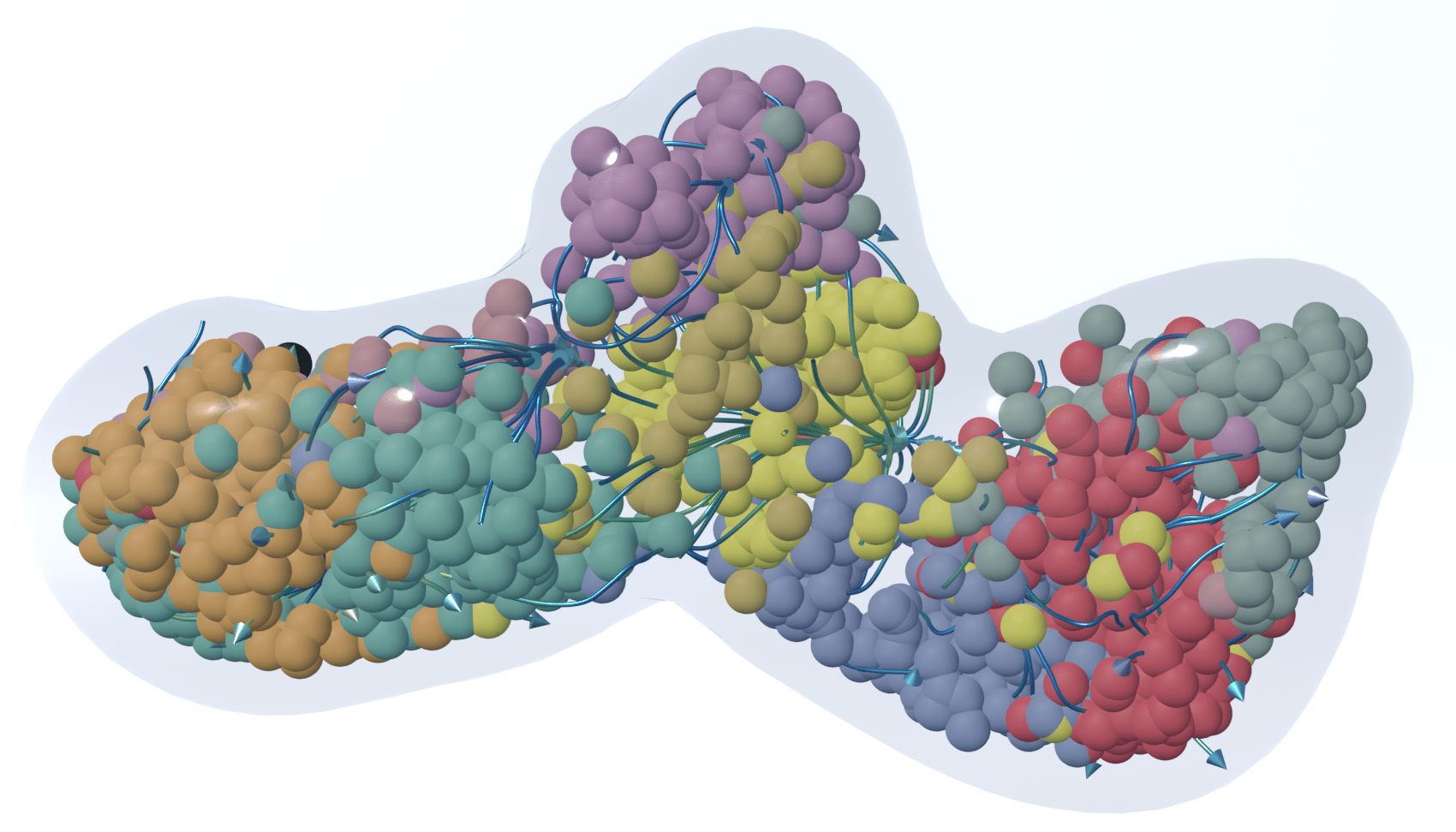}
  \raisebox{1cm}{\Large$+$}
	\includegraphics[trim=2cm 1cm 2cm 1cm, clip, width=0.18\textwidth]{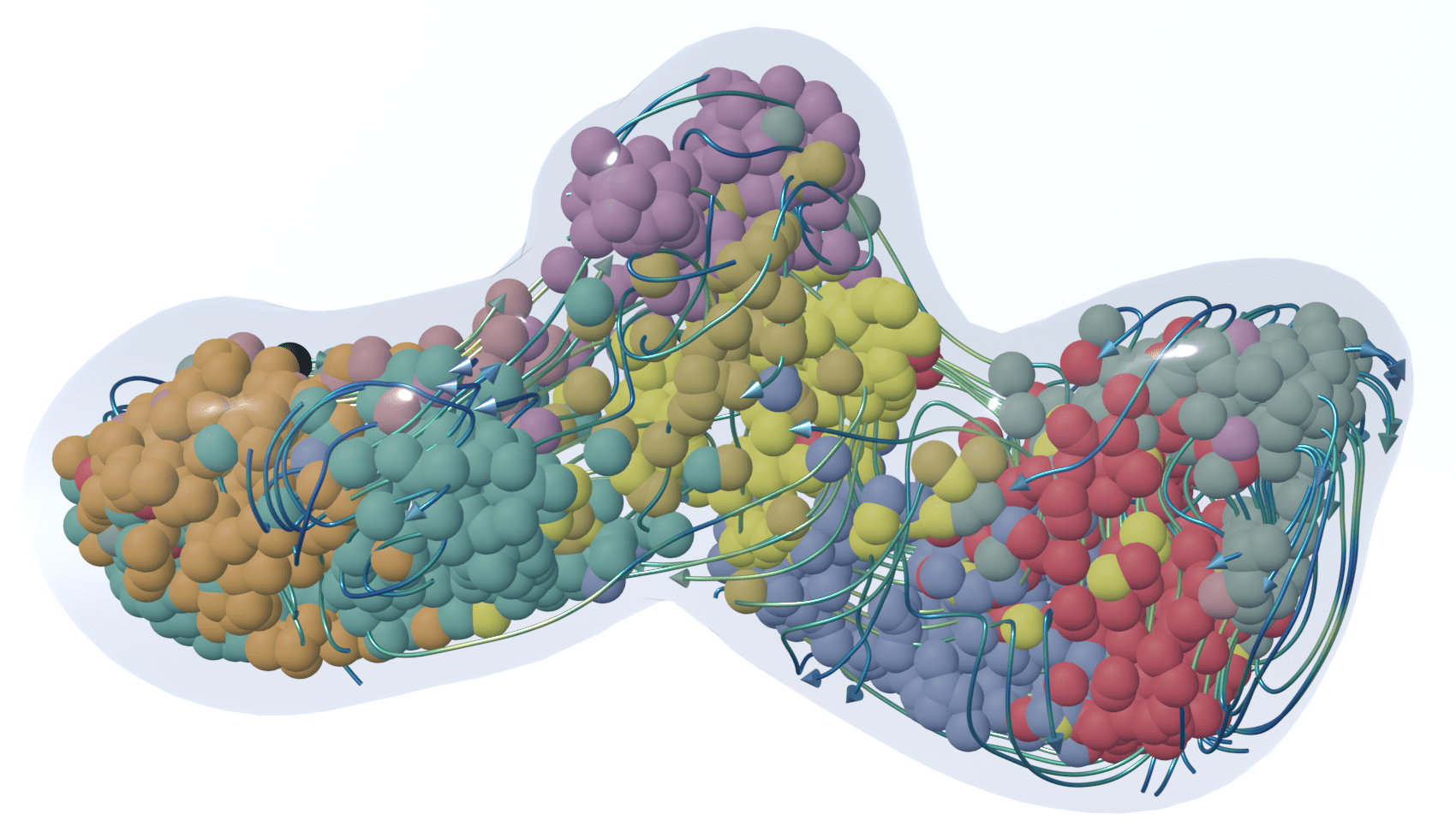}
  \raisebox{1cm}{\Large$+$}
	\includegraphics[trim=2cm 1cm 2cm 1cm, clip, width=0.18\textwidth]{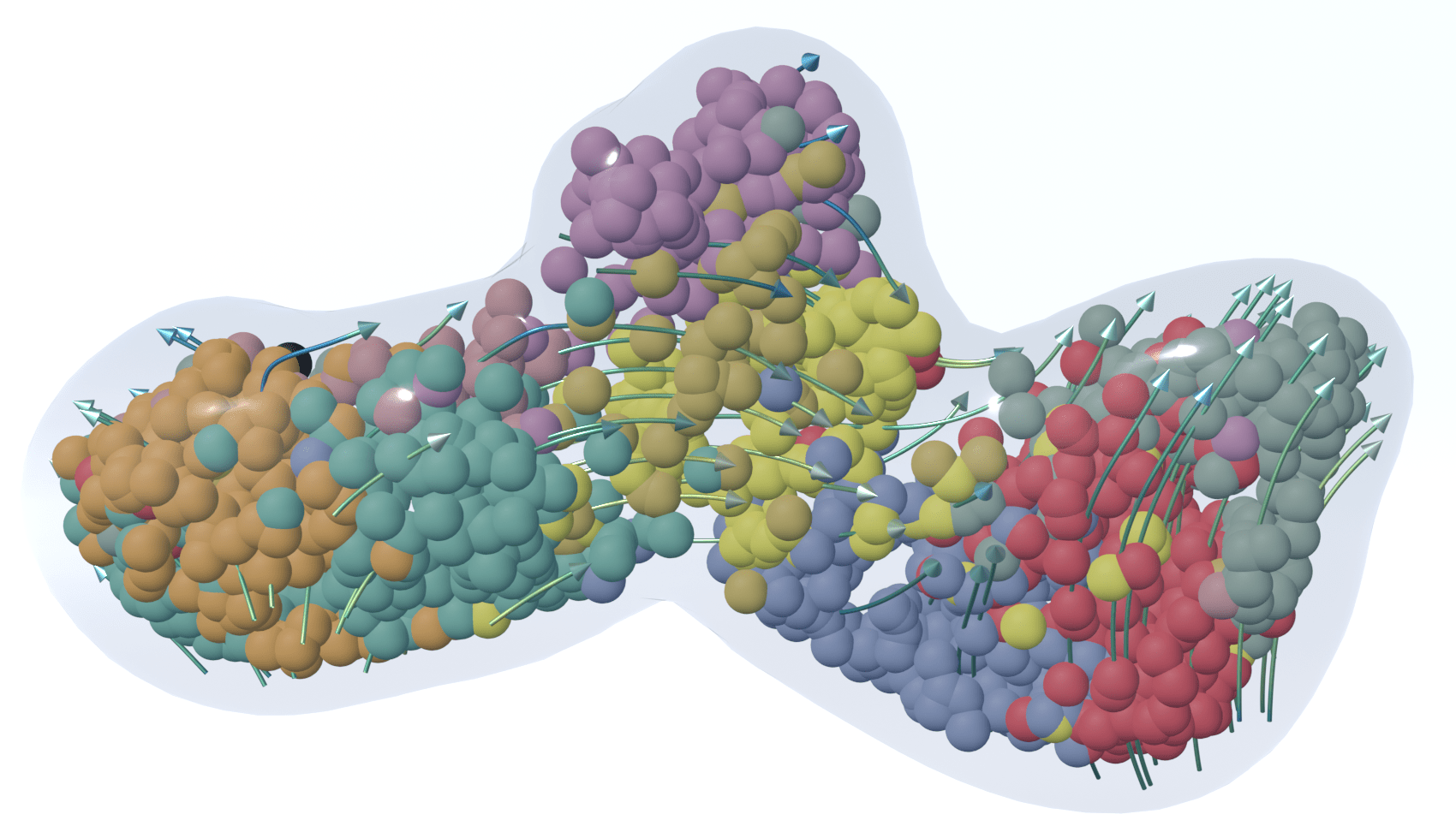}
	\caption{4-component Hodge decomposition on single-cell RNA velocity for the reprogramming Morris dataset from Cellrank \cite{lange2022cellrank}. There is no normal and tangential harmonic field, as the manifold has no cavities and tunnels.}
	\label{fig.hd.RNA.reprogramming}
\end{figure*}

\begin{table}
	\caption{Comparisons for the solid ball model between the components obtained using their formula and numerically by our proposed framework measured by the Hodge $L^2$-norm induced by \eqref{eq.l2innerProd.discrete.t} on the tangent support.}
	\label{table.l2error.hd.3dball}
	\begin{tabular}{lccccc}
		\hline\hline
		Grid size  & Component & exact $L^2 norm$  & $L^2 norm$ & Relative Error\\
		\hline
		\multirow{5}{*}{30} 
		& $d\alpha_n$ & 1.575355 & 1.597348 & 0.228610\\
		& $\delta\beta_t$ & 1.289066 & 1.289065 & 0.001302\\
		& $h_n$ & 0 & 0.000000 & N/A \\
		& $h_t$ & 0 & 0.000000 & N/A \\
		& $\eta$ & 1.766446 & 1.746584 & 0.206464\\
		\hline
		\multirow{5}{*}{50} 
		& $d\alpha_n$ & 1.582517 & 1.644876 & 0.357045\\
		& $\delta\beta_t$ & 1.292557 & 1.292557 & 0.000572\\
		& $h_n$ & 0 & 0.000000 & N/A \\
		& $h_t$ & 0 & 0.000000 & N/A \\
		& $\eta$ & 1.770509 & 1.712731 & 0.330820\\
		\hline
		\multirow{5}{*}{70} 
		& $d\alpha_n$ & 1.584136 & 1.607020 & 0.193919\\
		& $\delta\beta_t$ & 1.293478 & 1.293478 & 0.000349\\
		& $h_n$ & 0 & 0.000000 & N/A \\
		& $h_t$ & 0 & 0.000000 & N/A \\
		& $\eta$ & 1.771525 & 1.750793 & 0.175687\\
		\hline\hline
	\end{tabular}
\end{table}

\begin{table}
	\caption{Comparisons for the spherical shell model between the components obtained using their formula and numerically by our proposed framework measured by the Hodge $L^2$-norm induced by \eqref{eq.l2innerProd.discrete.t} on the tangent support.}
	\label{table.l2error.hd.3dshell}
	\begin{tabular}{lccccc}
		\hline\hline
		Grid size  & Component & exact $L^2 norm$  & $L^2 norm$ & Relative Error\\
		\hline
		\multirow{5}{*}{30} 
		& $d\alpha_n$ & 0 & 0.029324 & N/A \\
		& $\delta\beta_t$ & 0 & 0.011664 & N/A \\
		& $h_n$ & 3.568942 & 3.565251 &  0.045494 \\
		& $h_t$ & 0 & 0.000000 & N/A \\
		& $\eta$ & 0 & 0.159183 & N/A \\
		\hline
		\multirow{5}{*}{50} 
		& $d\alpha_n$ & 0 & 1.204563 & N/A \\
		& $\delta\beta_t$ & 0 & 0.003957 & N/A\\
		& $h_n$ & 3.561405 & 3.312227 & 0.380796 \\
		& $h_t$ & 0 & 0.000000 & N/A \\
		& $\eta$ & 0 & 0.511636 & N/A \\
		\hline
		\multirow{5}{*}{70} 
		& $d\alpha_n$ & 0 & 0.320805 & N/A \\
		& $\delta\beta_t$ & 0 & 0.001973 & N/A \\
		& $h_n$ & 3.549885 & 3.534490 & 0.093233 \\
		& $h_t$ & 0 & 0.000000 & N/A \\
		& $\eta$ & 0 & 0.078381 & N/A \\
		\hline\hline
	\end{tabular}
\end{table}

\begin{table}
	\caption{Comparisons for the torus model between the components obtained using their formula and numerically by our proposed framework measured by the Hodge $L^2$-norm induced by \eqref{eq.l2innerProd.discrete.t} on the tangent support.}
	\label{table.l2error.hd.3dtorus}
	\begin{tabular}{lccccc}
		\hline\hline
		Grid size  & Component & exact $L^2 norm$  & $L^2 norm$ & Relative Error\\
		\hline
		\multirow{5}{*}{30} 
		& $d\alpha_n$ & 0 & 0.010600 & N/A \\
		& $\delta\beta_t$ & 0 & 0.005725 & N/A \\
		& $h_n$ & 0 & 0.000000 & N/A \\
		& $h_t$ & 2.276355 & 2.276297 & 0.007139 \\
		& $\eta$ & 0 & 0.010906 & N/A \\
		\hline
		\multirow{5}{*}{50} 
		& $d\alpha_n$ & 0 & 0.003610 & N/A \\
		& $\delta\beta_t$ & 0 & 0.001969 & N/A \\
		& $h_n$ & 0 & 0.000000 & N/A \\
		& $h_t$ & 2.289593 & 2.289584 & 0.002785 \\
		& $\eta$ & 0 & 0.004874 & N/A \\
		\hline
		\multirow{5}{*}{70} 
		& $d\alpha_n$ & 0 & 0.002607 & N/A \\
		& $\delta\beta_t$ & 0 & 0.000984 & N/A \\
		& $h_n$ & 0 & 0.000000 & 0.002143 \\
		& $h_t$ & 2.292675 & 2.292670 & N/A \\
		& $\eta$ & 0 & 0.004047 & N/A \\
		\hline\hline
	\end{tabular}
\end{table}

\subsection{Sample applications on the analysis of scRNA velocity}


The Hodge decomposition has been applied to the study of single-cell RNA velocities, providing a reliable and robust way of extracting the dynamic information of cells \cite{su2024hodge}. The results show that the decomposed components distinctly reveal key cell dynamics features, such as cell cycle, bifurcation, and cell lineage differentiation. Here we illustrate some examples of the Hodge decomposition of RNA velocity fields to demonstrate that the algorithm also applies to abstract manifolds. The decomposition for RNA velocity fields is presented in Fig.~\ref{fig.hd.RNA.cellcycle} for the cell cycle dataset \cite{mahdessian2021spatiotemporal}, Fig.~\ref{fig.hd.RNA.pancreas} for the pancreas dataset from Scvelo \cite{bergen2020generalizing}, and Fig.~\ref{fig.hd.RNA.reprogramming} for the reprogramming Morris dataset from Cellrank \cite{lange2022cellrank}. The method offers a new analytical tool for examining the cell dynamics for single-cell datasets, and these decomposed components can be further studied to enhance the understanding of cell dynamics for single-cell datasets. For readers interested in further details, refer to \cite{su2024hodge}. 


\section{Conclusion}
Topology-preserving Hodge decomposition in the Eulerian representation is an important and challenging task in data science with various applications. In this work, we have presented a framework to compute the complete $5$-component orthogonal decomposition of 2D and 3D vector fields on domains embedded in regular Cartesian grids. We leveraged the correspondence between vector fields and differential forms with certain boundary conditions. The domain boundaries are encoded by isosurfaces of level-set functions. Compared to the methods on simplicial meshes \cite{razafindrazaka2019consistent,zhao20193d}, our framework greatly simplifies the data structure and discrete operators by using the vertices, edges, faces, and cells of the grid. This significantly improves the efficiency of our algorithms for evolving level set functions, which have applications in machine learning \cite{su2024persistent}. Extracting topological invariants in the Hodge decomposition is very demanding.    Our adaptation DEC ensures that the discrete Hodge decomposition preserves the crucial topological structure while conforming to specified boundary conditions. All components can be calculated by solving sparse linear systems with rank efficiency eliminated by accounting for the corresponding spaces of harmonic fields that depend only on the underlying manifold topology. The \l-orthogonality of the decomposed $5$ discrete components is also rigorously guaranteed by the linear algebra formulation. This framework has the potential to benefit various downstream applications due to the ubiquity of vector field analysis on domains in Euclidean spaces.

{\it Limitations and future work.} Our framework has so far been evaluated only on compact domains in 2D/3D Euclidean spaces. Extending it to higher dimensions should be relatively straightforward. However, a limitation lies in the domain representation: it cannot capture sharp features unless the grid resolution is sufficiently fine because the domain is modeled as a region bounded by an isocurve or an isosurface of a level-set function. Geometric details are also constrained by grid resolution. Therefore, one possible future extension is to incorporate adaptive data structures similar to octrees. Moreover, we only considered the Hodge decomposition with each component under just either normal or tangential boundary conditions. Future work could involve generalizing this to mixed boundary conditions as in \cite{zhao20193d} for the mesh setting. Another interesting variation worth exploring is the implementation of high-order Galerkin-type Hodge stars rather than diagonal ones, and studying their impact on the convergence rate. A possible speedup for the saddle point problem in Eq.~\eqref{eq.constrained.minimization} is through the Schur complement reduction~\cite{benzi2005numerical} initialized with our direct approach.


\section*{Acknowledgments}
This work was supported in part by NIH grants  R01GM126189, R01AI164266, and R01AI146210, NSF grants DMS-2052983,  DMS-1761320, and IIS-1900473,  NASA grant 80NSSC21M0023,  MSU Research Foundation,  Bristol-Myers Squibb 65109, and Pfizer. 
\bibliographystyle{ACM-Reference-Format}
\bibliography{refs}


\begin{thebibliography}{50}


\ifx \showCODEN    \undefined \def \showCODEN     #1{\unskip}     \fi
\ifx \showDOI      \undefined \def \showDOI       #1{#1}\fi
\ifx \showISBNx    \undefined \def \showISBNx     #1{\unskip}     \fi
\ifx \showISBNxiii \undefined \def \showISBNxiii  #1{\unskip}     \fi
\ifx \showISSN     \undefined \def \showISSN      #1{\unskip}     \fi
\ifx \showLCCN     \undefined \def \showLCCN      #1{\unskip}     \fi
\ifx \shownote     \undefined \def \shownote      #1{#1}          \fi
\ifx \showarticletitle \undefined \def \showarticletitle #1{#1}   \fi
\ifx \showURL      \undefined \def \showURL       {\relax}        \fi
\providecommand\bibfield[2]{#2}
\providecommand\bibinfo[2]{#2}
\providecommand\natexlab[1]{#1}
\providecommand\showeprint[2][]{arXiv:#2}

\bibitem[Arnold(2018)]%
        {arnold2018finite}
\bibfield{author}{\bibinfo{person}{Douglas~N Arnold}.}
  \bibinfo{year}{2018}\natexlab{}.
\newblock \bibinfo{booktitle}{\emph{Finite element exterior calculus}}.
\newblock \bibinfo{publisher}{SIAM}.
\newblock


\bibitem[Batty et~al\mbox{.}(2007)]%
        {batty2007fast}
\bibfield{author}{\bibinfo{person}{Christopher Batty},
  \bibinfo{person}{Florence Bertails}, {and} \bibinfo{person}{Robert Bridson}.}
  \bibinfo{year}{2007}\natexlab{}.
\newblock \showarticletitle{A fast variational framework for accurate
  solid-fluid coupling}.
\newblock \bibinfo{journal}{\emph{ACM Transactions on Graphics (TOG)}}
  \bibinfo{volume}{26}, \bibinfo{number}{3} (\bibinfo{year}{2007}),
  \bibinfo{pages}{100--es}.
\newblock


\bibitem[Bell and Hirani(2012)]%
        {bell2012pydec}
\bibfield{author}{\bibinfo{person}{Nathan Bell} {and} \bibinfo{person}{Anil~N
  Hirani}.} \bibinfo{year}{2012}\natexlab{}.
\newblock \showarticletitle{PyDEC: software and algorithms for discretization
  of exterior calculus}.
\newblock \bibinfo{journal}{\emph{ACM Transactions on Mathematical Software
  (TOMS)}} \bibinfo{volume}{39}, \bibinfo{number}{1} (\bibinfo{year}{2012}),
  \bibinfo{pages}{1--41}.
\newblock


\bibitem[Benzi et~al\mbox{.}(2005)]%
        {benzi2005numerical}
\bibfield{author}{\bibinfo{person}{Michele Benzi}, \bibinfo{person}{Gene~H
  Golub}, {and} \bibinfo{person}{J{\"o}rg Liesen}.}
  \bibinfo{year}{2005}\natexlab{}.
\newblock \showarticletitle{Numerical solution of saddle point problems}.
\newblock \bibinfo{journal}{\emph{Acta numerica}}  \bibinfo{volume}{14}
  (\bibinfo{year}{2005}), \bibinfo{pages}{1--137}.
\newblock


\bibitem[Bergen et~al\mbox{.}(2020)]%
        {bergen2020generalizing}
\bibfield{author}{\bibinfo{person}{Volker Bergen}, \bibinfo{person}{Marius
  Lange}, \bibinfo{person}{Stefan Peidli}, \bibinfo{person}{F~Alexander Wolf},
  {and} \bibinfo{person}{Fabian~J Theis}.} \bibinfo{year}{2020}\natexlab{}.
\newblock \showarticletitle{Generalizing RNA velocity to transient cell states
  through dynamical modeling}.
\newblock \bibinfo{journal}{\emph{Nat. Biotechnol.}}  \bibinfo{volume}{38}
  (\bibinfo{year}{2020}), \bibinfo{pages}{1408--1414}.
\newblock


\bibitem[Bhatia et~al\mbox{.}(2012)]%
        {bhatia2012helmholtz}
\bibfield{author}{\bibinfo{person}{Harsh Bhatia}, \bibinfo{person}{Gregory
  Norgard}, \bibinfo{person}{Valerio Pascucci}, {and}
  \bibinfo{person}{Peer-Timo Bremer}.} \bibinfo{year}{2012}\natexlab{}.
\newblock \showarticletitle{The Helmholtz-Hodge decomposition—a survey}.
\newblock \bibinfo{journal}{\emph{IEEE Transactions on visualization and
  computer graphics}} \bibinfo{volume}{19}, \bibinfo{number}{8}
  (\bibinfo{year}{2012}), \bibinfo{pages}{1386--1404}.
\newblock


\bibitem[Cang and Wei(2017)]%
        {cang2017topologynet}
\bibfield{author}{\bibinfo{person}{Zixuan Cang} {and} \bibinfo{person}{Guo-Wei
  Wei}.} \bibinfo{year}{2017}\natexlab{}.
\newblock \showarticletitle{TopologyNet: Topology based deep convolutional and
  multi-task neural networks for biomolecular property predictions}.
\newblock \bibinfo{journal}{\emph{PLoS computational biology}}
  \bibinfo{volume}{13}, \bibinfo{number}{7} (\bibinfo{year}{2017}),
  \bibinfo{pages}{e1005690}.
\newblock


\bibitem[Cantarella et~al\mbox{.}(2002)]%
        {cantarella2002vector}
\bibfield{author}{\bibinfo{person}{Jason Cantarella}, \bibinfo{person}{Dennis
  DeTurck}, {and} \bibinfo{person}{Herman Gluck}.}
  \bibinfo{year}{2002}\natexlab{}.
\newblock \showarticletitle{Vector calculus and the topology of domains in
  3-space}.
\newblock \bibinfo{journal}{\emph{The American Mathematical Monthly}}
  \bibinfo{volume}{109}, \bibinfo{number}{5} (\bibinfo{year}{2002}),
  \bibinfo{pages}{409--442}.
\newblock


\bibitem[Chen et~al\mbox{.}(2022)]%
        {chen2022persistent}
\bibfield{author}{\bibinfo{person}{Jiahui Chen}, \bibinfo{person}{Yuchi Qiu},
  \bibinfo{person}{Rui Wang}, {and} \bibinfo{person}{Guo-Wei Wei}.}
  \bibinfo{year}{2022}\natexlab{}.
\newblock \showarticletitle{Persistent Laplacian projected Omicron BA. 4 and
  BA. 5 to become new dominating variants}.
\newblock \bibinfo{journal}{\emph{Computers in Biology and Medicine}}
  \bibinfo{volume}{151} (\bibinfo{year}{2022}), \bibinfo{pages}{106262}.
\newblock


\bibitem[Chen et~al\mbox{.}(2021)]%
        {chen2021evolutionary}
\bibfield{author}{\bibinfo{person}{Jiahui Chen}, \bibinfo{person}{Rundong
  Zhao}, \bibinfo{person}{Yiying Tong}, {and} \bibinfo{person}{Guo-Wei Wei}.}
  \bibinfo{year}{2021}\natexlab{}.
\newblock \showarticletitle{Evolutionary de rham-hodge method}.
\newblock \bibinfo{journal}{\emph{Discrete and continuous dynamical systems.
  Series B}} \bibinfo{volume}{26}, \bibinfo{number}{7} (\bibinfo{year}{2021}),
  \bibinfo{pages}{3785}.
\newblock


\bibitem[De~Goes et~al\mbox{.}(2016)]%
        {de2016vector}
\bibfield{author}{\bibinfo{person}{Fernando De~Goes}, \bibinfo{person}{Mathieu
  Desbrun}, {and} \bibinfo{person}{Yiying Tong}.}
  \bibinfo{year}{2016}\natexlab{}.
\newblock \showarticletitle{Vector field processing on triangle meshes}.
\newblock In \bibinfo{booktitle}{\emph{ACM SIGGRAPH 2016 Courses}}.
  \bibinfo{pages}{1--49}.
\newblock


\bibitem[Desbrun et~al\mbox{.}(2006)]%
        {desbrun2006discrete}
\bibfield{author}{\bibinfo{person}{Mathieu Desbrun}, \bibinfo{person}{Eva
  Kanso}, {and} \bibinfo{person}{Yiying Tong}.}
  \bibinfo{year}{2006}\natexlab{}.
\newblock \showarticletitle{Discrete differential forms for computational
  modeling}.
\newblock In \bibinfo{booktitle}{\emph{ACM SIGGRAPH 2006 Courses}}.
  \bibinfo{pages}{39--54}.
\newblock


\bibitem[Edelsbrunner et~al\mbox{.}(2008)]%
        {edelsbrunner2008persistent}
\bibfield{author}{\bibinfo{person}{Herbert Edelsbrunner}, \bibinfo{person}{John
  Harer}, {et~al\mbox{.}}} \bibinfo{year}{2008}\natexlab{}.
\newblock \showarticletitle{Persistent homology-a survey}.
\newblock \bibinfo{journal}{\emph{Contemporary mathematics}}
  \bibinfo{volume}{453}, \bibinfo{number}{26} (\bibinfo{year}{2008}),
  \bibinfo{pages}{257--282}.
\newblock


\bibitem[Fedkiw et~al\mbox{.}(2001)]%
        {fedkiw2001visual}
\bibfield{author}{\bibinfo{person}{Ronald Fedkiw}, \bibinfo{person}{Jos Stam},
  {and} \bibinfo{person}{Henrik~Wann Jensen}.} \bibinfo{year}{2001}\natexlab{}.
\newblock \showarticletitle{Visual simulation of smoke}. In
  \bibinfo{booktitle}{\emph{Proceedings of the 28th annual conference on
  Computer graphics and interactive techniques}}. \bibinfo{pages}{15--22}.
\newblock


\bibitem[Friedrichs(1955)]%
        {friedrichs1955differential}
\bibfield{author}{\bibinfo{person}{Kurt~Otto Friedrichs}.}
  \bibinfo{year}{1955}\natexlab{}.
\newblock \showarticletitle{Differential forms on Riemannian manifolds}.
\newblock \bibinfo{journal}{\emph{Communications on Pure and Applied
  Mathematics}} \bibinfo{volume}{8}, \bibinfo{number}{4}
  (\bibinfo{year}{1955}), \bibinfo{pages}{551--590}.
\newblock


\bibitem[Grbi{\'c} et~al\mbox{.}(2022)]%
        {grbic2022aspects}
\bibfield{author}{\bibinfo{person}{Jelena Grbi{\'c}}, \bibinfo{person}{Jie Wu},
  \bibinfo{person}{Kelin Xia}, {and} \bibinfo{person}{Guo-Wei Wei}.}
  \bibinfo{year}{2022}\natexlab{}.
\newblock \showarticletitle{Aspects of topological approaches for data
  science}.
\newblock \bibinfo{journal}{\emph{Foundations of data science (Springfield,
  Mo.)}} \bibinfo{volume}{4}, \bibinfo{number}{2} (\bibinfo{year}{2022}),
  \bibinfo{pages}{165}.
\newblock


\bibitem[Guo et~al\mbox{.}(2005)]%
        {guo2005efficient}
\bibfield{author}{\bibinfo{person}{Qinghong Guo}, \bibinfo{person}{Mrinal~K
  Mandal}, {and} \bibinfo{person}{Micheal~Y Li}.}
  \bibinfo{year}{2005}\natexlab{}.
\newblock \showarticletitle{Efficient Hodge--Helmholtz decomposition of motion
  fields}.
\newblock \bibinfo{journal}{\emph{Pattern Recognition Letters}}
  \bibinfo{volume}{26}, \bibinfo{number}{4} (\bibinfo{year}{2005}),
  \bibinfo{pages}{493--501}.
\newblock


\bibitem[Hodge(1989)]%
        {hodge1989theory}
\bibfield{author}{\bibinfo{person}{William Vallance~Douglas Hodge}.}
  \bibinfo{year}{1989}\natexlab{}.
\newblock \bibinfo{booktitle}{\emph{The theory and applications of harmonic
  integrals}}.
\newblock \bibinfo{publisher}{CUP Archive}.
\newblock


\bibitem[Keros and Subr(2023)]%
        {keros2023spectral}
\bibfield{author}{\bibinfo{person}{Alexandros Keros} {and}
  \bibinfo{person}{Kartic Subr}.} \bibinfo{year}{2023}\natexlab{}.
\newblock \showarticletitle{Spectral coarsening with hodge laplacians}. In
  \bibinfo{booktitle}{\emph{ACM SIGGRAPH 2023 Conference Proceedings}}.
  \bibinfo{pages}{1--11}.
\newblock


\bibitem[Lange et~al\mbox{.}(2022)]%
        {lange2022cellrank}
\bibfield{author}{\bibinfo{person}{Marius Lange}, \bibinfo{person}{Volker
  Bergen}, \bibinfo{person}{Michal Klein}, \bibinfo{person}{Manu Setty},
  \bibinfo{person}{Bernhard Reuter}, \bibinfo{person}{Mostafa Bakhti},
  \bibinfo{person}{Heiko Lickert}, \bibinfo{person}{Meshal Ansari},
  \bibinfo{person}{Janine Schniering}, \bibinfo{person}{Herbert~B Schiller},
  \bibinfo{person}{Dana Pe'er}, {and} \bibinfo{person}{Fabian~J. Theis}.}
  \bibinfo{year}{2022}\natexlab{}.
\newblock \showarticletitle{CellRank for directed single-cell fate mapping}.
\newblock \bibinfo{journal}{\emph{Nat. Methods}}  \bibinfo{volume}{19}
  (\bibinfo{year}{2022}), \bibinfo{pages}{159--170}.
\newblock


\bibitem[Lim(2020)]%
        {lim2020hodge}
\bibfield{author}{\bibinfo{person}{Lek-Heng Lim}.}
  \bibinfo{year}{2020}\natexlab{}.
\newblock \showarticletitle{Hodge Laplacians on graphs}.
\newblock \bibinfo{journal}{\emph{Siam Review}} \bibinfo{volume}{62},
  \bibinfo{number}{3} (\bibinfo{year}{2020}), \bibinfo{pages}{685--715}.
\newblock


\bibitem[Liu et~al\mbox{.}(2015)]%
        {liu2015model}
\bibfield{author}{\bibinfo{person}{Beibei Liu}, \bibinfo{person}{Gemma Mason},
  \bibinfo{person}{Julian Hodgson}, \bibinfo{person}{Yiying Tong}, {and}
  \bibinfo{person}{Mathieu Desbrun}.} \bibinfo{year}{2015}\natexlab{}.
\newblock \showarticletitle{Model-reduced variational fluid simulation}.
\newblock \bibinfo{journal}{\emph{ACM Transactions on Graphics (TOG)}}
  \bibinfo{volume}{34}, \bibinfo{number}{6} (\bibinfo{year}{2015}),
  \bibinfo{pages}{1--12}.
\newblock


\bibitem[Liu et~al\mbox{.}(2023)]%
        {liu2023algebraic}
\bibfield{author}{\bibinfo{person}{Jian Liu}, \bibinfo{person}{Jingyan Li},
  {and} \bibinfo{person}{Jie Wu}.} \bibinfo{year}{2023}\natexlab{}.
\newblock \showarticletitle{The algebraic stability for persistent Laplacians}.
\newblock \bibinfo{journal}{\emph{arXiv preprint arXiv:2302.03902}}
  (\bibinfo{year}{2023}).
\newblock


\bibitem[Liu et~al\mbox{.}(2022)]%
        {liu2022biomolecular}
\bibfield{author}{\bibinfo{person}{Jian Liu}, \bibinfo{person}{Ke-Lin Xia},
  \bibinfo{person}{Jie Wu}, \bibinfo{person}{Stephen Shing-Toung Yau}, {and}
  \bibinfo{person}{Guo-Wei Wei}.} \bibinfo{year}{2022}\natexlab{}.
\newblock \showarticletitle{Biomolecular topology: Modelling and analysis}.
\newblock \bibinfo{journal}{\emph{Acta Mathematica Sinica, English Series}}
  \bibinfo{volume}{38}, \bibinfo{number}{10} (\bibinfo{year}{2022}),
  \bibinfo{pages}{1901--1938}.
\newblock


\bibitem[Mahdessian et~al\mbox{.}(2021)]%
        {mahdessian2021spatiotemporal}
\bibfield{author}{\bibinfo{person}{Diana Mahdessian},
  \bibinfo{person}{Anthony~J Cesnik}, \bibinfo{person}{Christian Gnann},
  \bibinfo{person}{Frida Danielsson}, \bibinfo{person}{Lovisa Stenstr{\"o}m},
  \bibinfo{person}{Muhammad Arif}, \bibinfo{person}{Cheng Zhang},
  \bibinfo{person}{Trang Le}, \bibinfo{person}{Fredric Johansson},
  \bibinfo{person}{Rutger Schutten}, \bibinfo{person}{Anna Backstrom},
  \bibinfo{person}{Ulrika Axelsson}, \bibinfo{person}{Peter Thul},
  \bibinfo{person}{Nathan~H. Cho}, \bibinfo{person}{Oana Carja},
  \bibinfo{person}{Mathias Uhlén}, \bibinfo{person}{Adil Mardinoglu},
  \bibinfo{person}{Charlotte Stadler}, \bibinfo{person}{Cecilia Lindskog},
  \bibinfo{person}{Burcu Ayoglu}, \bibinfo{person}{Manuel~D. Leonetti},
  \bibinfo{person}{Fredrik Ponten}, \bibinfo{person}{Devin~P. Sullivan}, {and}
  \bibinfo{person}{Emma Lundberg}.} \bibinfo{year}{2021}\natexlab{}.
\newblock \showarticletitle{Spatiotemporal dissection of the cell cycle with
  single-cell proteogenomics}.
\newblock \bibinfo{journal}{\emph{Nature}}  \bibinfo{volume}{590}
  (\bibinfo{year}{2021}), \bibinfo{pages}{649--654}.
\newblock


\bibitem[Morrey(1956)]%
        {morrey1956variational}
\bibfield{author}{\bibinfo{person}{Charles~B Morrey}.}
  \bibinfo{year}{1956}\natexlab{}.
\newblock \showarticletitle{A variational method in the theory of harmonic
  integrals, II}.
\newblock \bibinfo{journal}{\emph{American Journal of Mathematics}}
  \bibinfo{volume}{78}, \bibinfo{number}{1} (\bibinfo{year}{1956}),
  \bibinfo{pages}{137--170}.
\newblock


\bibitem[Ng et~al\mbox{.}(2009)]%
        {ng2009efficient}
\bibfield{author}{\bibinfo{person}{Yen~Ting Ng}, \bibinfo{person}{Chohong Min},
  {and} \bibinfo{person}{Fr{\'e}d{\'e}ric Gibou}.}
  \bibinfo{year}{2009}\natexlab{}.
\newblock \showarticletitle{An efficient fluid--solid coupling algorithm for
  single-phase flows}.
\newblock \bibinfo{journal}{\emph{J. Comput. Phys.}} \bibinfo{volume}{228},
  \bibinfo{number}{23} (\bibinfo{year}{2009}), \bibinfo{pages}{8807--8829}.
\newblock


\bibitem[Nguyen et~al\mbox{.}(2020a)]%
        {nguyen2020review}
\bibfield{author}{\bibinfo{person}{Duc~Duy Nguyen}, \bibinfo{person}{Zixuan
  Cang}, {and} \bibinfo{person}{Guo-Wei Wei}.}
  \bibinfo{year}{2020}\natexlab{a}.
\newblock \showarticletitle{A review of mathematical representations of
  biomolecular data}.
\newblock \bibinfo{journal}{\emph{Physical Chemistry Chemical Physics}}
  \bibinfo{volume}{22}, \bibinfo{number}{8} (\bibinfo{year}{2020}),
  \bibinfo{pages}{4343--4367}.
\newblock


\bibitem[Nguyen et~al\mbox{.}(2019)]%
        {nguyen2019mathematical}
\bibfield{author}{\bibinfo{person}{Duc~Duy Nguyen}, \bibinfo{person}{Zixuan
  Cang}, \bibinfo{person}{Kedi Wu}, \bibinfo{person}{Menglun Wang},
  \bibinfo{person}{Yin Cao}, {and} \bibinfo{person}{Guo-Wei Wei}.}
  \bibinfo{year}{2019}\natexlab{}.
\newblock \showarticletitle{Mathematical deep learning for pose and binding
  affinity prediction and ranking in D3R Grand Challenges}.
\newblock \bibinfo{journal}{\emph{Journal of computer-aided molecular design}}
  \bibinfo{volume}{33} (\bibinfo{year}{2019}), \bibinfo{pages}{71--82}.
\newblock


\bibitem[Nguyen et~al\mbox{.}(2020b)]%
        {nguyen2020mathdl}
\bibfield{author}{\bibinfo{person}{Duc~Duy Nguyen}, \bibinfo{person}{Kaifu
  Gao}, \bibinfo{person}{Menglun Wang}, {and} \bibinfo{person}{Guo-Wei Wei}.}
  \bibinfo{year}{2020}\natexlab{b}.
\newblock \showarticletitle{MathDL: mathematical deep learning for D3R Grand
  Challenge 4}.
\newblock \bibinfo{journal}{\emph{Journal of computer-aided molecular design}}
  \bibinfo{volume}{34} (\bibinfo{year}{2020}), \bibinfo{pages}{131--147}.
\newblock


\bibitem[Petronetto et~al\mbox{.}(2009)]%
        {petronetto2009meshless}
\bibfield{author}{\bibinfo{person}{Fabiano Petronetto}, \bibinfo{person}{Afonso
  Paiva}, \bibinfo{person}{Marcos Lage}, \bibinfo{person}{Geovan Tavares},
  \bibinfo{person}{H{\'e}lio Lopes}, {and} \bibinfo{person}{Thomas Lewiner}.}
  \bibinfo{year}{2009}\natexlab{}.
\newblock \showarticletitle{Meshless helmholtz-hodge decomposition}.
\newblock \bibinfo{journal}{\emph{IEEE transactions on visualization and
  computer graphics}} \bibinfo{volume}{16}, \bibinfo{number}{2}
  (\bibinfo{year}{2009}), \bibinfo{pages}{338--349}.
\newblock


\bibitem[Poelke(2017)]%
        {poelke2017hodge}
\bibfield{author}{\bibinfo{person}{Konstantin Poelke}.}
  \bibinfo{year}{2017}\natexlab{}.
\newblock \emph{\bibinfo{title}{Hodge-type decompositions for piecewise
  constant vector fields on simplicial surfaces and solids with boundary}}.
\newblock \bibinfo{thesistype}{Ph.\,D. Dissertation}.
\newblock


\bibitem[Poelke and Polthier(2016)]%
        {poelke2016boundary}
\bibfield{author}{\bibinfo{person}{Konstantin Poelke} {and}
  \bibinfo{person}{Konrad Polthier}.} \bibinfo{year}{2016}\natexlab{}.
\newblock \showarticletitle{Boundary-aware Hodge decompositions for piecewise
  constant vector fields}.
\newblock \bibinfo{journal}{\emph{Computer-Aided Design}}  \bibinfo{volume}{78}
  (\bibinfo{year}{2016}), \bibinfo{pages}{126--136}.
\newblock


\bibitem[Polthier and Preu{\ss}(2000)]%
        {polthier2000variational}
\bibfield{author}{\bibinfo{person}{Konrad Polthier} {and} \bibinfo{person}{Eike
  Preu{\ss}}.} \bibinfo{year}{2000}\natexlab{}.
\newblock \showarticletitle{Variational approach to vector field
  decomposition}. In \bibinfo{booktitle}{\emph{Data Visualization 2000:
  Proceedings of the Joint EUROGRAPHICS and IEEE TCVG Symposium on
  Visualization in Amsterdam, The Netherlands, May 29--30, 2000}}. Springer,
  \bibinfo{pages}{147--155}.
\newblock


\bibitem[Polthier and Preu{\ss}(2003)]%
        {polthier2003identifying}
\bibfield{author}{\bibinfo{person}{Konrad Polthier} {and} \bibinfo{person}{Eike
  Preu{\ss}}.} \bibinfo{year}{2003}\natexlab{}.
\newblock \showarticletitle{Identifying vector field singularities using a
  discrete Hodge decomposition}. In \bibinfo{booktitle}{\emph{Visualization and
  mathematics III}}. Springer, \bibinfo{pages}{113--134}.
\newblock


\bibitem[Razafindrazaka et~al\mbox{.}(2019)]%
        {razafindrazaka2019consistent}
\bibfield{author}{\bibinfo{person}{Faniry~H Razafindrazaka},
  \bibinfo{person}{Konstantin Poelke}, \bibinfo{person}{Konrad Polthier}, {and}
  \bibinfo{person}{Leonid Goubergrits}.} \bibinfo{year}{2019}\natexlab{}.
\newblock \showarticletitle{A Consistent Discrete 3D Hodge-type Decomposition:
  implementation and practical evaluation}.
\newblock \bibinfo{journal}{\emph{arXiv preprint arXiv:1911.12173}}
  (\bibinfo{year}{2019}).
\newblock


\bibitem[Ribando-Gros et~al\mbox{.}(2024)]%
        {ribando2022graph}
\bibfield{author}{\bibinfo{person}{Emily Ribando-Gros}, \bibinfo{person}{Rui
  Wang}, \bibinfo{person}{Jiahui Chen}, \bibinfo{person}{Yiying Tong}, {and}
  \bibinfo{person}{Guo-Wei Wei}.} \bibinfo{year}{2024}\natexlab{}.
\newblock \showarticletitle{Combinatorial and Hodge Laplacians: Similarity and
  difference}.
\newblock \bibinfo{journal}{\emph{SIAM Rev.}} \bibinfo{volume}{66},
  \bibinfo{number}{3} (\bibinfo{year}{2024}), \bibinfo{pages}{575--601}.
\newblock


\bibitem[Sawhney and Crane(2020)]%
        {sawhney2020monte}
\bibfield{author}{\bibinfo{person}{Rohan Sawhney} {and} \bibinfo{person}{Keenan
  Crane}.} \bibinfo{year}{2020}\natexlab{}.
\newblock \showarticletitle{Monte Carlo geometry processing: A grid-free
  approach to PDE-based methods on volumetric domains}.
\newblock \bibinfo{journal}{\emph{ACM Transactions on Graphics}}
  \bibinfo{volume}{39}, \bibinfo{number}{4} (\bibinfo{year}{2020}).
\newblock


\bibitem[Schwarz(2006)]%
        {schwarz2006hodge}
\bibfield{author}{\bibinfo{person}{G. Schwarz}.}
  \bibinfo{year}{2006}\natexlab{}.
\newblock \bibinfo{booktitle}{\emph{Hodge Decomposition - A Method for Solving
  Boundary Value Problems}}.
\newblock \bibinfo{publisher}{Springer Berlin Heidelberg}.
\newblock
\showISBNx{9783540494034}
\urldef\tempurl%
\url{https://books.google.com/books?id=6-17CwAAQBAJ}
\showURL{%
\tempurl}


\bibitem[Shonkwiler(2009)]%
        {shonkwiler2009poincare}
\bibfield{author}{\bibinfo{person}{Clayton Shonkwiler}.}
  \bibinfo{year}{2009}\natexlab{}.
\newblock \emph{\bibinfo{title}{Poincar{\'e} duality angles on Riemannian
  manifolds with boundary}}.
\newblock \bibinfo{thesistype}{Ph.\,D. Dissertation}.
  \bibinfo{school}{University of Pennsylvania}.
\newblock


\bibitem[Shonkwiler(2013)]%
        {shonkwiler2013poincare}
\bibfield{author}{\bibinfo{person}{Clayton Shonkwiler}.}
  \bibinfo{year}{2013}\natexlab{}.
\newblock \showarticletitle{Poincar{\'e} duality angles and the
  Dirichlet-to-Neumann operator}.
\newblock \bibinfo{journal}{\emph{Inverse Problems}} \bibinfo{volume}{29},
  \bibinfo{number}{4} (\bibinfo{year}{2013}), \bibinfo{pages}{045007}.
\newblock


\bibitem[Stam(2023)]%
        {stam2023stable}
\bibfield{author}{\bibinfo{person}{Jos Stam}.} \bibinfo{year}{2023}\natexlab{}.
\newblock \showarticletitle{Stable fluids}.
\newblock In \bibinfo{booktitle}{\emph{Seminal Graphics Papers: Pushing the
  Boundaries, Volume 2}}. \bibinfo{pages}{779--786}.
\newblock


\bibitem[Su et~al\mbox{.}(2024a)]%
        {su2024hodge}
\bibfield{author}{\bibinfo{person}{Zhe Su}, \bibinfo{person}{Yiying Tong},
  {and} \bibinfo{person}{Guo-Wei Wei}.} \bibinfo{year}{2024}\natexlab{a}.
\newblock \showarticletitle{Hodge decomposition of single-cell RNA velocity}.
\newblock \bibinfo{journal}{\emph{Journal of chemical information and
  modeling}} (\bibinfo{year}{2024}).
\newblock


\bibitem[Su et~al\mbox{.}(2024b)]%
        {su2024persistent}
\bibfield{author}{\bibinfo{person}{Zhe Su}, \bibinfo{person}{Yiying Tong},
  {and} \bibinfo{person}{Guo-Wei Wei}.} \bibinfo{year}{2024}\natexlab{b}.
\newblock \showarticletitle{Persistent de Rham-Hodge Laplacians in the Eulerian
  representation}.
\newblock \bibinfo{journal}{\emph{arXiv preprint arXiv:2408.00220}}
  (\bibinfo{year}{2024}).
\newblock


\bibitem[Tong et~al\mbox{.}(2003)]%
        {tong2003discrete}
\bibfield{author}{\bibinfo{person}{Yiying Tong}, \bibinfo{person}{Santiago
  Lombeyda}, \bibinfo{person}{Anil~N Hirani}, {and} \bibinfo{person}{Mathieu
  Desbrun}.} \bibinfo{year}{2003}\natexlab{}.
\newblock \showarticletitle{Discrete multiscale vector field decomposition}.
\newblock \bibinfo{journal}{\emph{ACM transactions on graphics (TOG)}}
  \bibinfo{volume}{22}, \bibinfo{number}{3} (\bibinfo{year}{2003}),
  \bibinfo{pages}{445--452}.
\newblock


\bibitem[Wang et~al\mbox{.}(2020)]%
        {wang2020persistent}
\bibfield{author}{\bibinfo{person}{Rui Wang}, \bibinfo{person}{Duc~Duy Nguyen},
  {and} \bibinfo{person}{Guo-Wei Wei}.} \bibinfo{year}{2020}\natexlab{}.
\newblock \showarticletitle{Persistent spectral graph}.
\newblock \bibinfo{journal}{\emph{International journal for numerical methods
  in biomedical engineering}} \bibinfo{volume}{36}, \bibinfo{number}{9}
  (\bibinfo{year}{2020}), \bibinfo{pages}{e3376}.
\newblock


\bibitem[Wang and Chern(2021)]%
        {wang2021computing}
\bibfield{author}{\bibinfo{person}{Stephanie Wang} {and}
  \bibinfo{person}{Albert Chern}.} \bibinfo{year}{2021}\natexlab{}.
\newblock \showarticletitle{Computing minimal surfaces with differential
  forms}.
\newblock \bibinfo{journal}{\emph{ACM Transactions on Graphics (TOG)}}
  \bibinfo{volume}{40}, \bibinfo{number}{4} (\bibinfo{year}{2021}),
  \bibinfo{pages}{1--14}.
\newblock


\bibitem[Yang et~al\mbox{.}(2021)]%
        {yang2021clebsch}
\bibfield{author}{\bibinfo{person}{Shuqi Yang}, \bibinfo{person}{Shiying
  Xiong}, \bibinfo{person}{Yaorui Zhang}, \bibinfo{person}{Fan Feng},
  \bibinfo{person}{Jinyuan Liu}, {and} \bibinfo{person}{Bo Zhu}.}
  \bibinfo{year}{2021}\natexlab{}.
\newblock \showarticletitle{Clebsch gauge fluid}.
\newblock \bibinfo{journal}{\emph{ACM Transactions on Graphics (TOG)}}
  \bibinfo{volume}{40}, \bibinfo{number}{4} (\bibinfo{year}{2021}),
  \bibinfo{pages}{1--11}.
\newblock


\bibitem[Zhao et~al\mbox{.}(2019)]%
        {zhao20193d}
\bibfield{author}{\bibinfo{person}{Rundong Zhao}, \bibinfo{person}{Mathieu
  Desbrun}, \bibinfo{person}{Guo-Wei Wei}, {and} \bibinfo{person}{Yiying
  Tong}.} \bibinfo{year}{2019}\natexlab{}.
\newblock \showarticletitle{3D Hodge decompositions of edge-and face-based
  vector fields}.
\newblock \bibinfo{journal}{\emph{ACM Transactions on Graphics (TOG)}}
  \bibinfo{volume}{38}, \bibinfo{number}{6} (\bibinfo{year}{2019}),
  \bibinfo{pages}{1--13}.
\newblock


\bibitem[Zomorodian and Carlsson(2004)]%
        {zomorodian2004computing}
\bibfield{author}{\bibinfo{person}{Afra Zomorodian} {and}
  \bibinfo{person}{Gunnar Carlsson}.} \bibinfo{year}{2004}\natexlab{}.
\newblock \showarticletitle{Computing persistent homology}. In
  \bibinfo{booktitle}{\emph{Proceedings of the twentieth annual symposium on
  Computational geometry}}. \bibinfo{pages}{347--356}.
\newblock


\end{thebibliography}
\newpage


\end{document}